\documentclass[11pt,english,reqno]{amsart}
\RequirePackage[utf8]{inputenc}
\usepackage{tipa}
\usepackage{amstext,amsopn,amscd,amsthm,tikz}
\usetikzlibrary{decorations.pathreplacing}
\usepackage{amsfonts,amsmath,amssymb,amsbsy}
\usepackage{layout}
\usepackage{newlfont}
\usepackage{rotating}
\usepackage{verbatim}
\usepackage{lipsum}
\usepackage{buczokoch}
\usepackage{version} 
\usepackage{bm}

\RequirePackage{graphicx}
\graphicspath{{Pictures/}}
\usepackage{float}

\usepackage{caption}
\usepackage{subcaption}

\usepackage{wrapfig}
\usepackage{xcolor}
\usepackage{lscape}

\usepackage{array}
\usepackage{diagbox}
\usepackage{colortbl}
\usepackage{hhline}
\usepackage{rotating,multirow,multicol}

\usepackage[bookmarks=false,colorlinks=true,linkcolor=magenta,urlcolor=magenta,citecolor=magenta,filecolor=magenta,linktocpage=true,pdfstartview={XYZ null null 0.62}]{hyperref}

\usepackage{titletoc}
\DeclareRobustCommand{\SkipTocEntry}[5]{}

\newcommand{\N}{\mathbb N}

\newcommand{\R}{\mathbb R}

\newcommand\dos{d_0^*}
\newcommand\alphat{{\widetilde{\alpha}}}
\newcommand{\la}{\lambda}
\newcommand\ga{{\gamma_\lambda}}
\newcommand\fla{{F^\lambda}}
\newcommand\flan{{F_n^{\lambda}}}
\newcommand\mula{{\mu_\la}}

\newcommand{\cB}{\cal{B}}
\newcommand{\cC}{\cal{C}}

\newcommand{\cE}{\cal{E}}

\newcommand{\cH}{\cal{H}}
\newcommand{\cI}{\cal{I}}

\newcommand{\cK}{\cal{K}}
\newcommand{\cL}{\cal{L}}

\newcommand{\cR}{\cal{R}}

\newcommand{\cV}{\cal{V}}

\newcommand{\eps}{\varepsilon}
\newcommand{\ep}{\eps}

\newcommand{\sgn}{\mathrm{sgn}}
\newcommand{\dimh}{\mathrm{dim_H}}

\newcommand\cEt{\widetilde{\cE}}

\usepackage{enumitem}
\usepackage{nccmath}
\numberwithin{equation}{section}

\newtheorem{theorem}{Theorem}[section]
\newtheorem{corollary}[theorem]{Corollary}

\newtheorem{lemma}[theorem]{Lemma}
\newtheorem{proposition}[theorem]{Proposition}

\theoremstyle{definition}
\newtheorem{definition}[theorem]{Definition}
\newtheorem{remark}[theorem]{Remark}

\title[The parametrized von Koch functions]{The multifractal nature of a parametrized family\\ of von Koch functions}

\author[Z. Buczolich, Y. Demichel and S. Seuret]{Zolt\'an Buczolich\textsuperscript{1}, Yann Demichel\textsuperscript{2} and St\'{e}phane Seuret\textsuperscript{3}}

\address{\textsuperscript{1} Department of Analysis, ELTE E\"otv\"os Lor\'and\\
University, P\'azm\'any P\'eter S\'et\'any 1/c, 1117 Budapest, Hungary}
\email{zoltan.buczolich@ttk.elte.hu}\urladdr{http://buczo.web.elte.hu, ORCID Id: 0000-0001-5481-8797}%

\address{\textsuperscript{2} Laboratoire MODAL'X, UMR CNRS 9023, UPL, Universit\'e Paris Nanterre, 200 avenue de la R\'e\-pu\-bli\-que, 92001 Nanterre, France.}
\email{ydemichel@parisnanterre.fr}\urladdr{https://www.parisnanterre.fr/m-yann-demichel, ORCID Id: 0009-0008-4306-9469}%

\address{\textsuperscript{3} St\'ephane Seuret, Univ Paris Est Creteil, Univ Gustave Eiffel, CNRS, LAMA UMR8050, F-94010 Creteil, France}
\email{seuret@u-pec.fr}
\urladdr{https://sites.math.u-pem.fr/sseuret/, ORCID Id: 0000-0002-7113-3156}%

\begin{document}

\begin{abstract} In a famous paper published in 1904, Helge von Koch introduced the curve that still serves nowadays as an iconic representation of fractal shapes. In fact, von Koch's main goal was the construction of a continuous but nowhere differentiable function, very similar to the snowflake, using elementary geometric procedures, and not analytical formulae. We prove that a parametrized family of functions (including and) generalizing von Koch's example enjoys a rich multifractal behavior, thus enriching the class of historical mathematical objects having surprising regularity properties. The analysis relies on the study of the orbits of an underlying dynamical system and on the introduction of self-similar measures and non-trivial iterated functions systems adapted to the problem.
\end{abstract}


\vspace*{-0.05cm}
\maketitle

\tableofcontents

\section{Introduction}\label{sec:intro}

In 1872, Weierstrass presented to the Royal Prussian Academy of Sciences the first explicit example of a family of real functions that are continuous but nowhere differentiable, see \cite{Weier72}, definitively refuting the common belief that a continuous curve must admit well-defined tangents at each of its points. From this date, several examples and parametrized families of such functions modeled by lacunary {Fourier-like} series have been successively proposed by Darboux \cite{darb75}, Dini \cite{dini77}, Cell\'{e}rier \cite{cell90}, or Takagi \cite{taka03}. Let us also mention the famous Riemann-Fourier series $\sum_{n\geq 1}\frac{\sin(n^2\pi x)}{n^2}$ which continues to attract the attention of mathematicians from different fields, but which is known to be differentiable at some rational points; see \cite{ger69}. Very disappointed by the fact that the construction of all these functions, as well as the proofs of their pathological irregularity, were based on analytical formulas rather than geometric and intuitive arguments, the Swedish mathematician Helge von Koch looked for different examples. In 1904 he succeeded in constructing a continuous but nowhere differentiable Jordan curve using a basic geometric method; see \cite{Koch04,Koch06}. The curve he exhibited gave birth to the so-called {\it snowflake curve} and became one of the most iconic curves of the 20th century. In its original form, this curve does not represent a function, but with a slight modification of the construction, von Koch finally obtained a function with the desired properties.

\medskip

\begin{figure}[h!]
\includegraphics[width=0.475\textwidth]{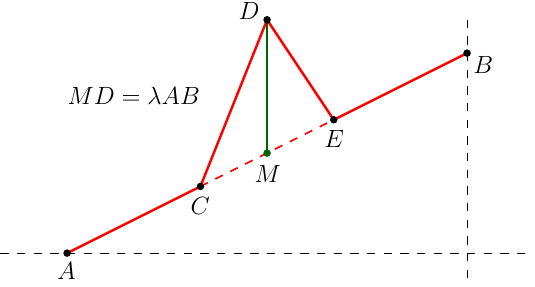}
\captionsetup{width=1.0\textwidth}
\caption{The basic geometric operation $\Omega^\la$ transforms a given line segment $AB$ into the polygonal line $ACDEB$ made up of $4$ consecutive line segments such that $|AC|=|CE|=|EB|$ and $|MD|=\la |AB|$.}
\label{fig:OmegaLambda}
\end{figure}

\noindent{\bf The generalized von Koch function $\fla$.} Roughly speaking, the graph of the von Koch function is obtained in the same way as the snowflake curve, but the triangles added at each stage of the procedure are straightened up vertically. Precisely, fix a parameter $\la>0$ and consider in the plane a line segment $AB$ that forms an arbitrary angle with the $x$-axis. Divide $AB$ into three equal parts $AC$, $CE$, $EB$, and erect a triangle $CDE$ with base $CE$ and median $MD$ parallel to the $y$-axis, directed towards the positive $y$ direction and with length $|MD|=\la |AB|$. Following the footsteps of von Koch, we call $\Omega^\la$ the operation that transforms the line segment $AB$ into the polygonal line $ACDEB$, see Figure \ref{fig:OmegaLambda}.

Call $F^\la_0$ the zero function on $[0,1]$ and $AB$ the segment that represents it. Then define $F^\la_1$ as the piecewise affine function whose graph is the polygonal line obtained by applying $\Omega^\la$ to $AB$, i.e.
\begin{equation*}
F^\la_1(x) =
\begin{cases}
 0 & \mbox{ if } x\in[0,\frac13], \\
 6\la(x-\frac13) & \mbox{ if } x\in[\frac13,\frac12], \\
 6\la(\frac23-x) & \mbox{ if } x\in[\frac12,\frac23], \\
 0 & \mbox{ if } x\in[\frac23,1].
\end{cases}
\end{equation*}

\begin{figure}[h!]
\includegraphics[width=0.495\textwidth]{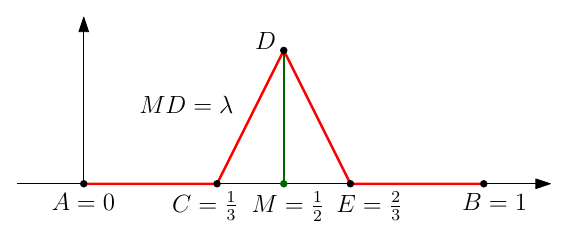}
\includegraphics[width=0.495\textwidth]{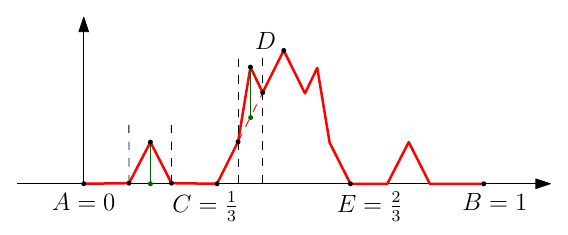}
\captionsetup{width=1.0\textwidth}
\caption{Iterative construction of the generalized von Koch function $\fla$. Left: Applying $\Omega^\la$ to the line segment $AB$ gives the graph of $F^\la_1$. Right: Applying $\Omega^\la$ to each line segment $AC$, $CD$, $DE$ and $EB$ gives the graph of $F^\la_2$.}
\label{fig:FunctionStep}
\end{figure}

The function $F^\la_1$ is then affine over the four maximal intervals $[0,\frac13]$, $[\frac13,\frac12]$, $[\frac12,\frac23]$ and $[\frac23,1]$. Applying now $\Omega^\la$ to each line segment corresponding to each of these intervals we obtain the graph of the function $F^\la_2$, see Figure \ref{fig:FunctionStep}.

\begin{figure}
\centering
\includegraphics[scale=0.55]{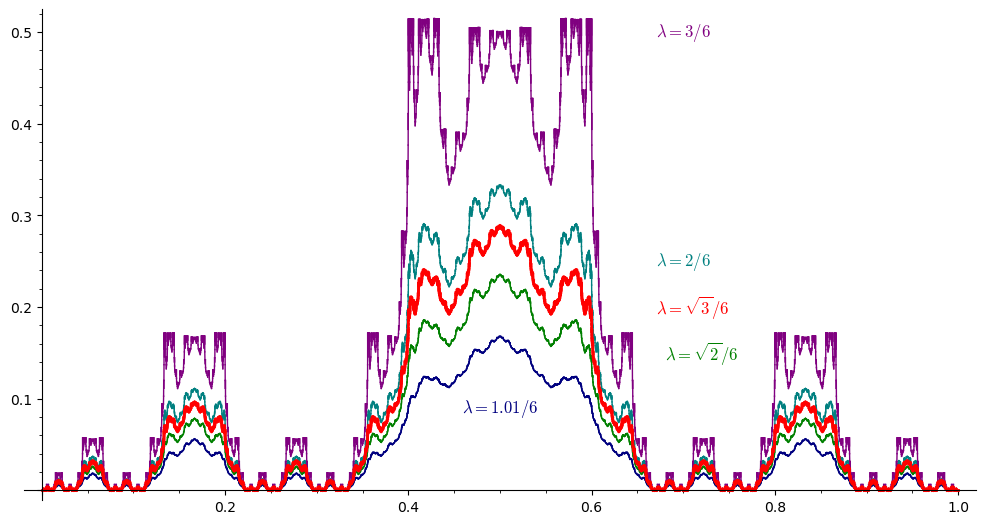}
\captionsetup{width=1\textwidth}
\caption{The generalized von Koch function $\fla$ for various parameters $\la$. The original function by von Koch corresponding to the choice $\la=\frac{\sqrt3}6$ is shown by the thick red curve.}
\label{fig:HVKsamples}
\end{figure}

The construction can be continued infinitely and produces a sequence of piecewise affine functions $\flan$, $n\geq0$ (note that $\flan$ is affine over exactly $4^n$ intervals). Since the sequence $(\flan)_{n\geq0}$ is non-decreasing, its limit function, denoted by $\fla$ and named the {\it generalized von Koch function with parameter $\la$}, exists (with eventually some infinite values); see Figure \ref{fig:HVKsamples} for some examples. The function originally considered by von Koch corresponds to the choice $\la=\frac{\sqrt3}6$, the thick red curve in the figure. In this case, $\Omega^{\frac{\sqrt3}6}$ transforms horizontal segments by replacing the median interval by the two sides of an equilateral triangle as for the classical snowflake curve.

In \cite{Koch04,Koch06}, von Koch proved by elementary but delicate arguments that $F^{\frac{\sqrt3}6}$ is nowhere differentiable on $[0,1]$. In fact, he only stated that this function has nowhere a finite derivative, leaving open the question of whether there might be points at which the derivative is infinite. He also did not go further in studying the pointwise regularity of $F^{\frac{\sqrt3}6}$, and it seems that this question has never been asked. Even with the rise of multifractal analysis, while most famous pathological functions constructed in the last century have been rediscovered and their multifractal spectrum determined, see e.g. \cite{Jaff97,JM}, von Koch's example has not been considered. Curiously, von Koch's function seems to have been forgotten by mathematicians, while the snowflake has become one of the most iconic mathematical objects. Perhaps it was thought to have the same self-similarity property and monofractal behavior as the snowflake.

\medskip

\noindent{\bf Main results.} The present paper offers a thorough analysis of the local and multiscale behavior of the generalized von Koch function $\fla$ when $\la\in (\frac{\sqrt2}{6},\frac56)$, especially its multifractal properties. In general, the multifractal analysis of a mathematical object $X$ consists in describing its local regularity properties together with its multiscale properties, and connecting them with each other. On one hand, one computes the multifractal spectrum of $X$, i.e. the mapping $d_X$ which associates with every $h\geq 0$ the Hausdorff dimension of the set of points having a pointwise regularity equal to $h$. On the other hand, one describes the global statistical properties of $X$ via the study of its $L^q$-scaling function $\tau_X$ (see Definitions \ref{defspectrum} and Section \ref{sec:mula} for more details on the definitions of pointwise regularity for functions and measures, respectively). Multifractals originate from the study of turbulence \cite{FrischParisi} and thermodynamics \cite{Halsey}, and have proved insightful in many mathematical domains, from real and harmonic analysis \cite{JAFF_FRISCH,BS-FP-1} to probability theory \cite{JaffardLevy,Vik-Law} and dynamical systems \cite{Pesin2,Shmerkin1}, with deep connection to geometric measure theory and metric number theory \cite{JM}. Recently, multifractality has been proven for solutions to PDEs connected to the binormal flow \cite{BanicaVega}, and investigating the multifractal nature of mathematical objects is a very active research area.

We perform the multifractal analysis of $\fla$ in a certain range of $\la$, proving that $\fla$ enjoys an extremely rich multifractal behavior, which probably von Koch did not suspect. Surprisingly, the multifractal study of $\fla$ is deeply dependent on different subranges of $\la$. A key fact, as well as a challenging issue, is that, unlike classical functions such as the Weierstrass series, there is no closed analytic representation for $\fla$. In the same vein, despite the resemblance of its graph to the so-called ``alternating Takagi function'', see e.g. \cite{All11}, the von Koch function is not a self-similar de Rham curve, and cannot be directly studied by using functional equations or as the attractor of an Iterated Function System (IFS). The analysis of $\fla$ requires the study of a dynamical system on $[0,1]$, the multifractal properties of auxiliary self-similar probability measures, and the construction of several IFSs tailored to the problem.

In this paper, we focus on the range of parameters $(\frac16,\frac56)$. In terms of the geometric construction of $\flan$, $\frac16$ corresponds to the critical value under which it may occur that after iterating the construction scheme $\Omega^\la$ on a segment with positive (resp. negative) slope, the four new segments all have a positive (resp. negative) slope (see Figure \ref{fig:OmegaLambdaBis} below). In the range $(\frac16,\frac56)$, such a phenomenon does not occur, and the functions $\fla$ follow a construction scheme similar to that of von Koch's original example. Also, due to this geometrical phenomenon, the strategy we develop, based on the study of the dynamical system $([0,1],\cB([0,1]),T)$ introduced below, is not adapted and the whole problem should be approached differently.

\begin{figure}[ht!]
\centering
\includegraphics[width=0.35\textwidth]{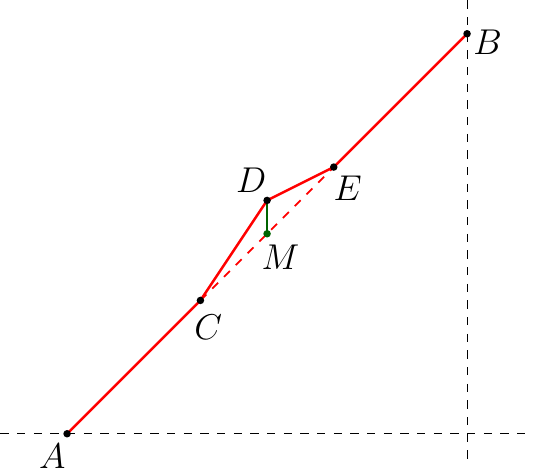}
\captionsetup{width=1\textwidth}
\caption{When $\la\in(0,\frac16)$, the operation $\Omega^\la$ may lead to new slopes all with the same sign, hence a monotonous part in the graph of $\flan$.}
\label{fig:OmegaLambdaBis}
\end{figure}

Observe that if $\la\leq \la'$ then $F^\la(x)\leq F^{\la'}(x)$ for every $x\in[0,1]$. Our first result gives the critical value $\la$ such that $F^\la(x)<+\infty$ for every $x\in[0,1]$.

\begin{theorem}\label{theo-finite}
When $\la\in [0,\frac56)$, the sequence $(\flan)_{n\geq1}$ converges uniformly to $\fla$, which is continuous but nowhere differentiable. When $\la \geq \frac56$, there are points $x\in[0,1]$ at which $\fla(x)=\lim_{n\to\infty} \flan(x)=+\infty$.
\end{theorem}

\begin{remark}\label{rk1}
Observe that $\fla(x)\geq0$ for every $x\in[0,1]$ and $\fla(x)=0$ if and only if $x$ belongs to the classical triadic Cantor set.
\end{remark}

It follows from its construction that $\fla$ satisfies some de Rham's type equations:
\begin{equation}\label{equafunc}
\forall\, x\in[0,1],\quad
\begin{cases}
\fla(1-x) & =\fla(x), \\
\fla(\frac{x}3) & =\frac13\fla(x).
\end{cases}
\end{equation}
However, in contrast to the snowflake, the graph of $\fla$ is not a self-similar set. This follows for instance from Theorem 1.1 of \cite{Moreira-Xi-Zhang}. We will show in an upcoming work that the graph of $\fla$ is not a self-affine set either; hence it is not possible to derive the pointwise Hölder exponent using this approach (see \cite{Allaart2020} for the multifractal analysis of one-dimensional self-affine functions).

\medskip

In \cite{Koch04,Koch06}, von Koch left open the question of the existence of infinite derivatives of its function at some points $x$. Since $\frac{\sqrt3}6\in (\frac16,\frac12)$, we answer the question positively with the following Theorem.

\begin{theorem}\label{theo-derivinfini}
Let $\cV^\la_{+\infty}$ (resp. $\cV^\la_{-\infty}$) the set of those points $x\in[0,1]$ where the derivative of $\fla$ exists and is equal to $+\infty$ (resp. $-\infty$).
\begin{itemize}[parsep=-0.15cm,itemsep=0.25cm,topsep=0.2cm,wide=0.175cm,leftmargin=0.5cm]
\item[$(i)$] If $\la\in (\frac16,\frac12)$ then $\dimh(\cV^\la_{+\infty})=\dimh(\cV^\la_{-\infty})>0$.
\item[$(ii)$] If $\la\in (\frac23,\frac56)$ then $\cV^\la_{+\infty}=\cV^\la_{-\infty}=\emptyset$.
\end{itemize}
\end{theorem}

This theorem highlights a clear dichotomy between small and large values of the parameter $\lambda$. The transition between these two regimes appears to occur within the interval $[\frac12,\frac23]$. However, understanding the sets $\mathcal V^\lambda_{+\infty}$ and $\mathcal V^\lambda_{-\infty}$ in this range seems rather delicate, and lies beyond the scope of the present work. We suspect that the structure of the sets of points where $\fla$ has infinite derivatives depends intricately on the value of $\la$, in a manner comparable to that observed for the Okamoto self-affine functions (see \cite{Okamoto,Allaart2016}).

\medskip

Let us now turn to the pointwise regularity of $\fla$. Intuitively, the larger the parameter $\la$, the more irregular $\fla$ should be, see Figure \ref{fig:HVKsamples} and statements of Theorem \ref{theo-derivinfini}.

Recall how the pointwise behavior of a locally bounded function is quantified.

\begin{definition}\label{defholderpoint}
Let $f\in L^\infty([0,1])$ and $x_0\in[0,1]$. For $\alpha\geq 0$, $f$ is said to belong to $\cC^\alpha(x_0)$ if there are a polynomial $P$ of degree less than $\lfloor\alpha\rfloor$ and two constants $C,\delta>0$ such that
\begin{equation}\label{defpoint}
\forall\,x\in[0,1],\quad  |x-x_0|<\delta \Longrightarrow  |f(x)-P(x-x_0)| \leq C |x-x_0|^\alpha.
\end{equation}
The pointwise H\"older exponent of $f$ at $x_0$ is
\begin{equation}\label{defholder}
h_f(x_0) = \sup\{\alpha\geq 0: f\in \cC^\alpha(x_0)\}.
\end{equation}
\end{definition}

\medskip

More specifically, the quantity of interest is the multifractal spectrum of $f$.
\begin{definition}\label{defspectrum}
Let $f\in L^\infty([0,1])$. Consider the level sets $E_f(\alpha)$ for $h_f$:
\begin{equation}\label{deflevelset}
E_f(\alpha)=\{x\in[0,1] : h_f(x)=\alpha\}.
\end{equation}
The multifractal spectrum of $f$ is the mapping $d_f:\R^+\to[0,1]\cup\{-\infty\}$ defined by
\begin{equation}\label{defspec}
d_f(\alpha)=\dimh(E_f(\alpha))
\end{equation}
where $\dimh$ stands for the Hausdorff dimension, with the convention $\dimh(\emptyset)=-\infty$. The function $f$ is said to be multifractal if the support $\{\alpha: d_f(\alpha)\not = -\infty \}$ of $d_f$  is larger than a single point.
\end{definition}

In Theorem \ref{spectruma} we determine the multifractal spectrum $d_\fla(\alpha)$ of $\fla$ when $\la\in(\frac{\sqrt 2}6,\frac56)$. This parameter range includes the classical von Koch function corresponding to $\la=\frac{\sqrt3}6$. To state this result, we need to introduce a dynamical system that plays a key role in the analysis of $\fla$, as well as in the definition of the self-similar measure $\mula$ that is intrinsically related to $\fla$ (these relationships will be explored in the following sections).

\medskip

\begin{figure}[!ht]
\centering
\includegraphics[width=0.475\textwidth]{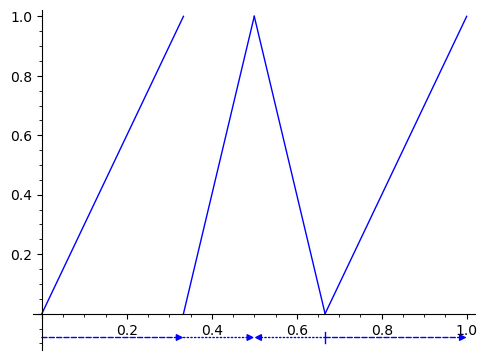}
\captionsetup{width=1\textwidth}
\caption{Illustration for the dynamics of $T$ associated with the construction of the function $\fla$.}
\label{fig:MapT}
\end{figure}

Consider the interval $[0,1]$, $\cB([0,1])$ the Borel $\sigma$-algebra over $[0,1]$, and the transformation $T:[0,1]\to [0,1]$ defined as the piecewise affine map defined by \eqref{def-T} and displayed in Figure \ref{fig:MapT}:

\begin{equation}\label{def-T}
T(x) =
\begin{cases}
 \ \ \  3x & \mbox{ if }  0  \leq x < \frac13, \\
6x -2 & \mbox{ if } \frac13 \leq x < \frac12, \\
4-6x & \mbox{ if } \frac12 \leq x < \frac23, \\
3x-2 & \mbox{ if } \frac23 \leq x \leq 1.
\end{cases}
\end{equation}

Given $x\in[0,1]$, the pointwise regularity of $\fla$ at $x$ will be determined by the orbit of $x$ under the action of $T$, see Sections \ref{sec:dynam} and \ref{sec:localholder}.

\medskip

Next, we build the IFS corresponding to the inverse branches of $T$. Consider the four mappings $S_i:[0,1]\to [0,1]$, $i=0,\ldots,3$, defined by
\begin{equation}\label{defSi}
S_0(x)=\frac{x}{3},\, S_1(x)=\frac{x}{6}+\frac13,\, S_2(x)=-\frac{x}{6}+\frac23 \text{ and } S_3(x)=\frac{x}{3}+\frac23.
\end{equation}
These mappings are contracting similarities with ratios $r_i\in(0,1)$ respectively given by
\begin{equation}\label{defri}
r_0=r_3=\frac13,\, r_1=r_2=\frac16.
\end{equation}
There exists a unique non-empty compact set $K$ satisfying $K=\bigcup_{i=0}^3 S_i(K)$. This set, called the attractor of the IFS, is clearly the interval $[0,1]$. Moreover, one easily checks that this IFS satisfies the Open Set Condition (OSC, see \cite{BRMICHPEY,LauNgai99} and \cite{Hochman14,BaFe21,Fen23} for recent developments), i.e. $S_i((0,1))\cap S_j((0,1))=\emptyset$ as soon as $i\neq j$.

\medskip

Let us now introduce the particular self-similar measure $\mula$ associated with this IFS and supported on its attractor $[0,1]$.

\begin{definition}\label{def:mula}
Let $\la\in(\frac16,\frac56)$ and let $\ga \geq 1$ be the unique real number such that
\begin{equation}\label{def-gamma}
\frac1{3^{\ga}} +\frac{6\la+1}{6^\ga} + \frac{6\la-1}{6^\ga} + \frac1{3^{\ga}}=1.
\end{equation}
Consider the probability vector $(p_{\la,i})_{i=0,\ldots,3}$ where
\begin{equation}\label{def-pila}
p_{\la,0}=p_{\la,3} = \frac1{3^{\ga}}, \, p_{\la,1} = \frac{6\la+1}{6^\ga} \text{ and } p_{\la,2} = \frac{6\la-1}{6^\ga}.
\end{equation}
Then $\mula$ is the unique probability measure satisfying
\begin{equation}\label{defmula}
\mula = \sum_{i=0}^3 p_{\la,i}\mula\circ S_i^{-1} .
\end{equation}
\end{definition}

The multifractal analysis of $\mula$  is a classical issue and is based on its $L^q$-spectrum $\tau_\mula$ and its Legendre transform $\tau^*_\mula$. At this point, we only need their definitions; see Section \ref{sec:mula} for comments and recalls on their properties.

\begin{definition}\label{def:Lqspectrum}
The $L^q$-spectrum of $\mula$ is the mapping $\tau_\mula:q\in\R\longrightarrow \R\cup\{-\infty,+\infty\}$ such that the value $\tau_\mula(q)$ is the unique solution to the equation
\begin{equation}\label{eq:Lqspectrum}
\sum_{i=0}^3 p_{\la,i}^q r_i ^{-\tau_\mula(q)} =1.
\end{equation}
Moreover, the Legendre transform of $\tau_\mula$ is the mapping $\tau^*_\mula:\alpha\in[0,+\infty)\to\R$ defined by
\begin{equation}\label{eq:taustar}
\tau^*_\mula(\alpha)=\inf_{q\in \R} (\alpha q -\tau_\mula(q)).
\end{equation}
\end{definition}

We are now ready to state our main result which asserts that the multifractal spectrum of $\fla$ is that of $\mula$ up to a translation and possibly a truncation.

\begin{theorem}\label{spectruma}
Let $\la\in(\frac{\sqrt2}6,\frac56)$ and consider $\mula$, $\tau_\mula$ and $\tau^*_\mula$.
\begin{itemize}[parsep=-0.15cm,itemsep=0.25cm,topsep=0.2cm,wide=0.175cm,leftmargin=0.5cm]
\item[$(i)$] The support of $d_\fla(\alpha)$ is $[\alpha_{\la,\min},1]:=\Big[1-\mfrac{\log(6\la+1)}{\log 6},1 \Big]$.
\item[$(ii)$]
For every $\alpha\in [\alpha_{\la,\min},1]$,
\begin{equation}\label{res-spectrum}
d_\fla(\alpha) = \tau^*_\mula(\alpha+\ga-1).
\end{equation}
In particular, $d_\fla$ is strictly concave on its support, the maximum of $d_\fla$ is 1 and is reached at the exponent
\begin{equation}\label{*hlexp}
\alpha_{\la,\cL} := 1-\frac{\log(36\la^2-1)}{4\log 3+2\log 6}.
\end{equation}
\item[$(iii)$]
Denote by $\mathrm{Gr}(d_\fla)$ the graph of $d_\fla$, and $\cK(\R^2)$ the set of non-empty compact sets of $\R^2$ endowed with the Hausdorff distance. The mapping $\la\in(\frac{\sqrt2}6,\frac56) \mapsto \mathrm{Gr}(d_\fla) \in \cK(\R^2)$ is continuous except at $\la=\frac13$.
\end{itemize}
\end{theorem}

\begin{figure}[ht]
\centering
\includegraphics[scale=0.7]{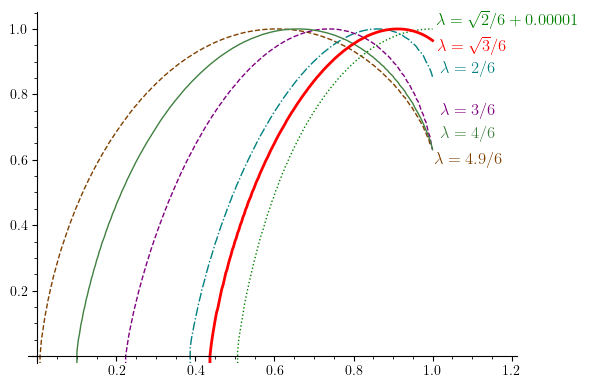}
\captionsetup{width=1\textwidth}
\caption{The multifractal spectrum $d_\fla$ of $\fla$ for different values of $\la\in(\frac{\sqrt2}6,\frac56)$. The spectrum of the original  von Koch function corresponding to the choice $\la=\frac{\sqrt3}6$ is shown by the thick red curve.}
\label{fig:HVKspectra}
\end{figure}

Figure \ref{fig:HVKspectra} displays the multifractal spectrum $d_\fla$ of $\fla$ for various parameter values $\la\in(\frac{\sqrt2}6,\frac56)$.

Although the statement of Theorem \ref{spectruma} is short, it appears that there are significant differences between two parts of the parameter range. Let us gather some remarks on this.

\medskip

\noindent{\bf Comments on Theorem \ref{spectruma}.}
\begin{itemize}[parsep=-0.15cm,itemsep=0.25cm,topsep=0.2cm,wide=0.175cm,leftmargin=0.85cm]
\item[$(a)$] From the statement, the multifractal spectrum of $\fla$ is increasing on $[\alpha_{\la,\min},\alpha_{\la,\cL}]$ and decreasing on $[\alpha_{\la,\cL},1]$. The value $d_\fla(1)$ is always strictly positive, while $d_\fla(\alpha_{\la,\min})$ is always $0$. In addition, Lebesgue-almost every $x\in [0,1]$ is such that $h_\fla(x)=\alpha_{\la,\cL}$. Finally, there is no point $x$ with a pointwise exponent strictly greater than $1$, reinforcing the initial result by von Koch that $\fla$ is nowhere differentiable.
\item[$(b)$] When $\la$ increases to $\frac56$, the abscissa of the leftmost point of the spectrum tends to $0$. This is consistent with Theorem \ref{theo-finite}, since $\fla$ does not converge when $\la=\frac56$. Also, the ordinate of the rightmost point is non-decreasing with $\la$ and converges to $1$ when
$\la\searrow\frac{\sqrt2}6$.
\item[$(c)$] The two mappings $\la\mapsto\alpha_{\la,\min}$ and $\la \mapsto \alpha_{\la,\cL}$ are decreasing with $\la$. This is natural since, as said before, the larger $\la$ is, the more irregular $\fla$ is, i.e. the lowest and the Lebesgue almost sure Hölder exponents both decrease as $\la$ grows.
\item[$(d)$] Surprisingly, the analysis reveals that there is a discontinuity in the mapping $\la\mapsto\mathrm{Gr}(d_\fla)$ at $\la=\frac13$. When $\la>\frac13$, $d_\fla(1)$ is constant, it equals $ \frac{\log 2}{\log3}$, and a phase transition occurs at $\la=\frac13$, where $d_\fla(1)$ jumps to $s$, to the unique solution to
\begin{equation}\label{def-dim-s}
2\cdot 3^{-s}+6^{-s}=1.
\end{equation}
Let us give the reason for this phenomenon. When $\la>\frac13$, the largest local dimension of $\mula$ equals $\frac{-\log p_{\la,0}}{\log 3}$. After some computations, see Section \ref{sec:spectrum}, we can show that the maximum of the support of $d_\fla(\alpha)=\tau^*_\mula(\alpha+\ga-1)$ equals $1$, and its value at $1$ is constant, it equals $\frac{\log2}{\log3}$. Intuitively, this follows from the fact the points of the triadic Cantor set (recall Remark \ref{rk1}) all satisfy $h_\fla(x)=1$. Then, precisely at $\la=\frac13$, $\frac{-\log p_{\la,0}}{\log 3}=\frac{-\log p_{\la,2}}{\log 6}$, and this coincidence implies that $\mula$ has the same local dimension not only on the triadic Cantor set, but on the larger inhomogeneous Cantor set obtained by iteratively keeping, from any interval $[a,b]$, three out of the four subintervals $[a,a+\frac13(b-a)]$, $[a+\frac13(b-a),a+\frac12(b-a)]$, $[a+\frac12(b-a),a+\frac23(b-a)]$, and $[a+\frac23(b-a),b]$, and removing either $(a+\frac13(b-a),a+\frac12(b-a))$ or $(a+\frac12(b-a),a+\frac23(b-a))$ (both of length $\frac16(b-a)$), depending on the monotonicity of the approximating function $\flan$ on $[a,b]$. This set has dimension given by the solution $s$ to \eqref {def-dim-s}, and now the set of points $x$ satisfying $h_\fla(x)=1$ has Hausdorff dimension $s\approx0.8533>\frac{\log2}{\log3}$. This justifies the ``jump'' of the graph of the spectrum $\mathrm{Gr}(d_\fla)$ at $\la=\frac13$ and item $(iii)$ of Theorem \ref{spectruma}.
\item[$(e)$] Let us emphasize that although the statement is the same for every $\la\in(\frac{\sqrt2}6,\frac56)$, the analysis is much more involved in the range $(\frac{\sqrt2}6,\frac13]$ (which includes the classical von Koch function) than in the range $(\frac13,\frac56)$. This is due to the fact that when $\la>\frac13$, the slopes of the functions $\flan$ approximating $\fla$ after $n$ steps of the construction tend to infinity (almost) everywhere, while this is not the case when $\la\leq \frac13$. This part deserves a specific treatment in Section \ref{sec:increasespec}.
\end{itemize}

\medskip

Let us finish this introduction by some remarks on the other range of parameters $\la\in(0,\frac{\sqrt 2}6]$. Precisely at $\la=\frac{\sqrt2}{6}$, the Lebesgue-almost sure exponent $\alpha_{\la,\cL}$ equals $1$, so when $\la$ goes below this value, it is expected that the decreasing part disappears. In \cite{BDS2}, we investigate the range $\la\in (\frac16,\frac{\sqrt2}6)$ and we prove that the equality \eqref{res-spectrum} does not hold when $\alpha$ is close to $1$, and that the right hand-side term is only a lower bound for $d_\fla(\alpha)$ when $\alpha \in [\alpha_{\la,\min},1]$ (we also prove that $d_\fla$ is continuous at $\alpha=1$). The exact spectrum $d_\fla$ for these $\la$s remains an open question.

When $\la$ becomes smaller than $\frac16$, our approach is no longer relevant, since the weight $p_{\la,2}$ in the construction of $\mula$ would become negative. As mentioned above, it corresponds to the critical value when it may happen that the situation of Figure \ref{fig:OmegaLambdaBis} occurs. Then the dynamical system $([0,1],\cB([0,1]),T)$ is not relevant any more to describe the pointwise behavior of $\fla$ any more, and alternative methods must then be developed.

We emphasize that the limit case $\la=\frac16$ is of particular interest, if only for historical reasons. In the third edition of his {\it Traité d'analyse} published in 1922, {\'E}. Picard presented his own example of a continuous but nowhere differentiable function, see \cite[pp. 242--247]{Picard}. This function is based on von Koch's construction except that right triangles are erected at each stage, so that its representative curve admits ``corners'' everywhere. In fact, this function coincides with $F^{\frac16}$. In 1944, G. Bouligand returned to the subject of continuous but nowhere differentiable functions in a book he wrote on the role of intuition in mathematics. He gave the same construction as Picard, see \cite[pp. 222--223]{Bou44}. Curiously, Bouligand gives no proof of his result, and no references to von Koch's or Picard's work.

\medskip

\noindent{\bf Open questions.} Let us indicate some open questions:
\begin{itemize}[parsep=-0.15cm,itemsep=0.25cm,topsep=0.15cm,wide=0.175cm,leftmargin=0.85cm]
\item[$(1)$] What can we say about the (points of) convergence of $\fla$ when $\la\geq\frac56$?
\item[$(2)$] What can we say about the set of points where $\fla$ has infinite derivative when $\la\in[\frac12,\frac23]$? We conjecture that the Hausdorff dimension of the sets $\cV^\la_{\pm\infty}$ depends on $\la$.
\item[$(3)$] What is the value of the multifractal spectrum $d_\fla$ when $\la\in(0,\frac{\sqrt2}6]$? \\ We conjecture that the concavity of $d_\fla$ is preserved at all values of $\la$.
\item[$(4)$] Does $\fla$ satisfy a multifractal formalism?
\item[$(5)$] What is the Hausdorff dimension of the graph of $\fla$?
\end{itemize}

\medskip

\noindent{\bf Outline of the paper.} The paper is organized as follows. In Section \ref{sec:dynam}, we study the dynamical system $([0,1],\cB([0,1]),T)$ and obtain a key decomposition of real numbers $x\in[0,1]$ useful for the analysis of $\fla$. Then, in Section \ref{sec:localholder} we study the slopes of the intermediary steps of the construction of $\fla$, resulting in estimates for the pointwise H\"older exponents.  {Theorem  \ref{theo-finite} is proved in this section.} Section \ref{sec:derivinfini} is devoted to the proof of Theorem \ref{theo-derivinfini}. In particular, its part $(ii)$ will follow from Proposition \ref{propderivinfini}. In Section \ref{sec:spectrum} we collect and prove various regularity properties of $\mula$. The connection between the local dimensions of $\mula$ and the pointwise H\"older exponents of $\fla$ is determined in Section \ref{sec:spectrumirr}. From this a lower estimate is deduced for $d_\fla$, which turns out to be sharp when $\la\in(\frac{\sqrt2}6,\frac56)$. This sharpness is obtained in very different ways depending on the value of $\la$ and on whether the increasing or the decreasing part of the spectrum is considered. In particular, the decreasing part is treated in Section \ref{sec:spectrumirr}, while the increasing part when $\la\in (\frac{\sqrt2}6,\frac13]$ is treated separately in Section \ref{sec:increasespec}. It requires an intricate proof based on the introduction of iterated function systems adapted to the problem.

\addtocontents{toc}{\vspace{0.2cm}}%
\section{The dynamical system and the associated decomposition of $x\in[0,1]$}\label{sec:dynam}

In this section, we study the dynamical system $([0,1],\cB([0,1]),T)$ associated with the recursive construction of $\fla$, where the transformation $T$ was given in \eqref{def-T}.

Let
\begin{equation}\label{chidef}
U(x) =
\begin{cases}
0 & \mbox{ if }       0 \leq x < \frac13, \\
1 & \mbox{ if } \frac13 \leq x < \frac12, \\
2 & \mbox{ if } \frac12 \leq x < \frac23, \\
3 & \mbox{ if } \frac23 \leq x \leq 1,
\end{cases}
\quad\text{and}\quad
\widetilde{U}(x) =
\begin{cases}
      0 & \mbox{ if }       0 \leq x < \frac13, \\
\frac13 & \mbox{ if } \frac13 \leq x < \frac12, \\
\frac23 & \mbox{ if } \frac12 \leq x < \frac23, \\
\frac23 & \mbox{ if } \frac23 \leq x \leq 1.
\end{cases}
\end{equation}

For every $x\in[0,1]$, consider its orbit $(T^n(x))_{n\geq 0}$ under the action of $T$ (with $T^0(x)=x$ by definition). Then, associate with $x$ the two sequences $(u_n(x))_{n\geq 0} \in \Omega = \{0,1,2,3\}^\N$ and $(\widetilde{u}_n(x))_{n\geq 0} \in \{0,\frac13,\frac23\}^\N$ defined as follows:
\begin{equation}\label{def-xn}
u_n(x)=U(T^n(x)) \quad\text{and}\quad\widetilde{u}_n(x)=\widetilde{U}(T^n(x)).
\end{equation}

The sequence $(u_n(x))_{n\geq 0}$ can be thought of as the ``digits'' of $x$ associated with $T$. However, since $T$ is not defined over intervals of the same length, and has a part with a negative slope, we need the sequence $(\widetilde{u}_n(x))_{n\geq 0}$ to obtain a canonical decomposition of $x$ with respect to $T$ using the lower or upper endpoints of the division intervals according to the monotonicity of $T^n$.

\begin{definition}\label{def-indices}
For $x\in[0,1]$, the number of occurrences of the digits of $x$ are given by
\begin{equation}\label{betadef}
\beta_i(x,n) =\#\{k< n : u_k(x)=i\} \ \text{ and }\  \beta_{i,j}(x,n) =\#\{k < n: u_k(x)=i \text{ or } u_k(x)=j\},
\end{equation}
where $i,j\in\{0,1,2,3\}$, $i\neq j$, and $n\geq0$. Notice that $\beta_{i,j}(x,n) = \beta_i(x,n)+\beta_j(x,n)$ and $\beta_{0,3}(x,n)+\beta_{1,2}(x,n) = n$. In particular, $\beta_i(x,0)=\beta_{i,j}(x,0)=0$.

\medskip

Next, we introduce the sequences $(\eps_n(x))_{n\geq0}$, $(\ell_n(x))_{n\geq0}$ and $(a_n(x))_{n\geq0}$ defined as follows:
\begin{equation}\label{lnxdef}
\eps_n(x)=(-1)^{\beta_2(x,n)} \text{, }  \, \ell_n(x)=3^{-\beta_{0,3}(x,n)}6^{-\beta_{1,2}(x,n)}
\end{equation}
\begin{equation*}
\text{and for $n\geq 1$, }
a_n(x)=\sum_{k=0}^{n-1}\widetilde{u}_k(x)\eps_k(x)\ell_k(x).
\end{equation*}
Then $\eps_0(x)=1$, $\ell_0(x)=1$ and by definition $a_0(x)=0$. Moreover, one can easily see that $\eps_{n+1}(x)=\eps_1(T^n(x))\eps_n(x)$.
\end{definition}

\medskip

According to their digits, two subsets of points $x\in[0,1]$ are of particular interest for us.
\begin{definition} \hfill
\begin{itemize}[parsep=-0.15cm,itemsep=0.25cm,topsep=0.2cm,wide=0.175cm,leftmargin=0.75cm]
\item[$(i)$] The countable set $\cE$ consists of those points with digits all equal to $0$ or all equal to $3$ ultimately:
\begin{equation}\label{defE}
\cE= \{x\in [0,1] : u_n(x)   \mbox{ is eventually constant  $0$, or constant $3$} \}.
\end{equation}
\item[$(ii)$] The set $\cEt$ consists of those points with a finite number of 1s and 2s in their decomposition:
\begin{equation}\label{defEtilde}
\cEt = \{x\in [0,1]: \mbox{the sequence $(\beta_{1,2}(x,n))_{n\geq 0}$ is eventually constant} \}.
\end{equation}
\end{itemize}
\end{definition}
\noindent Notice that $\cE\subset\cEt$, and that $\cE$ and $\cEt $ are (strongly) $T$ invariant, i.e. $T(\cE)=\cE$ and $T(\cEt)=\cEt$.

\medskip

We denote by $I=(a,b)$ the open interval consisting of $x$ between $a$ and $b$, even when $b<a$. The corresponding closed interval $[a,b]$ will be denoted by $\overline{I}$, and $|I|$ will stand for the common length of $I$ and $\overline{I}$.

\begin{proposition} \label{dynprop}
Let $x\in [0,1]\setminus\cE$. Then, for every $n\geq0$,
\begin{equation}\label{xint}
x\in I_n(x):=\big(a_n(x),b_n(x)\big) \mbox{ where } b_n(x)=a_n(x)+\eps_n(x)\ell_n(x).
\end{equation}
The sign of the slope of $T^n$ on $I_n(x)$ is $\eps_n(x)$, $T^n(a_n(x))=0$, and $T^n$ maps $I_n(x)$ onto $(0,1)$.
\end{proposition}

\begin{proof}
We proceed by induction on $n\geq0$. Since $I_0(x)=(0,1)$ and $T^0(x)=x$ the result is true when $n=0$. Assume that the statements of Proposition \ref{dynprop} hold for some fixed $n\geq0$.

First, suppose that $\eps_n(x)=1$. Then $T^n$ maps monotone increasingly $I_n(x)=(a_n(x), a_n(x)+\ell_n(x))$ onto $(0,1)$, and $T^n(a_n(x))=0$. Using the definition \eqref{lnxdef} of $a_n(x)$, we get $a_{n+1}(x)=a_n(x)+\widetilde{u}_n(x)\eps_n(x)\ell_n(x)$. In addition to the definition of $I_n(x)$ in \eqref{xint}, we obtain
\begin{equation}\label{inxplus10}
I_{n+1}(x) = \big(a_n(x)+\widetilde{u}_n(x)\eps_n(x)\ell_n(x),a_n(x)+\widetilde{u}_n(x)\eps_n(x)\ell_n(x)+\eps_{n+1}(x)\ell_{n+1}(x)\big).
\end{equation}
The definition of $T$ and \eqref{lnxdef} imply that $T^{n+1}(a_{n+1}(x))=0$ and $\eps_{n+1}(x)=\eps_1(T^n(x))\eps_n(x)=\eps_1(T^n(x))$ is the sign of the slope of $T^{n+1}$ on $I_{n+1}(x)$. The definitions of $\widetilde{U}(x)$ and $\widetilde{u}_n(x)$ in \eqref{chidef} and \eqref{def-xn} and the linearity of $T^n$ on $I_n(x)$ imply that $x\in (a_{n+1}(x),b_{n+1}(x))$.

Suppose now that $\eps_n(x)=-1$. Then $T^n$ maps monotone decreasingly $(a_n(x)-\ell_n(x),a_n(x))$ onto $(0,1)$, and $T^n(a_n(x))=0$. Again, recalling \eqref{lnxdef}, we see that $T^{n+1}(a_{n+1}(x))=0$ and $\eps_{n+1}(x)=\eps_1(T^n(x))\eps_n(x)=-\eps_1(T^n(x))$ is the sign of the slope of $T^{n+1}$ on $I_{n+1}(x)$.
The proof of \eqref{xint} is analogous to the case $\eps_n(x)=1$.
\end{proof}

Let us emphasize that when $x\in \cE$, $I_n(x)$ is well-defined and we check that $x\in \overline{I}_n(x):= \overline{I_n(x)}$.
In addition, observe that $\cE$ coincides with the points of non-differentiability $\flan$ which are the points $a_n(x)$ and $b_n(x)$ for every $x\in[0,1]$ and $n\geq0$.

Also, we see that $\cEt$ is a countable union of affine copies of the triadic Cantor set and in particular $\dimh(\cEt) = \frac{\log 2}{\log 3}$.

For further ease of notation, we introduce the (finite) set $\cI_n$ of all open intervals of generation $n\geq0$, namely
\begin{equation}\label{def-In}
\cI_n= \{I_n(x) : x\in [0,1]\setminus \cE\}.
\end{equation}

\begin{corollary}\label{corodecomp-x}
Every $x\in [0,1]$ is uniquely written
\begin{equation}\label{decomp-x}
x = \sum_{k=0}^\infty \widetilde{u}_k(x) \eps_k(x)\ell_k(x).
\end{equation}
\end{corollary}

\begin{proof}
When $x\notin \cE$, using that  $\ell_n(x)\to 0$ it is a direct consequence of \eqref{xint}. The fact that \eqref{decomp-x} also holds when $x\in \cE$ is left to the reader, and it will not be used in the sequel.
\end{proof}

\addtocontents{toc}{\vspace{0.2cm}}%
\section{Estimates on the pointwise regularity of $\fla$ when $\la\in(\frac16,\frac56)$}\label{sec:localholder}

\subsection{The approximating functions $\flan$}\label{sec:studyFn}

\noindent

This subsection contains technical lemmas that are unavoidable in proceeding to the analysis of $\fla$. Our first key lemma is an analytic transcription of the operation $\Omega^\la$, explaining how $F^{\la}_{n+1}$ is locally obtained from $\flan$ (recall Figure \ref{fig:OmegaLambda} that explains the splitting process).

\begin{lemma}\label{lemslch}
Let $\la\geq \frac16$ and $n\geq0$. Consider an interval $[a,b]$ on which the function $\flan$ is affine with slope $m\in\R$. Set $\ell=b-a>0$. Then, the function $F^{\la}_{n+1}$ is affine on the four consecutive subintervals
\begin{equation}\label{*Ik}
I_0=\big(a,a+\mfrac{\ell}3\big), \ I_1=\big(a+\mfrac{\ell}3,a+\mfrac{\ell}2\big), \ I_2=\big(a+\mfrac{\ell}2,a+\mfrac{2\ell}3\big) \text{ and } I_3=\big(a+\mfrac{2\ell}{3},a+\ell\big),
\end{equation}
with slopes
\begin{equation}\label{*slom}
m, \ m+6\la\sqrt{1+m^2}>0, \ m-6\la\sqrt{1+m^2}<0 \text{ and $m$ respectively.}
\end{equation}
\end{lemma}

\begin{proof}
We only deal with the case $m\geq0$, the other case is left to the reader. Identifying the line segment $AB$ with the graph of $\flan$ we can make use of Figure \ref{fig:OmegaLambda} to make some easy calculations. Since we work with slopes after rescaling and translating we can assume moreover that $(a,b)=(0,1)$ (so $\ell=1$) and that the projection of $AB$ onto the $x$-axis is $[0,1]$ with $A$ being the origin. In addition, by the definition of $\Omega^\la$, it follows that $A=(a,0)=(0,0)$, $B=(b,m)=(1,m)$, $C=(\frac13,\frac{m}3)$, $D=(\frac12,\frac{m}2+\la\sqrt{1+m^2})$, $E=(\frac23,\frac{2m}3)$ and $M=(\frac12,\frac{m}2)$. On $I_0=(0,\frac13)$ the slope of $AC$ equals $m$. On $I_3=(\frac23,1)$ the slope of $EB$ also equals $m$. On $I_1=(\frac13,\frac12)$ the slope of $CD$ is given by
\begin{equation*}
\frac{(\frac{m}2+\la\sqrt{1+m^2})-\frac{m}3}{\frac12-\frac13} = m + 6\la\sqrt{1+m^2}.
\end{equation*}
On $I_2=(\frac12,\frac23)$ the slope of $DE$ is given by
\begin{equation*}
\frac{\frac{2m}3-(\frac{m}2+\la\sqrt{1+m^2})}{\frac23-\frac12} = m-6\la\sqrt{1+m^2}.
\end{equation*}
Finally, for all $\la\geq\tfrac16$ and all $m\in\R$:
\begin{equation*}
\left\{
 \begin{array}{ll}
 m+6\la\sqrt{1+m^2}\geq m+\sqrt{1+m^2} > m+|m|\geq 0, \\
 m-6\la\sqrt{1+m^2}\leq m-\sqrt{1+m^2} < m-|m|\leq 0.
 \end{array}
\right.
\end{equation*}
This concludes the proof.
\end{proof}

Notice that a slope $m>0$ (resp. $m<0$) gives birth to four consecutive new slopes with the sequence of signs $+$, $+$, $-$, $+$ (resp. $-$, $+$, $-$, $-$) whatever $\la\geq\frac16$ is. But when $\la<\frac16$, the sign of the third (resp. the second) slope may change according to the value of the initial slope $m$. Precisely, when $m>0$ is large enough, then the four new slopes are all positive, so the new piecewise part in the iterative construction of $\fla$ is monotonous increasing. This is consistent with the fact that the function $\fla$ is more regular for smaller $\la$.

\medskip

We focus now on the pointwise behavior of $\flan$ at points $x\in[0,1]\setminus\cE$. Denote by $m_n(x)$ the slope of $\flan$ in $x$. One has $m_n(x)=0$ if and only if $\beta_{1,2}(x,n)=0$, and in this case $\eps_n(x)=1$. In particular, this is true for $n=0$. The next lemma connects the dynamical system $([0,1],\cB([0,1]),T)$ to the construction of the function $\fla$. The sign of $x$ is denoted by $\sgn(x)$: $\sgn(x)=1$ if $x>0$, $\sgn(x)=-1$ if $x<0$, and $\sgn(x)=0$ if~$x=0$.

\begin{lemma}\label{lemInx}
Let $\la\geq \frac16$ and $x\in[0,1]\setminus\cE$. For every  $n\geq0$, $I_n(x)$ is exactly the open interval of maximal length containing $x$ and on which the function $\flan$ is affine. In particular,
\begin{equation}\label{slopeFnmnx}
m_n(x)=\frac{\flan(x)-\flan(a_n(x))}{x-a_n(x)}
\end{equation}
and $\sgn(m_n(x))=\sgn(\beta_{1,2}(x,n))\eps_n(x)$.
\end{lemma}

We can rephrase the lemma by saying that the intervals $I\in\cI_n$ of generation $n$ are exactly those on which $\flan$ is affine.

\begin{proof}
We proceed by induction on $n\geq0$. Since $I_0(x)=(0,1)$ and $F^{\la}_0(x)=0$ the result is true when $n=0$.
Suppose now that the statement of Lemma \ref{lemInx} holds for some fixed $n\geq0$. Our strategy is to use Lemma \ref{lemslch} with $(a,b)=I_n(x)$ and \eqref{xint} to identify the interval $I_{n+1}(x)$ with a certain $I_k$ from \eqref{*Ik}. Using \eqref{xint}, we  consider the two cases $\eps_n(x)=1$ and $\eps_n(x)=-1$, and then distinguish four cases according to the value of $u_n(x)\in\{0,1,2,3\}$. Straightforward but useful calculations following from the definitions \eqref{betadef} and \eqref{lnxdef} are collected in the table below.

\renewcommand{\arraystretch}{1.5}
\begin{table}[ht!]
\centering
\begin{tabular}{|c||c||c|c|c|c||c||c|}
\hline
$u_n(x)$ & $\widetilde{u}_n(x)$ & $\beta_2(x,n+1)$ & $\beta_{0,3}(x,n+1)$ & $\beta_{1,2}(x,n+1)$ & $\eps_{n+1}(x)$ & $\ell_{n+1}(x)$ \\
\hhline{|=||=||=|=|=|=||=||=|}
0 & 0 & $\beta_2(x,n)$ & $\beta_{0,3}(x,n)+1$ & $\beta_{1,2}(x,n)$ & $\eps_n(x)$ & $\frac13\ell_n(x)$ \\
\hline
1 & $\frac13$ & $\beta_2(x,n)$ & $\beta_{0,3}(x,n)$ & $\beta_{1,2}(x,n)+1$ & $\eps_n(x)$ & $\frac16\ell_n(x)$ \\
\hline
2 & $\frac23$ & $\beta_2(x,n)+1$ & $\beta_{0,3}(x,n)$ & $\beta_{1,2}(x,n)+1$ & $-\eps_n(x)$ & $\frac16\ell_n(x)$ \\
\hline
3 & $\frac23$ & $\beta_2(x,n)$ & $\beta_{0,3}(x,n)+1$ & $\beta_{1,2}(x,n)$ & $\eps_n(x)$ & $\frac13\ell_n(x)$ \\
\hline
\end{tabular}
\end{table}

\noindent{\it Case 1:} Suppose that $\eps_n(x)=1$. We apply Lemma \ref{lemslch} with $a=a_n(x)$ and $b=a_n(x)+\ell_n(x)$.

\begin{itemize}[parsep=-0.15cm,itemsep=0.25cm,topsep=0.25cm,wide=0.175cm,leftmargin=0.5cm]
\item[$\bullet$] Suppose that $u_n(x)=0$. Using \eqref{inxplus10} with the second line of the table above and using notation $I_j$, $j=0,\ldots,3$ from \eqref{*Ik} we get
\begin{align*}
I_{n+1}(x) = \big(a_n(x),a_n(x)+ \tfrac13\ell_n(x)\big) = I_0.
\end{align*}
In accordance with Lemma \ref{lemslch}, this latter interval is a maximal open interval on which $F^{\la}_{n+1}$ is affine, contains $x$, and the slope of $F^{\la}_{n+1}$ is $m_{n+1}(x) = m_n(x)$. Thus,
\begin{align*}
\sgn(m_{n+1}(x))=\sgn(m_n(x))=\sgn(\beta_{1,2}(x,n))\eps_n(x) = \sgn(\beta_{1,2}(x,n+1))\eps_{n+1}(x).
\end{align*}

\item[$\bullet$] Suppose that $u_n(x)=1$. Using \eqref{inxplus10} with the third line of the table above we obtain
\begin{align*}
I_{n+1}(x)
& = \big(a_n(x)+\tfrac13\ell_n(x),a_n(x)+\tfrac13\ell_n(x)+\tfrac16\ell_n(x)\big) \\
& = \big(a_n(x)+ \tfrac13\ell_n(x),a_n(x)+ \tfrac12\ell_n(x)\big) = I_1.
\end{align*}
By Lemma \ref{lemslch}, $I_1$ is a maximal open interval on which $F^{\la}_{n+1}$ is affine, contains $x$, and the slope of $F^{\la}_{n+1}$ is $m_{n+1}(x)>0$. Hence,
\begin{align*}
\sgn(m_{n+1}(x))= 1 = \sgn(\beta_{1,2}(x,n)+1)\eps_n(x) = \sgn(\beta_{1,2}(x,n+1))\eps_{n+1}(x).
\end{align*}

\item[$\bullet$] Suppose that $u_n(x)=2$. Using \eqref{inxplus10} with the fourth line of the table above we obtain
\begin{align*}
I_{n+1}(x)
& = \big(a_n(x)+\tfrac23\ell_n(x),a_n(x)+\tfrac23\ell_n(x)-\tfrac16\ell_n(x)\big) \\
& = \big(a_n(x)+ \tfrac23\ell_n(x),a_n(x)+ \tfrac12\ell_n(x)\big). 
\end{align*}
Observe that in this case there is an inversion in the relative position of the extremities of $I_2$ due to the negativeness of the slope. Hence, in accordance with Lemma \ref{lemslch}, $I_{n+1}(x)=I_2$ is still a maximal open interval on which $F^{\la}_{n+1}$ is affine, contains $x$, and the slope of $F^{\la}_{n+1}$ is $m_{n+1}(x)<0$. Moreover,
\begin{align*}
\sgn(m_{n+1}(x))= -1 = \sgn(\beta_{1,2}(x,n)+1)(-\eps_n(x)) = \sgn(\beta_{1,2}(x,n+1))\eps_{n+1}(x).
\end{align*}

\item[$\bullet$] Finally, suppose that $u_n(x)=3$. Using \eqref{inxplus10} with the fifth line of the table above we get
\begin{align*}
I_{n+1}(x)
& = \big(a_n(x)+\tfrac23\ell_n(x),a_n(x)+\tfrac23\ell_n(x)+\tfrac13\ell_n(x)\big) \\
& = \big(a_n(x)+\tfrac23\ell_n(x),a_n(x)+\ell_n(x)\big) = I_3.
\end{align*}
By Lemma \ref{lemslch}, $I_3$ is a maximal open interval on which $F^{\la}_{n+1}$ is affine, contains $x$, and the slope of $F^{\la}_{n+1}$ is $m_{n+1}(x)=m_n(x)$. Thus,
\begin{align*}
\sgn(m_{n+1}(x))=\sgn(m_n(x))=\sgn(\beta_{1,2}(x,n))\eps_n(x) = \sgn(\beta_{1,2}(x,n+1))\eps_{n+1}(x).
\end{align*}
\end{itemize}

\smallskip
\noindent{\it Case 2:} The case $\eps_n(x)=-1$ is treated similarly by applying Lemma \ref{lemslch} with $a=a_n(x)-\ell_n(x)$ and $b=a_n(x)$, and is left to the reader.
\end{proof}

Next lemma provides a recursive relation for the sequences of slopes $(m_n(x))_{n\geq0}$ when $x\notin\cE$.

\begin{lemma}\label{lemsloprec}
Let $\la\geq \frac16$ and $x\in[0,1]\setminus\cE$. Suppose that there exists a finite $p\geq0$ such that $\beta_{1,2}(x,p)=0$ and $\beta_{1,2}(x,p+1)=1$. Then, $m_0(x)=\cdots = m_p(x)=0$,
\begin{equation}\label{slopinit}
m_{p+1}(x)
= \left\{
\begin{array}{ll}
\phantom{-}6\la & \text{if $u_p(x)=1$,} \\
-6\la & \text{if $u_p(x)=2$,}
\end{array}
\right.
\end{equation}
and, for all $n\geq p+1$, the slopes $m_n(x)$ are recursively given by
\begin{equation}\label{sloprec}
m_{n+1}(x)
= \left\{
\begin{array}{ll}
m_n(x) & \text{if $u_n(x)\in\{0,3\}$,} \\
\sgn(m_n(x)) \big(|m_n(x)|+6\la\sqrt{1+m^2_n(x)}\big) & \text{if $u_n(x)=1$,} \\
\sgn(m_n(x)) \big(|m_n(x)|-6\la\sqrt{1+m^2_n(x)}\big) & \text{if $u_n(x)=2$.}
\end{array}
\right.
\end{equation}
\end{lemma}

\begin{proof}
The proof is analogous to that of Lemma \ref{lemInx}. If $u_n(x)\in\{0,3\}$, $m_{n+1}(x)=m_n(x)$ for both cases $\eps_n(x)=1$, and $\eps_n(x)=-1$. In particular, since $F^{\la}_0(x)=0$ and $m_0(x)=0$, one has $m_{n+1}(x)=0$ as long as $u_n(x)\in\{0,3\}$, that is, as long as $\beta_{1,2}(x,n+1)=0$. Hence, the first nonzero slope appears in the first step $p\geq0$ such that $u_p(x)\in\{1,2\}$: for this $p$, $\beta_{1,2}(x,n)=0$ for all $n\leq p$ but $\beta_{1,2}(x,p+1)=1$. Then, the value of $m_{p+1}(x)$ is computed using Lemma \ref{lemslch} with $(a,b)=I_p(x)$ and $m=m_p(x) = 0$. The result \eqref{slopinit} comes from the fact that in this case $x\in I_{u_p(x)}(x)$. Finally, let us prove \eqref{sloprec}. For all $n\geq p+1$, $m_n(x)\neq 0$, $\beta_{1,2}(x,n)>0$, and from Lemma \ref{lemInx}, $\sgn(m_n(x))=\eps_n(x)$. We distinguish two cases:

\smallskip

$\bullet$ If $u_n(x)=1$ then
\begin{align*}
m_{n+1}(x)
& = \left\{
\begin{array}{ll}
m_n(x) + 6\la\sqrt{1+m^2_n(x)} & \hbox{if $\eps_n(x)=1$,} \\
m_n(x) - 6\la\sqrt{1+m^2_n(x)} & \hbox{if $\eps_n(x)=-1$}
\end{array}
\right. \\
& = \sgn(m_n(x))|m_n(x)| + \eps_n(x)(6\la\sqrt{1+m^2_n(x)}) \\
& = \sgn(m_n(x))(|m_n(x)| + 6\la\sqrt{1+m^2_n(x)}).
\end{align*}

$\bullet$ If $u_n(x)=2$ then
\begin{align*}
m_{n+1}(x)
& = \left\{
\begin{array}{ll}
m_n(x) - 6\la\sqrt{1+m^2_n(x)} & \hbox{if $\eps_n(x)=1$,} \\
m_n(x) + 6\la\sqrt{1+m^2_n(x)} & \hbox{if $\eps_n(x)=-1$}
\end{array}
\right. \\
& = \sgn(m_n(x))|m_n(x)| - \eps_n(x)(6\la\sqrt{1+m^2_n(x)}) \\
& = \sgn(m_n(x))(|m_n(x)| - 6\la\sqrt{1+m^2_n(x)}).
\end{align*}
This concludes the proof.
\end{proof}

\begin{remark}\label{rkslopeupperbound}
We deduce from Lemma \ref{lemsloprec} that for every $x\in[0,1]\setminus\cE$, the following upper bound holds for the successive slopes: $|m_{n+1}(x)|\leq (1+6\la)|m_n(x)|+6\la$.
\end{remark}

Given a very large slope $|m_n(x)|$ above a point $x$ at step $n$, the next slope $|m_{n+1}(x)|$ can be much smaller and even close to $0$ when $\la$ gets close to $\frac16$. When $\la$ is large enough, this phenomenon does not occur and it may be possible to determine asymptotics for $|m_n(x)|$, at least when $x\notin\cEt$, which is enough for us. The next lemma makes these observations precise.

\begin{lemma}\label{lemslopbound}
Let $\la\geq \frac16$ and $x\in[0,1]\setminus\cE$. If $m_n(x)\neq0$, then $|m_n(x)|\in[m^-_n,m^+_n]$ with $m^-_n=\sqrt{36\la^2-1}$ and $m^+_n\asymp (6\la+1)^n$. Moreover, when $\la>\frac13$, for every $x\notin \cEt$, $\lim_{n\to\infty}|m_n(x)|=+\infty$.
\end{lemma}

\begin{proof}
We first focus on the lower bound $m^-_n$ for any $|m_n(x)|\neq0$. Applying Lemma \ref{lemsloprec}, one sees that in step $n\geq1$, the minimal absolute value of the nonzero slope of $\flan$ is greater than $6\la\sqrt{1+m^2} - |m|$ where $m$ is one of the $4^{n-1}$ slopes of $F^{\la}_{n-1}$. Straightforward calculations show that the function $t\in\R \mapsto 6\la\sqrt{1+t^2} - |t|$ is minimal in $\pm 1/\sqrt{36\la^2-1}$ with minimal value $\sqrt{36\la^2-1}$, hence the value for $m^-_n$.

Let us now estimate the upper bound $m^+_n$. We call $\overline{x}$ the real number whose digits are all equal to $1$, i.e. $u_n(\overline{x})=1$ for every $n$, and consider the function $f(t)=t+6\la\sqrt{1+t^2}$. Applying Lemma \ref{lemsloprec}, it is clear from the construction that $m_n(\overline{x})$ is always the maximal absolute slope of $\flan$. Hence, $m^+_1=m_1(\overline{x})=6\la$ and $m^+_n=m_n(\overline{x})=f^{n-1}(6\la)$. We have
\begin{equation*}\label{ineqmax}
0 < (1+6\la)m \leq f(m)=m+6\la\sqrt{1+m^2} \leq (1+6\la)m+6\la,
\end{equation*}
and, since $f$ is increasing, by iterating this inequality we get for any $n\geq 1$ that
\begin{equation*}\label{ineqmaxbis}
0 < (1+6\la)^n m \leq f^n(m) \leq (1+6\la)^n m + 6\la \left(\mfrac{1-(1+6\la)^n}{1-(1+6\la)}\right).
\end{equation*}
Therefore, putting $m=6\la$ to get the maximum of $f^{n-1}(m)$, we conclude that
\begin{equation*}\label{ineqmaxter}
0 < 6\la(1+6\la)^{n-1} \leq m^+_n \leq (1+6\la)^n-1,
\end{equation*}
which proves the claim about $m^+_n$.

Finally, assume that $\la>\frac13$ and fix $\eps>0$ such that $\la > \frac{2+\eps}6>\frac13$. Let $x\notin \cEt$, which implies that $\beta_{1,2}(x,n) \to\infty$. Proceeding as for $m^-_n$, we have, as soon as $m_n(x)\neq0$ and $u_n(x)=1$ or $u_n(x)=2$,
\begin{align*}
|m_{n+1}(x)|-(1+\eps)|m_n(x)|
\geq 6\la\sqrt{1+m^2_n(x)} -(2+\eps)|m_n(x)|
\geq (6\la-(2+\eps))|m_n(x)|> 0,
\end{align*}
hence
\begin{equation*}\label{ineqmaxfour}
|m_{n+1}(x)|>(1+\eps)|m_n(x)|>0.
\end{equation*}
Each time $u_n(x)=1$ or $u_n(x)=2$, the above inequality can be applied and when $u_n(x)=0$ or $u_n(x)=3$, $m_{n+1}(x)=m_n(x)$.
Since $x\notin\cEt$, the sequence of its digits $u_n(x)$ has to take at least one of the two values $1$ or $2$ for infinitely many $n$, hence $|m_n(x)|\to\infty$.
\end{proof}

We end this section by stating an upper bound for the difference $F^{\la}_{n+1}(z)-\flan(z)$, $z\in[0,1]$, which will be used, in particular, in the proof of Theorem \ref{theo-finite}. This bound is actually valid for every $\la>0$.

\renewcommand{\arraystretch}{1.75}
\begin{lemma}\label{lemmajosup}
Let $\la>0$. Then, for all $n\geq0$, it holds that
\begin{align}\label{majosup}
\sup_{z\in [0,1]} \big|F^{\la}_{n+1}(z)-\flan(z)\big| \leq \la
\left(\mfrac16+\max\big(\la,\mfrac16\big)\right)^n.
\end{align}
\end{lemma}

\begin{proof} Let $n\geq0$. Consider an interval $I_n(x)=(a_n(x),b_n(x))$  of generation $n$ for some $x\in[0,1]\setminus \cE$, on which $\flan$ is affine with slope $m_n(x)$. Assume that $m_n(x)\geq0$ (the case $m_n(x)\leq0$ is treated in a similar way). It follows from the geometric construction of $\flan$ that
\begin{equation}\label{majincFn}
S_n(x):= \sup_{z\in I_n(x)} \big|F^{\la}_{n+1}(z)-\flan(z)\big|= \la L_n(x),
\end{equation}
where $L_n(x)=\ell_n(x)\sqrt{1+m^2_n(x)}$ is the length of the line segment $A_n B_n$, denoted by $|A_n B_n|$, see Figure \ref{fig:MajoSup} below. Defining $L_m(x)$ as $L_n(x)$ by using $I_m(x)$ and $F_m^\la$, we claim that, for all $k\geq0$,
\begin{align}\label{boundLn}
L_{n+k}(x)
\leq\left(\mfrac16+\max\big(\la,\mfrac16\big)\right)^k L_n(x).
\end{align}

\medskip

First, let $L'_{n+1}(x)$ be the length of the line segment $A'_{n+1}B'_{n+1}$ on the graph of $F_{n+1}^\la$ over $I_n(x)$ for which the absolute value of the slope is maximal, see Figure \ref{fig:MajoSup} again.

\begin{figure}[ht]
\includegraphics[width=0.6\textwidth]{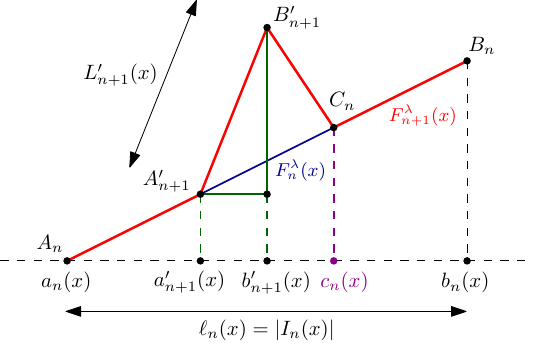}
\captionsetup{width=1\textwidth}
\caption{An upper bound for $\big|F^{\la}_{n+1}(z)-\flan(z)\big|$ over $I_n(x)$.}
\label{fig:MajoSup}
\end{figure}

\noindent Note that the calculation of the maximal slope of $\flan$ that we performed in Lemma \ref{lemslch} is in fact true whenever $\la>0$. Hence, the slope of $F_{n+1}^\la$ over the interval $(a'_{n+1}(x),b'_{n+1}(x))$ is $m'_{n+1}(x)=m_n(x)+6\la\sqrt{1+m^2_n(x)}$ and
\begin{align}\label{majoLn}
(L'_{n+1}(x))^2
& = \big(b'_{n+1}(x)-a'_{n+1}(x)\big)^2\big(1+(m'_{n+1}(x))^2\big) \nonumber\\
& = \big(\mfrac16\ell_n(x)\big)^2\big(1+\big(m_n(x)+6\la\sqrt{1+m^2_n(x)}\big)^2\big) \nonumber\\
& = \mfrac1{36}\ell^2_n(x)\big(1+m^2_n(x)+12\la m_n(x)\sqrt{1+m^2_n(x)}+36\la^2(1+m^2_n(x))\big).
\end{align}

\noindent $\bullet$ Assume that $\la\geq\frac13$. We claim that
\begin{equation}\label{geoLn}
L_{n+1}(x)\leq L_{n+1}'(x) \leq(\tfrac16+\la)L_n(x).
\end{equation}
Indeed, it follows from \eqref{majoLn} that
\begin{align}\label{minLnprime}
(L'_{n+1}(x))^2 \geq \la^2\ell^2_n(x)(1+m^2_n(x))
\geq \big(\mfrac13\ell_n(x)\big)^2(1+m^2_n(x)) =|A_n A'_{n+1}|^2=
|C_n B_n|^2.
\end{align}
Obviously, $L'_{n+1}(x)\geq |B'_{n+1}C_n|$. Thus, independent of whether
$I_{n+1}(x)$ equals $(a_n(x),a'_{n+1}(x))$, $(a'_{n+1}(x),b'_{n+1}(x))$,
$(b_{n+1}'(x),c_n(x))$, or $(c_n(x),b_n(x))$,
we have
\begin{equation}\label{*eqadd}
L'_{n+1}(x)\geq L_{n+1}(x).
\end{equation}
Next, for the upper bound, one has from \eqref{majoLn} that
\begin{align}\label{maxLnprime}
(L'_{n+1}(x))^2
& \leq \mfrac1{36}\ell^2_n(x)(1+12\la+36\la^2) (1+m^2_n(x))=  \big(\mfrac{1+6\la}{6}\big)^2L^2_n(x).
\end{align}
Combining \eqref{*eqadd} and \eqref{maxLnprime} we obtain \eqref{geoLn}. An induction yields \eqref{boundLn}.

\medskip

\noindent $\bullet$ Assume that $\la\in[\frac16,\frac13)$. We claim that
\begin{equation}\label{geoLnBis}
L_{n+1}(x) \leq (\tfrac16+\la)L_n(x).
\end{equation}
If on Figure~\ref{fig:MajoSup} the line segment $A_{n+1}'B_{n+1}'$ is longer than $A_nA'_{n+1}$, that is \eqref{minLnprime} holds, then we can use the previous argument. If not then $L_{n+1}(x)$ is less or equal than the length of $A_nA'_{n+1}$ which is one third of $L_n(x)$. Using that $\la\geq\frac16$ we obtain
\begin{equation*}
L_{n+1}(x)\leq \mfrac13 L_n(x) \leq\left(\mfrac16+\la\right)L_n(x)
\end{equation*}
which yields \eqref{geoLnBis}. An induction gives \eqref{boundLn} again.

\medskip

\noindent $\bullet$ Assume that $\la\in(0,\frac16)$. We claim that
\begin{equation}\label{geoLnTer}
L_{n+1}(x) \leq \tfrac13 L_n(x).
\end{equation}
Indeed, in this last case one obtain in the same way as for \eqref{maxLnprime} that
\begin{align}\label{maxLnprimeTer}
(L'_{n+1}(x))^2
\leq \big(\mfrac13 L_n(x)\big)^2 = |A_nA'_{n+1}|^2.
\end{align}
Thus the line segment $A_{n+1}'B_{n+1}'$ is shorter than $A_nA'_{n+1}$ and we obtain \eqref{maxLnprimeTer} by using the previous argument saying that $L_{n+1}(x)$ is less or equal than the length of $A_nA'_{n+1}$, which is one third of $L_n(x)$. Since $\mfrac16+\max\big(\la,\mfrac16\big) = \frac13$, we deduced  \eqref{boundLn} by induction again.

\medskip

Finally, since $S_n(x)\leq\la L_n(x)$ and the interval $[0,1]$ may be partitioned into a finite number of intervals of the form $I_n(x)$, we obtain \eqref{majosup} using \eqref{majincFn} with $n=0$ and $L_0(x)=1$.
\end{proof}

\renewcommand{\arraystretch}{1.5}
\subsection{Proof of Theorem \ref{theo-finite}}\label{sec:proofTheoConverge}

\noindent

For every $x\in[0,1]$, it follows from Lemma \ref{lemmajosup} that the series $\sum \big|F^{\la}_{n+1}(x)-\flan(x)\big|$ is uniformly convergent for all values of $\la\in(0,\frac56)$. Hence, the existence and continuity of the function $\fla$ (recall that every function $\flan$ is linear, thus continuous).

\medskip

Next, let us prove that $\fla$ is nowhere differentiable. Let $x_0\in[0,1]$.

\begin{figure}[ht]
\includegraphics[width=0.6\textwidth]{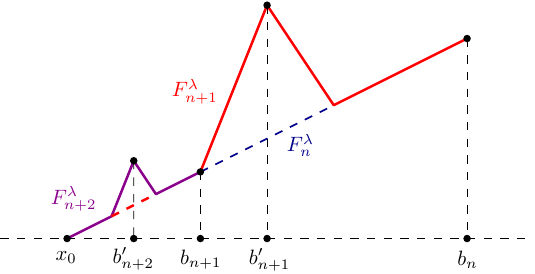}
\captionsetup{width=1\textwidth}
\caption{Non-differentiability of $\fla$ at a point $x_0\in\cE$.}
\label{fig:NondiffA}
\end{figure}

First, assume that $x_0\in\cE$. Then, there exists a minimal integer $n_0\geq0$ such that $x_0$ is a breakpoint of $\flan$ for all $n\geq n_0$. Assume $x_0\neq1$ and denote by $(x_0,b_n)$ the interval of generation $n\geq n_0$ with extremity $x_0$ on which $\flan$ is affine with slope $m_n$. Set $b'_{n+1} = x_0+ \frac12(b_n-x_0)$. By construction of the sequence $(\flan)_{n\geq n_0}$, we have $m_n=m_{n_0}$ for all $n\geq n_0$, see Figure \ref{fig:NondiffA}.

\noindent Assume that $m_{n_0}\geq0$ (the case $m_{n_0}\leq0$ can be treated similarly). Then, for all $n\geq n_0$,
\begin{equation*}
\fla(b_n)-\fla(x_0) = F^{\la}_{n+1}(b_n)-F^{\la}_{n+1}(x_0) = m_{n_0}(b_n-x_0)
\end{equation*}
and
\begin{equation*}
\fla(b'_{n+1})-\fla(x_0) = F^{\la}_{n+1}(b'_{n+1})-F^{\la}_{n+1}(x_0) = (m_{n_0}+2\la\sqrt{1+m^2_{n_0}})(b'_{n+1}-x_0).
\end{equation*}
Since $b_n,b'_{n+1}\to x_0$ and $m_{n_0}\neq m_{n_0}+2\la\sqrt{1+m^2_{n_0}}$, the function $\fla$ cannot be differentiable at $x_0$. Finally, the case $x_0=1$ is a consequence of the relation $\fla(1-x)=\fla(x)$ (see \eqref{equafunc}).

\begin{figure}[ht]
\includegraphics[width=0.6\textwidth]{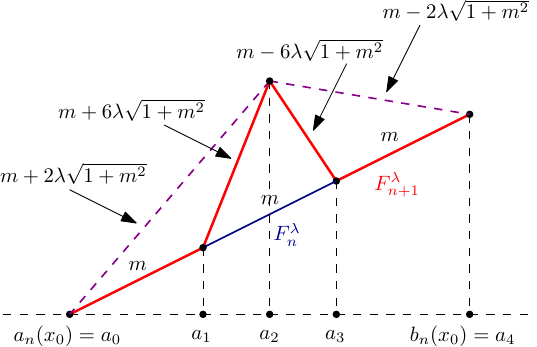}
\captionsetup{width=1\textwidth}
\caption{Non-differentiability of $\fla$ at a point $x_0\notin\cE$.}
\label{fig:NondiffB}
\end{figure}

Assume now that $x_0\notin\cE$. Assuming $\fla$ is differentiable at $x_0$, the mapping
\begin{equation*}
D:(h,h')\in(0,\infty)^2\to \frac{\fla(x_0+h)-\fla(x_0-h')}{h+h'}
\end{equation*}
admits a finite limit at $(0,0)$ equal to $(\fla)'(x_0)$. Consider the interval $I_n(x_0)=(a_n(x_0),b_n(x_0))$ of generation $n$ that contains $x_0$ on which $\flan$ is affine with slope $m_n(x_0)$. Suppose, for example, that $m_n(x_0)\geq0$ (the case $m_n(x_0)\leq0$ may be treated similarly). Set $m=m_n(x_0)$, $a_0 = a_n(x_0)$, $a_1=a_0+\frac13\ell_n(x_0)$, $a_2=a_0+\frac12\ell_n(x_0)$, $a_3=a_0+\frac23\ell_n(x_0)$ and $a_4=a_0+\ell_n(x_0)=b_n(x_0)$, see Figure \ref{fig:NondiffB}.

\medskip

We claim that we can always select points $a_{i_1},a_{j_1}\in[a_0,x_0)$ and $a_{i_2},a_{j_2}\in(x_0,a_4]$ such that
\begin{equation}\label{eq:diffD}
\left|\frac{F^{\la}_{n+1}(a_{j_2})-F^{\la}_{n+1}(a_{j_1})}{a_{j_2}-a_{j_1}} - \frac{F^{\la}_{n+1}(a_{i_2})-F^{\la}_{n+1}(a_{i_1})}{a_{i_2}-a_{i_1}}\right| \geq 2\la>0.
\end{equation}
We indicate in the table below how to select $a_{i_1},a_{i_2}$ and $a_{j_1},a_{j_2}$ according to the position of $x_0$ inside $I_n(x_0)$, and the values for the corresponding increments $\rho_{i_1,i_2}$ and $\rho_{j_1,j_2}$ of $F^{\la}_{n+1}$ needed in \eqref{eq:diffD}.

\renewcommand{\arraystretch}{1.5}
\begin{table}[ht!]
\centering
\begin{tabular}{|c||c|c|c||c|c|c||c||c|}
\hline
$x_0\in$  & $a_{i_1}$ & $a_{i_2}$ & $\rho_{i_1,i_2}$ & $a_{j_1}$ & $a_{j_2}$ & $\rho_{j_1,j_2}$ & $|\rho_{j_1,j_2}-\rho_{i_1,i_2}|$ \\
\hhline{|=||=|=|=||=|=|=||=||}
$(a_0,a_1)$ & $a_0$ & $a_1$ & $m$ & $a_0$ & $a_2$ & $m+2\la\sqrt{1+m^2}$ & $2\la\sqrt{1+m^2}$ \\
\hline
$(a_1,a_2)$ & $a_1$ & $a_2$ & $m+6\la\sqrt{1+m^2}$ & $a_1$ & $a_3$ & $m$ & $6\la\sqrt{1+m^2}$ \\
\hline
$(a_2,a_3)$ & $a_1$ & $a_3$ & $m$ & $a_2$ & $a_3$ & $m-6\la\sqrt{1+m^2}$ & $6\la\sqrt{1+m^2}$ \\
\hline
$(a_3,a_4)$ & $a_2$ & $a_4$ & $m-2\la\sqrt{1+m^2}$ & $a_3$ & $a_4$ & $m$ & $2\la\sqrt{1+m^2}$ \\
\hline
\end{tabular}
\end{table}

\medskip

\noindent The lower bound in \eqref{eq:diffD} then holds, since $6\la\sqrt{1+m^2}\geq 2\la\sqrt{1+m^2}\geq 2\la$.

\medskip

Since $F^{\la}_{n+1}$ and $\fla$ coincides at points $a_k$, we can rewrite \eqref{eq:diffD} as
\begin{equation*}\label{eq:diffDbis}
\left| D(a_{j_2}-x_0,x_0-a_{j_1}) - D(a_{i_2}-x_0,x_0-a_{i_1}) \right| \geq 2\la.
\end{equation*}
Since $a_k\to x_0$, the function $D$ does not satisfy Cauchy's condition at $(0,0)$ and then cannot have a limit at this point. Hence, a contradiction.

\medskip

Finally, let $\la\geq \frac56$, and consider again the special point $\overline{x}\in[0,1]$ whose digits are all $u_n(x)=1$. We call $\overline{x}_n=0.11\cdots11$ its approximation with $n$ digits equal to $1$ and $y_n=\fla(\overline{x}_n)=\flan(\overline{x}_n)$. By construction, for all $n\geq1$, one has $\ell_n(\overline{x})=\frac1{6^n}$, and $\overline{x}_n=\frac13(1+\frac16+\cdots+\frac1{6^n})$. The slopes $m_n(\overline{x})$ satisfy $m_{n+1}(\overline{x})=m_n(\overline{x})+6\la \sqrt{1+m^2_n(\overline{x})}$. From $m_1(\overline{x})=6\la\geq 5 $, we deduce that $m_n(\overline{x})$ tends to $\infty$, and in particular
\begin{equation*}
m_{n+1}(\overline{x}) \geq (1+6\la)m_n(\overline{x}) \geq 6m_n(\overline{x}).
\end{equation*}
Again, since $m_1(\overline{x}) \geq 5 $, we conclude that for some positive constant $C>0$, $m_n(\overline{x}) \geq C\cdot 6^n$. Finally, observe that $y_{n+1} = y_n + (\overline{x}_{n+1}-\overline{x}_n)m_{n+1}(\overline{x})$. Combining this, the fact that $\overline{x}_{n+1}-\overline{x}_n=\frac1{3\cdot 6^{n+1}}$ with the previous inequality on $m_n(\overline{x})$, yields $y_{n+1} - y_n \geq C'$ for some positive constant $C'>0$. In particular, $\lim_{n\to \infty} y_n=+\infty$, and $\fla(x)=+\infty$. This concludes the proof of Theorem \ref{theo-finite}.

\begin{remark}\label{rk2}
Observe that the last part of the proof holds true for any $x\in[0,1]$ whose digits are all ultimately $1$, and actually for many more points. It would be interesting to characterize the set of convergence of $\fla$ in this range of parameters.
\end{remark}

\subsection{The increments of $\fla$}\label{sec:studyF}

\noindent

We begin this section with a lemma which will be useful in the sequel. It aims to compare $I_n(x)$ and $I_n(y)$ when $x$ and $y$ are close enough.

\begin{lemma}\label{lem-neighbors}
Let $\la\in(\frac16,\frac56)$ and $n\geq0$. Let $x,y\in[0,1]\setminus\cE$ such that $|x-y|\leq \frac12\ell_n(x)$. Then, either $I_n(x)=I_n(y)$ or $I_n(x)$ and $I_n(y)$ are contiguous. Moreover, $\frac12\ell_n(y)\leq \ell_n(x) \leq 2\ell_n(y)$ and
\begin{equation}\label{majmny}
|m_n(y)|\leq (1+6\la)|m_n(x)|+ 6\lambda\leq 6(|m_n(x)|+1).
\end{equation}
\end{lemma}

\begin{proof}
If $y\in I_n(x)$ then $I_n(y)=I_n(x)$, hence $\ell_n(y)=\ell_n(x)$ and $m_n(y)=m_n(x)$ so the result is obvious. Assume that $y\notin I_n(x)$.
Without limiting generality, we can assume that $x<y$. Consider the unique integer $0\leq p\leq n-1$ such that $I_p(x)=I_p(y)$ but $I_{p+1}(x)\neq I_{p+1}(y)$. The construction process implies that
\begin{equation}\label{*eqfel}
\frac12\ell_{p+1}(y)\leq \ell_{p+1}(x)\leq  2 \ell_{p+1}(y).
\end{equation}
Obviously $\ell_n(x)\leq \ell_{p+1}(x)\leq \ell_p(x)/3$ and any interval of generation $p+1$ included in $I_p(x) $ is of length at least $\ell_p(x)/6$. Since by assumption $|x-y|\leq \frac12\ell_n(x)$, $x$ and $y$ must belong to consecutive intervals of generation $p+1$. Any interval $I\in\cI_n$ contained in $I_{p+1}(x)$ is of length at most $\ell_{p+1}(x)/3^{n-p-1}$ and the ones at the extremities of $I_{p+1}(x)$ are exactly of this length. So, if there is a point $y\not \in I_p(x)$ at distance $|x-y|<\ell_n(x)/2<\ell_n(x)$, then $I_n(x)$ must be the rightmost subinterval of $I_{p+1}(x)$ of generation $n$, and $\ell_n(x)=\ell_{p+1}(x)/3^{n-p-1}$. Similarly, we see that $I_n(y)$ should be the leftmost subinterval in $\cI_n$ of the generation $n$ included in $I_{p+1}(y)$ of length $\ell_{p+1}(y)/3^{n-p-1}$. Then it follows from \eqref{*eqfel} that $\frac12\ell_n(y)\leq \ell_n(x) \leq 2\ell_n(y)$.
Next, given the relative position of $I_n(x)$ and $I_n(y)$, the construction process implies that
\begin{equation}\label{*eqfem}
m_k(y)=m_k(x) \text{ for all $k\in\{0,\ldots,p\}$ but  $m_{p+1}(y)\neq m_{p+1}(x)$}.
\end{equation}
It remains to prove \eqref{majmny} for $n=p+1$. Noting that $m_p(y)=m_p(x)$, it is a direct consequence of Lemma \ref{lemsloprec} (see also Remark \ref{rkslopeupperbound}). Since $x$ and $y$ are close to the intersection point of $I_k(x)$ and $I_k(y)$ for $k=p+1,\ldots,n$, one concludes that $m_k(x)=m_{p+1}(x)$ and $m_k(y)=m_{p+1}(y)$ for $k=p+1,\ldots,n$.
\end{proof}

Next, we compute $h_\fla(x_0)$ for $x_0\notin\cE$. To this end, first we give an upper bound of the local increments of $\fla$ using the slopes of $\flan$.

\begin{proposition}\label{prop:lipabove}
Let $\la\in(\frac16,\frac56)$ and $n\geq0$. Let $x,y\in[0,1]$ such that
\begin{equation}\label{*eqlnx}
\frac1{12}\ell_n(x)\leq |x-y|\leq\frac12\ell_n(x).
\end{equation}
Assume that $x\notin\cE$. Then we have
\begin{equation}\label{lipabove}
|\fla(x)-\fla(y)|\leq C(|m_n(x)|+1)|x-y|,
\end{equation}
where the constant $C>0$ does not depend on $n$, nor on the points $x,y$.
\end{proposition}

\begin{proof} First, we deal with the case $y\in I_n(x)$. Since $\flan$ is affine on $I_n(x)$ with slope $m_n(x)$,
\begin{align}\label{majincF}
|\fla(x)-\fla(y)|
& \leq |\flan(x)-\flan(y)| + |(\fla-\flan)(x)-(\fla-\flan)(y)| \nonumber \\
& \leq |m_n(x)||x-y| + 2 \sup_{z\in I_n(x)} \big|\fla(z)-\flan(z)\big|.
\end{align}
Writing $\fla(z)-\flan(z) = \sum_{k=n}^{\infty} F^{\la}_{k+1}(z)-F^{\la}_k(z)$, it follows from \eqref{majincFn} and \eqref{boundLn} that
\begin{equation}\label{majsupin}
\sup_{z\in I_n(x)} \big|\fla(z)-\flan(z)\big|\leq \la\sum_{k=n}^{\infty} L_k(x) \leq \la\sum_{k=0}^\infty L_{k+n}(x)\leq C_1 L_n(x)
\end{equation}
with $C_{1}>0$. Moreover,
\begin{equation}\label{majln}
L_n(x) = \ell_n(x)\sqrt{1+m^2_n(x)} \leq \ell_n(x)(1+|m_n(x)|) \leq 12|x-y|(1+|m_n(x)|)
\end{equation}
which yields \eqref{lipabove} by combining \eqref{*eqlnx}, \eqref{majincF}, \eqref{majsupin} and \eqref{majln}.

Suppose now that $y\notin I_n(x)$. From Lemma \ref{lem-neighbors}, the interval $I_n(y)$ is contiguous to $I_n(x)$. We assume that $I_n(y)$ precedes $I_n(x)$ with $m_n(x)\geq0$ (the other cases are treated similarly). We proceed as for the first case but using the intermediate point $a_n(x)$ by writing
\begin{align}\label{majincFbis}
|\fla(x)-\fla(y)| \leq |\fla(x)-\fla(a_n(x))|+|\fla(y)-\fla(a_n(x))|.
\end{align}
Noticing that $\fla(a_n(x))= \flan(a_n(x))$, we obtain, with  \eqref{*eqlnx}, \eqref{majsupin} and \eqref{majln},
\begin{align}\label{majincFbisx}
|\fla(x)-\fla(a_n(x))|
& \leq |\flan(x)-\flan(a_n(x))| + |\fla(x)-\flan(x)| \nonumber \\
& \leq |m_n(x)||x-a_n(x)| + C_1 L_n(x) \nonumber \\
& \leq |m_n(x)| \ell_n(x) + C_1 \ell_n(x)(1+|m_n(x)|) \nonumber \\
& \leq C'_1 (1+|m_n(x)|)|x-y|.
\end{align}
Similarly, we get
\begin{align}\label{majincFbisy}
|\fla(y)-\fla(a_n(x))|
& \leq |m_n(y)| \ell_n(y) + C_1 \ell_n(y)(1+|m_n(y)|).
\end{align}

\noindent Using Lemma \ref{lem-neighbors} to compare $\ell_n(y)$ with $\ell_n(x)$, and $m_n(y)$ with $m_n(x)$, \eqref{majincFbisy} implies that
\begin{align}\label{majincFbisybis}
|\fla(y)-\fla(a_n(x))| \leq (1+6C_1)\ell_n(x)(1+|m_n(x)|) \leq C''_1 |x-y|(1+|m_n(x)|).
\end{align}
Finally, by combining \eqref{majincFbis}, \eqref{majincFbisx} and \eqref{majincFbisybis} we obtain the desired inequality \eqref{lipabove}.
\end{proof}

\begin{proposition}\label{prop:lipbelow}
Let $\la\in(\frac16,\frac56)$ and $n\geq0$. Then, for all $x\in[0,1]\setminus\cE$,
\begin{equation}\label{lipbelow}
|\fla(x)-\fla(a_n(x))|\geq |m_n(x)||x-a_n(x)|.
\end{equation}
\end{proposition}

\begin{proof} Recall the definition of $\cE$ given in Definition \ref{def-indices}. The result is obvious when $m_n(x)=0$. Suppose that $m_n(x)\neq0$. Since $I_n(x)=(a_n(x),a_n(x)+\eps_n(x)\ell_n(x))$, it follows from Lemma \ref{lemInx} that $\sgn(m_n(x))=\sgn(x-a_n(x))$. Thus, $\flan(x)-\flan(a_n(x))>0$ and then $\fla(x)\geq \flan(x)> \flan(a_n(x))=\fla(a_n(x))$. Therefore, we obtain successively
\begin{align*}\label{minincF}
|\fla(x)-\fla(a_n(x))|
& = \fla(x)-\fla(a_n(x)) = \fla(x)-\flan(a_n(x)) \nonumber \\
& \geq \flan(x)-\flan(a_n(x)) = m_n(x)(x-a_n(x)) = |m_n(x)||x-a_n(x)|. \nonumber
\end{align*}
That concludes the proof.
\end{proof}

\subsection{Estimates on the pointwise H\"older exponents of $\fla$}\label{sec:basicresults}

\noindent

The first proposition concerns the pointwise H\"older exponent of $\fla$.

\begin{proposition}\label{propexphmax}
Let $\la\in(\frac16,\frac56)$. For every $x_0\in[0,1]$ we have $h_\fla(x_0)\leq1$, and $h_\fla(x_0)=1$ if $x\in\cEt$.
\end{proposition}

\begin{proof}
By Theorem \ref{theo-finite} the function $\fla$ is nowhere differentiable, hence $h_\fla(x_0)\leq1$ for all $x_0\in[0,1]$.

Let $x_0\in\cEt$.

First, assume that $x_0\in\cE$, so $x_0$ is a breakpoint of $\flan$ for some $n\geq1$. We call $[b',x_0],[x_0,b]\in\cI_n$ the two consecutive intervals of generation $n$ with extremity $x_0$, on which $\flan$ is affine. Let us deal with the interval $[x_0,b]$, the case $[b',x_0]$ is similar. Observe that from its construction, $\fla$ enjoys a sort of ``limited self-similarity'' property (see \eqref{equafunc}):
\begin{equation}\label{limssim}
\fla(x_0+3(x-x_0))-\fla(x_0)=3(\fla(x)-\fla(x_0))
\end{equation}
holds for all $x\in(x_0,x_0+\frac{b-x_0}3]$. Denoting by $M>0$ the maximum of $\fla(x)$ on $[0,1]$, we obtain
\begin{equation}\label{estF}
\left|\frac{\fla(x)-\fla(x_0)}{x-x_0}\right|\leq \frac{3\cdot 2M}{b-x_0}
\end{equation}
for all $x\in(x_0+\frac{b-x_0}3,b]$. Using repeatedly \eqref{limssim}, for any $n\geq 1$, if $x\in\big(x_0+\tfrac{b-x_0}{3^n},x_0+\tfrac{b-x_0}{3^{n-1}}\big]$ then \eqref{estF} holds. This implies that the absolute values of the right Dini derivatives of $\fla$ at $x_0$ are bounded. A similar argument shows the same result for the left Dini derivatives of $\fla$. Thus $h_{\fla}(x_0)\geq 1$, hence $h_{\fla}(x_0)=1$.

Assume now that $x_0\in\cEt\setminus\cE$. By construction, there exists $n_{x_0}\geq1$ such that for every $n\geq n_{x_0}$, $m_n(x_0)=m\in\R$. For every $y$ such that $|y-x_0|\leq \ell_{n_{x_0}}(x_0)/12$, there exists (at least) one integer $n$ such that \eqref{*eqlnx} holds true (this follows from the fact that $\ell_n(x)/6\leq \ell_{n+1}(x)\leq \ell_n(x)/3$). Then Proposition \ref{prop:lipabove} applies, and $|\fla(x_0)-\fla(y)|\leq C(|m_n(x_0)|+1)|x_0-y|= C(m+1)|x_0-y|$. This gives $h_\fla(x_0)\geq 1$, which allows us to conclude.
\end{proof}

In particular, it follows that the pointwise H\"older exponent of $\fla$ at every $x\in [0,1]$ is conveniently given by the formula
\begin{equation}\label{defpointbis}
h_\fla(x)= \liminf_{y\to x}\frac{\log |\fla(y)-\fla(x)|}{\log |y-x|}.
\end{equation}

Now, with Propositions \ref{prop:lipabove}, \ref{prop:lipbelow} and \ref{propexphmax} we are able to obtain a first expression for $h_\fla(x)$.

\begin{proposition}\label{propholda}
Let $\la\in(\frac16,\frac56)$. For every $x\in[0,1]\setminus\cEt$, we have
\begin{equation}\label{expholda}
h_\fla(x) = 1 - \limsup_{n\to\infty}\frac{\log |m_n(x)|}{|\log\ell_n(x)|}.
\end{equation}
\end{proposition}

\begin{proof} Fix such an $x$.
Since $x\not \in \cEt$ we have $m_{n}(x)=0$ only for finitely many $n$s.
We set $\alpha_n(x) =\log (\ell_n(x)|m_n(x)|)/\log\ell_n(x)$ and $\alpha(x)=\liminf_{n\to\infty} \alpha_n(x)$. We have to prove that $h_\fla(x)=\alpha(x)$, which is equivalent to \eqref{expholda}.

We first show that $h_\fla(x)\geq\alpha(x)$. Let $\eps>0$. By the definition of $\alpha(x)$, there exists an integer $N\geq 1$ such that, for all $n\geq N$, it holds that $\alpha_n(x)>\alpha(x)-\eps$. This provides $|m_n(x)|\leq \ell_n(x)^{\alpha(x)-\eps-1}$. Let $\delta = \frac12\ell_N(x)>0$. Let $y\in[0,1]\setminus\cEt$ such that $|x-y|<\delta$. We can find $n\geq N$ such that $\frac1{12}\ell_n(x)\leq |x-y|\leq \frac12\ell_n(x)$. Hence, by using Proposition \ref{prop:lipabove}, we obtain, with absolute constants $C,C'>0$,
\begin{equation*}
|\fla(x)-\fla(y)|\leq C|m_n(x)||x-y|\leq C' |x-y|^{\alpha{(x)}-\eps}.
\end{equation*}
Therefore, $h_\fla(x)\geq\alpha(x)-\eps$. Since this is true for every $\eps>0$ we get $h_\fla(x)\geq\alpha(x)$.

Next, we prove that $h_\fla(x)\leq\alpha(x)$. Let $\eps>0$ and $\delta>0$. Let $N\geq 1$ be an integer such that $\ell_N(x)\leq \delta$. By the definition of $\alpha(x)$, there exists an integer $n\geq N$ such that $\alpha_n(x)<\alpha(x)+\eps$. This provides $|m_n(x)|\geq \ell_n(x)^{\alpha(x)+\eps-1}$. Hence, by using Proposition \ref{prop:lipbelow}, we obtain the existence of a point $y=a_n(x)$ satisfying $|x-y|\leq \ell_n(x)<\delta$ and such that
\begin{equation*}
|\fla(x)-\fla(y)|\geq |m_n(x)||x-y|\geq |x-y|^{\alpha(x)+\eps}.
\end{equation*}
Therefore, $h_\fla(x)\leq\alpha(x)+\eps$. Since this is true for every $\eps>0$ we get $h_\fla(x)\leq\alpha(x)$.
\end{proof}

The next step is to express $h_\fla(x)$ in terms of the frequencies of the digits of $x$ through successive slopes $m_n(x)$. Since a precise result will be obtained when $\lim_{n\to\infty}|m_n(x)|=\infty$, we introduce the following particular set of points.

\begin{definition}\label{*defI}
We denote by $\cI$ the set of those $x\in [0,1]\setminus\cEt$ for which
\begin{equation*}\label{defmninfinity}
\lim_{n\to\infty}|m_n(x)| = +\infty.
\end{equation*}
\end{definition}

\begin{proposition}\label{propholdabis}
Let $\la\in(\frac16,\frac56)$. For every $x\in[0,1]\setminus\cEt$, we have
\begin{equation}\label{hexpa}
h_\fla(x) \leq 1-\max \bigg(\limsup_{n\to\infty}\frac{\beta_1(x,n)\log(6\la+1)+\beta_2(x,n)\log(6\la-1)}{\beta_{0,3}(x,n)\log3 +\beta_{1,2}(x,n)\log 6},0\bigg),
\end{equation}
and equality holds if $x\in\cI$.
\end{proposition}

\begin{proof} Fix $x\in[0,1]\setminus\cEt$. If the $\limsup$ in \eqref{hexpa}
is negative or zero, then by Lemma \ref{propexphmax} $h_\fla(x)\leq 1$, and the result is obvious. Assume that the $\limsup$ is positive. By Proposition \ref{propholda}, we have to prove that
\begin{equation}\label{hexpb}
\limsup_{n\to\infty}\frac{\log |m_n(x)|}{|\log\ell_n(x)|} \geq
\limsup_{n\to\infty}\frac{\beta_1(x,n)\log(6\la+1)+\beta_2(x,n)\log(6\la-1)}{\beta_{0,3}(x,n)\log3 +\beta_{1,2}(x,n)\log 6},
\end{equation}
and that the equality holds when $x\in\cI$.

Let us introduce the following functions $(f_i)_{i=0,\ldots,3}$ and $(g_i)_{i=0,\ldots,3}$:
\begin{equation}\label{iterfuncfg}
f_i(t)=\left\{
\begin{array}{ll}
\hspace{8mm} t & \hbox{if $i\in\{0,3\}$,} \\
(6\la+1)t & \hbox{if $i=1$,} \\
(6\la-1)t & \hbox{if $i=2$,}
\end{array}
\right.
\text{ and } \
g_i(t)=\left
\{\begin{array}{ll}
\hspace{12mm} 1 & \hbox{if $i\in\{0,3\}$,} \vspace{0.15cm} \\
\vspace{0.2cm}
\displaystyle\frac{6\la\sqrt{1+1/t^2}+1}{6\la+1} & \hbox{if $i=1$,} \\
\displaystyle\frac{6\la\sqrt{1+1/t^2}-1}{6\la-1} & \hbox{if $i=2$}.
\end{array}
\right.
\end{equation}
We also set $g_i(0)=1$ for all $i$. Observe that $g_i(t)\to1$ when $t\to+\infty$. Since $x\notin\cEt$ and the $\limsup$ in \eqref{hexpa} is positive, we may assume without loss of generality that the integer $p$ associated with $x$ in Lemma \ref{lemsloprec} is equal to $p=0$
and $\beta_{1}(x,1)=1$, that is $u_{1}(x)=1$. Then $m_0(x)=0$ and $|m_1(x)|=6\la$. For all $n\geq 1$,
\begin{equation}\label{itermn}
|m_{n+1}(x)|= f_{u_n(x)}(|m_n(x)|)g_{u_n(x)}(|m_n(x)|).
\end{equation}
Due to the functions $g_i$, the sequence of slopes $(m_n(x))_{n\geq0}$ is not a multiplicative process. Our strategy is to compare it with the ``idealized'' sequence $(\tm_n(x))_{n\geq0}$ we would obtain if there were no $g_i$, and show that they have the same asymptotics when $x\in\cI$. Thus, let us define the sequence $(\tm_n(x))_{n\geq0}$ by
$\tm_0(x)=0$, $\tm_1(x)=6\la$ and for all $n\geq 1$,
\begin{equation}\label{itermnideal}
\tm_{n+1}(x)= f_{u_n(x)}(\tm_n(x)).
\end{equation}
Elementary computation gives that for all $n\geq1$,
\begin{equation}\label{seqideal}
\tm_n(x) = 6\la(6\la+1)^{\beta_1(x,n)-1}(6\la-1)^{\beta_2(x,n)}>0.
\end{equation}
Recalling $\ell_n(x)=3^{-\beta_{0,3}(x,n)}6^{-\beta_{1,2}(x,n)}$, see \eqref{lnxdef}, we obtain
\begin{equation}\label{hexpbideal}
\limsup_{n\to\infty}\frac{\log \tm_n(x)}{|\log\ell_n(x)|} =
\limsup_{n\to\infty}\frac{\beta_1(x,n)\log(6\la+1)+\beta_2(x,n)\log(6\la-1)}{\beta_{0,3}(x,n)\log3 +\beta_{1,2}(x,n)\log 6}.
\end{equation}
Therefore, \eqref{hexpb} will be a consequence of
\begin{equation}\label{equivmn}
\limsup_{n\to\infty}\frac{\log |m_n(x)|}{|\log\ell_n(x)|}\geq\limsup_{n\to\infty}\frac{\log \widetilde{m}_n(x)}{|\log\ell_n(x)|}.
\end{equation}
We easily prove by induction that, for all $n\geq 2$,
\begin{equation}\label{itermnbis}
|m_n(x)|= \tm_n(x)\prod_{k=1}^{n-1} g_{u_k(x)}(|m_k(x)|).
\end{equation}
It follows that
\begin{equation}\label{itermnter}
\frac{\log |m_n(x)|}{|\log\ell_n(x)|} = \frac{\log \tm_n(x)}{|\log\ell_n(x)|}
+ \bigg[\frac1{|\log\ell_n(x)|}\sum_{k=1}^{n-1} \log g_{u_k(x)}(|m_k(x)|)\bigg].
\end{equation}
Straightforward calculations yield $g_2(x)\geq g_1(x)\geq g_0(x)=g_3(x)=1$. It follows that the term inside the brackets above is non-negative, which in turn yields \eqref{equivmn}.

Next, using \eqref{lnxdef} and the observation after \eqref{betadef}, we get
\begin{equation*}
|\log \ell_n(x)| = \beta_{0,3}(x,n)\log3 +\beta_{1,2}(x,n)\log 6 \geq (\beta_{0,3}(x,n) +\beta_{1,2}(x,n))\log 3 = n \log 3.
\end{equation*}
Therefore, the following upper bound holds:
\begin{align*}\label{majoiterg}
0\leq \frac1{|\log\ell_n(x)|} \sum_{k=1}^{n-1} \log g_{u_k(x)}(|m_k(x)|) \leq \frac1{\log 3} \Big[\frac1{n}\sum_{k=1}^{n-1} \log \Big(\mfrac{6\la\sqrt{1+1/m_k^2(x)}-1}{6\la-1}\Big)\Big].
\end{align*}
Assuming that $x\in\cI$, i.e. $|m_n(x)|\to+\infty$, we can apply Cesàro's convergence lemma to show that the term inside brackets tends to $0$. With \eqref{itermnter}, this implies that \eqref{equivmn} is an equality and concludes the proof.
\end{proof}

Observe that the $\limsup$ in \eqref{hexpa} always belongs to $[0,1]$ when $\la\in[\frac13,\frac56]$ since
\begin{align*}
0\leq \beta_1(x,n) \log(6\la+1)+\beta_2(x,n)\log(6\la-1)
\leq \beta_{1,2}(x,n)\log 6.
\end{align*}

\addtocontents{toc}{\vspace{0.2cm}}%
\section{The set of points where $\fla$ has infinite derivative}\label{sec:derivinfini}

\noindent

In this section we are interested in the points $x$ at which $\fla$ has an infinite derivative, and we prove Theorem \ref{theo-derivinfini}.

\subsection{Proof of Theorem \ref{theo-derivinfini} $(i)$}\label{sec:proofTheoDerinf1}

\noindent

Assuming that $\la\in(\frac16,\frac12)$ is fixed, we prove in this subsection that there are points $x\in[0,1]$ at which $\fla$ has infinite derivatives and that $\dimh(\cV^\la_{\pm\infty})>0$.

\begin{proposition}\label{*lemmapi}
Let $\la\in(\frac16,\frac12)$ and $x\in[0,1]\setminus\cE$. There exists $M_\la>0$ such that if $|m_n(x)|> M_\la$ then
\begin{equation}\label{*eqlemmapib}
\fla(I_n(x))\subset\flan(I_n(x))=(\flan(a_n(x)),\flan(b_n(x))).
\end{equation}
\end{proposition}

To prove \eqref{*eqlemmapib} it is sufficient to prove that
if $|m_n(x)|> M_\la$ then
\begin{equation}\label{*eqlemmapi}
\fla(x)\in\flan(I_n(x))=(\flan(a_n(x)),\flan(b_n(x))).
\end{equation}
Indeed, for any $z\in I_n(x)$ we have $I_n(x)=I_n(z)$. Hence \eqref{*eqlemmapi} applied with $z$ instead of $x$ gives $\fla(z)\in\fla_n(I_n(z))=\fla_n(I_n(x))$.

Since $\la<\frac12$ we can select a sufficiently large $M'_\la>1$ for which
\begin{equation}\label{*mllv}
\frac12>\la\sqrt{1+\mfrac1{(M'_\la)^2}} = \frac12-\eps'_\la \text{ with } \eps'_\la>0.
\end{equation}
We claim that
\begin{equation}\label{*mllclaim}
F^\la_{n+1}(I_{n+1}(x))=\big(F^\la_{n+1}(a_{n+1}(x)),F^\la_{n+1}(b_{n+1}(x))\big)\subset \flan(I_n(x)).
\end{equation}

\medskip

We start with a lemma where only the breakpoints of $F_{n+1}^\la$ are considered.

\begin{lemma}\label{lem-inclusion}
If $|m_n(x)|\geq M'_\la$, then $F_{n+1}^\la(a_{n+1}(x))$ and $F_{n+1}^\la(b_{n+1}(x))$
both belong to the interval $[\fla(a_n(x)), \fla(b_n(x))]$.
\end{lemma}

This lemma means that when the slope of $\flan$ is large, the three breakpoints of $F^\la_{n+1}$ inside $I_n(x)$ have images that are located strictly inside the original interval $(\fla(a_n(x)),\fla(b_n(x)))$. This is the crucial difference between ``small'' and ``large'' $\la$, since this does not occur when $\la>\frac12$, as it can be deduced from the first two lines of the calculations in \eqref{*one12} below,  see also  Figures \ref{fig:NondiffB}, \ref{fig:noinfderiv1} and \ref{fig:noinfderiv2}.

\begin{proof}
We assume without loss of generality that $\flan$ is increasing on $I_n(x)$, so that $\fla(a_n(x))< \fla(b_n(x))$ and $m_n(x)>0$, the other case is symmetric.

The claim is trivial when $a_{n+1}(x)=a_n(x)$ or $b_{n+1}(x)=b_n(x)$. To conclude, it is enough to take care of the three breakpoints of $F^\la_{n+1}$ located inside $I_n(x)$ (recall Figure \ref{fig:NondiffB}).

The first breakpoint is $a_n(x)+\eps_n(x)\mfrac{\ell_n(x)}3$ for which one has
\begin{equation}\label{*one13}
F^\la_{n+1}\big(a_n(x)+\eps_n(x)\mfrac{\ell_n(x)}3\big) = \flan(a_n(x))+\mfrac13\ell_n(x)m_n(x) \in \flan(I_n(x)).
\end{equation}
Similarly, for the third breakpoint one has
\begin{equation}\label{*one23}
F^\la_{n+1}\big(a_n(x)+\eps_n(x)\mfrac{2\ell_n(x)}3\big) = \flan(a_n(x))+\mfrac23\ell_n(x)m_n(x) \in \flan(I_n(x)).
\end{equation}
The critical case is the second breakpoint located at the middle of $I_n(x)$. Since $m_n(x)\geq M'_\la$ by \eqref{*mllv} we obtain that
\begin{align}\label{*one12}
F^\la_{n+1}\big(a_n(x) +   \eps_n(x)\mfrac{\ell_n(x)}2\big)
& = \flan(a_n(x))+\mfrac13\ell_n(x)m_n(x) + \mfrac16\ell_n(x)m_{n+1}(x) \nonumber \\
& = \flan(a_n(x)) + \ell_n(x)m_n(x)\Big(\mfrac13 +  \mfrac16\Big( 1+6\la\sqrt{1+\mfrac1{m_n^2(x)}}\Big)\Big)  \nonumber \\
& \leq \flan(a_n(x)) +\ell_n(x)m_n(x)(1-\eps'_\la) \nonumber \\
& < \flan(a_n(x)) + \ell_n(x)m_n(x) = \flan(b_n(x)),
\end{align}
hence the result.
\end{proof}

We are now able to prove Proposition \ref{*lemmapi}.

\begin{proof}[Proof of Proposition \ref{*lemmapi}]

First assume that $\la\in[\frac13,\frac12)$, that is $6\la-1\geq 1$. Recall that by \eqref{*slom}, for every $n\geq0$ and $x\in[0,1]\setminus\cE$, either  $m_{n+1}(x)=m_n(x)$ or $|m_{n+1}(x)|\geq (6\la-1)|m_n(x)|$. In this case the sequence of slopes is non-decreasing. Hence,   $|m_{n+1}(x)|\geq |m_n(x)|\geq M'_\la$ and Lemma~\ref{lem-inclusion} can be used on $I_{n+1}(x)$ as well. Subsequently, if $\la\geq\frac13$ then the sets $\flan(I_n(x))$ are nested and
\begin{equation*}
\{\fla(x)\}= \bigcap_{k=0}^{\infty} F^\la_{n+k}(I_{n+k}(x))\subset \flan(I_n(x)).
\end{equation*}

\medskip

\noindent Assume now that $\la\in(\frac16,\frac13)$, that is $6\la-1\in(0,1)$. The situation is more intricate since the slopes $|m_{n+k}(x)|$ may decrease when $k$ increases.

Set
\begin{equation}\label{*eqMl}
M_\la = \frac{M'_\la}{(6\la-1)^l} >M'_\la>1
\end{equation}
with a sufficiently large integer $l$, to be determined later in \eqref{*lval}.

\noindent We start with three observations:
\begin{itemize}[parsep=-0.15cm,itemsep=0.25cm,topsep=0.2cm,wide=0.175cm,leftmargin=0.5cm]
\item[$\bullet$] When $u_n(x)\in\{1,2\}$ then $x\in(a_n(x)+\frac13\ell_n(x),a_n(x) + \frac23\ell_n(x))$. Therefore, recalling \eqref{*mllv},
\begin{align}\label{*A2*a}
\flan(b_n(x)) - F^\la_{n+1}(b_{n+1}(x))
& \geq \ell_n(x) \Big(|m_n(x)| - \Big( \mfrac{|m_n(x)|}{2} +\la\sqrt{1+m^2_n(x)} \Big) \Big) \nonumber \\
& = \ell_n(x)|m_n(x)|\Big( 1 - \frac12 - \la\sqrt{1+\mfrac1{m^2_n(x)}}  \Big) \nonumber\\
& \geq \ell_n(x)|m_n(x)|\Big( \frac12 - \la\sqrt{1+\mfrac1{(M'_\la)^2}}  \Big) = \ell_n(x)|m_n(x)|\eps'_\la.
\end{align}
\item[$\bullet$] When $u_n(x)\in\{0,1,3\}$ then $|m_{n+1}(x)|\geq |m_n(x)|\geq M'_\la$ and Lemma \ref{lem-inclusion} can be used on $I_{n+1}(x)$ as well.
\item[$\bullet$] When $u_n(x)=2$ then
\begin{equation}\label{*A2*b}
|m_{n+1}(x)|=|m_n(x)| \Big (6\la\sqrt{1+\mfrac1{m_n^2(x)}} - 1\Big) >|m_n(x)|(6\la-1)>0.
\end{equation}
In this case, $|m_{n+1}(x)|$ can be smaller than $|m_n(x)|$ since $6\la-1\in(0,1)$, this happens especially when $|m_n(x)|$ is very large, since in this situation $|m_{n+1}(x)|\approx (6\la-1)|m_n(x)|$.
\end{itemize}

\noindent Suppose that $n$ is such that $|m_n(x)|> M_\la$.

If $u_{n+k}(x)\in\{0,1,3\}$ for all $k\geq0$, then Lemma \ref{lem-inclusion} applies to every $I_{n+k}(x)$ and the sets $F^\la_{n+k}(I_{n+k}(x))$ are nested. Hence \eqref{*eqlemmapi}, and then \eqref{*eqlemmapib}, hold true for $x$.

Assume there exists $k_1$ such that $u_{n+k}(x)\in\{0,1,3\}$ for $0\leq k<k_1$, but $u_{n+k_1}(x)=2$. The idea is that $M_\la$ is so large that even if the slopes decrease after $k_1$ and that Lemma \ref{lem-inclusion} does not apply anymore, the slopes $|m_{n+k}(x)|$ will remain larger than $M'_\la$ for sufficiently many integers $k$ to ensure that the oscillations of $\fla$ will stay controlled on $I_{n+k_1}(x)$.

Observe first that by definition of $k_1$,
\begin{equation}\label{*nestko}
|m_{n+k_1}(x)| \geq |m_n(x) |\geq M_\la \text{ and }
F^\la_{n+k_1}(I_{n+k_1}(x))\subset \flan(I_n(x)).
\end{equation}
By \eqref{*A2*b}, since $u_{n+k_1}(x)=2$,
\begin{equation}\label{*mnlest}
|m_{n+k_1+1}(x)| > |m_{n+k_1}(x)| (6\la-1)\geq M_\la(6\la-1) > M'_\la.
\end{equation}
Moreover, by \eqref{*A2*a} and \ref{*nestko}  recalling that $L_n(x)=\ell_n(x)\sqrt{1+m^2_n(x)}<2 \ell_n(x)|m_n(x)|$ if $|m_n(x)|>M_\la>1$, one has
\begin{align}\label{*fnknfnbn}
F^\la_{n+k_1+1}(x) \leq F^\la_{n+k_1+1}(b_{n+k_1+1}(x))
& \leq F^\la_{n+k_1}(b_{n+k_1}(x)) -\eps'_\la\ell_{n+k_1}(x)|m_{n+k_1}(x)| \nonumber \\
& \leq \flan(b_n(x)) -\mfrac12\eps'_\la L_{n+k_1}(x).
\end{align}
This implies that
\begin{equation}\label{eq-iteration1}
F^\la_{n+k_1+1}(I_{n+k_1+1}(x))\subset F^\la_{n+k_1}(I_{n+k_1}(x))\subset \flan(I_n(x)).
\end{equation}

\noindent Repeating \eqref{*mnlest} and \eqref{eq-iteration1}
as long as $|m_{n+k_1+k}(x)|\geq M'_\la$, the same argument gives
\begin{equation}\label{eq-inclu2}
F^\la_{n+k_1+k+1}(I_{n+k_1+k+1}(x))\subset F^\la_{n+k_1+k}(I_{n+k_1+k}(x))\subset F^\la_{n+k_1}(I_{n+k_1}(x)).
\end{equation}

Let us now look at the slopes $|m_{n+k_1+k}(x)|$. Recall that $(6\la-1)\in(0,1)$. Hence
\begin{equation*}\label{27a}
\text{$|m_{n+k_1+k+1}(x)|>|m_{n+k_1+k}(x)|>(6\la-1)|m_{n+k_1+k}(x)|$ when $u_{n+k_1+k}(x)\in\{0,1,3\}$}.
\end{equation*}
If $u_{n+k_1+k}(x)=2$, we can use \eqref{*mnlest}. Therefore, for $k\leq l-1$, we have
\begin{equation}
|m_{n+k_1+k+1}(x) | > |m_{n+k_1+1}(x)| (6\la-1)^k\geq M_\la(6\la-1)^{k+1} \geq M'_\la.
\end{equation}
By \eqref{*eqMl}, this applies in particular to $k=l-1$, and $F^\la_{n+k_1+l}(I_{n+k_1+l}(x))\subset F^\la_{n+k_1}(I_{n+k_1}(x))$.

We will prove now that the oscillations of $\fla$ on $I_{n+k_1+l}(x)$ are so small that they force $\fla(x)$ to belong to the initial interval $(\flan(a_n(x)),\flan(b_n(x)))$. Using \eqref{boundLn} and \eqref{majincFn} with $n+k_1$ instead of $n$ and with $z=x$, we obtain, for every $k\geq0$,
\begin{equation}\label{*Fnke}
|F^\la_{n+k_1+l+k+1}(x)-F^\la_{n+k_1+l+k}(x)|\leq \la L_{n+k_1+k+l}(x) \leq \la\Big(\mfrac16+\la\Big)^{k+l}L_{n+k_1}(x).
\end{equation}
Since $F^\la_{n+k_1+l+k}(x)\to\fla(x) $ monotone increasingly as $k\to+\infty$, we obtain
\begin{equation}
|\fla(x)- F^\la_{n+k_1+l}(x)| \leq \la\Big(\mfrac16+\la\Big)^l L_{n+k_1}(x) \sum_{k=0}^{\infty}\Big(\mfrac16+\la\Big)^k<\mfrac12\eps'_\la L_{n+k_1}(x),
\end{equation}
where $l$ is chosen so large that the last inequality holds, that is
\begin{equation}\label{*lval}
\la\Big(\mfrac16+\la\Big)^l\Big(\mfrac6{5-6\la}\Big)< \mfrac12\eps'_\la.
\end{equation}
From \eqref{*fnknfnbn} it finally follows that $\fla(x)\leq \flan(b_n(x))$, which implies \eqref{*eqlemmapi}.
\end{proof}

\begin{figure}[h!]
\centering
\includegraphics[width=0.5\textwidth]{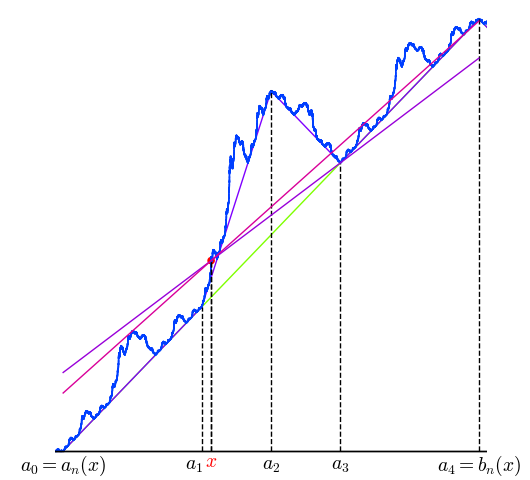}
\captionsetup{width=1\textwidth}
\caption{Infinite derivative of $\fla$ at a point $x$ when $\la\in(\frac16,\frac12)$.}
\label{fig:noinfderiv1}
\end{figure}

We can now prove item $(i)$ of Theorem \ref{theo-derivinfini}. Our strategy is to identify subsets $\cV_{\pm}\subset \cV^\la_{\pm\infty}$ whose Hausdorff dimension may be computed. Precisely, set
\begin{align*}
\cV_+ & = \{x\in [0,1] :  (u_{3n}(x),u_{3n+1}(x),u_{3n+2}(x))\in \{(0,1,0), (0,0,1)\} \text{ for all $n\geq0$}\} \\
\mbox{ and }
\cV_- & = \{x\in [0,1] : (u_0(x),u_1(x),u_2(x))=(2,0,0) \text{ and } \\
& \hspace{25mm} (u_{3n}(x),u_{3n+1}(x),u_{3n+2}(x))\in\{(0,1,0),(0,0,1)\} \text{ for all $n\geq1$}\}.
\end{align*}

\begin{proposition}\label{propderivinfini}
Let $\la\in(\frac16,\frac12)$. Then, $\fla$ has a positive (resp. negative) infinite derivative at every point $x\in\cV_+$ (resp.  $x\in\cV_-$) and $\dimh(\cV_+)=\dimh(\cV_-) =\mfrac{\log 2}{\log 54}>0$.
\end{proposition}

\begin{proof} Since the two cases are similar we detail only the case of $\cV_+$. It is clear that $m_n(x)\to+\infty$ for $x\in  \cV_+$. Hence we can suppose that $n$ is so large that $m_n(x) \geq M_\la$ (recall Proposition \ref{*lemmapi}). Observe also that for every $n$,  $a_n(x) < b_n(x) $.

Suppose that $x,y\in I_n(x)=(a_n(x),b_n(x))$ but
$y\not\in I_{n+1}(x)=(a_{n+1}(x),b_{n+1}(x))$, in particular $x\neq y$. We bound from below the difference quotient $(\fla(x)-\fla(y))/(x-y)$ to prove that $\fla$ has a positive infinite derivative at $x$.

\medskip

\noindent{\it Case 1: $y<x$}. Recall that the digits of $x$ are $(u_0(x),u_1(x),\ldots)$.  Since $y<x$ it is impossible to have $u_n(x)=0$, so necessarily $u_n(x)=1$ and $u_n(y)=0$. One has then $a_{n+1}(y)=a_n(x)=a_n(y)$ and $b_{n+1}(y)=a_{n+1}(x) = a_n(x) + \frac13\ell_n(x)$. By Proposition \ref{*lemmapi},
\begin{equation*}
\fla(y)\leq \flan(a_{n+1}(x) )=\fla(a_{n+1}(x)).
\end{equation*}
From the definition of $\cV_+$, necessarily $u_{n+1}(x)=0$ and either $u_{n+2}(x)=1$, or ($u_{n+2}(x)=0$ and $u_{n+3}(x)=1$), or ($u_{n+2}(x)=u_{n+3}(x)=0$ and $u_{n+4}(x)=1$). In all situations,
\begin{equation*}
\fla(x)\geq \flan(x)\geq \fla(a_{n+1}(x))+m_n(x)(x-a_{n+1}(x))\geq \fla(a_{n+1}(x)) +\frac{\ell_n(x)m_n(x)}{6 \cdot 3  \cdot 3 \cdot 3}.
\end{equation*}
Thus,
\begin{equation}\label{*eqdiffqa}
\frac{\fla(x) - \fla(y)}{x-y}\geq  \frac{\ell_n(x) m_n(x)}{162\ell_n(x)}= \frac1{162} m_n(x).
\end{equation}

\medskip

\noindent{\it Case 2: $y>x$}. We consider several subcases according to the digits of $x$ and $y$.

\smallskip
{(a)} If $u_n(x)=0$, then $(a_{n+1}(x), b_{n+1}(x))=(a_n(x), a_n(x)+ \frac13\ell_n(x))$. Moreover, $a_{n+1}(y)\geq b_{n+1}(x)$. Hence,
\begin{equation}\label{*becsl}
\fla(y)\geq F^\la_{n+1}(y) \geq F^\la_{n+1}(a_{n+1}(y)) \geq F^\la_{n+1}( b_{n+1}(x) ) = \flan(a_{n}(x)) + \mfrac13\ell_n(x)m_n(x).
\end{equation}

\smallskip
{(a1)} If $u_{n+1}(x) =0$, then $(a_{n+2}(x), b_{n+2}(x) )=(a_n(x),a_n(x) +\frac19\ell_n(x))$, and Proposition  \ref{*lemmapi} gives
\begin{equation*}
\fla(x)\leq  F^\la_{n+1}( b_{n+1}(x))= \flan( b_{n+2}(x) )= \flan( a_n(x)) +\mfrac19\ell_n(x)m_n(x).
\end{equation*}
This yields the estimate
\begin{equation}\label{*eqdiffqb}
\frac{\fla(y)- \fla(x)}{y-x}\geq \frac{\ell_n(x) m_n(x)(\frac13 - \frac19)}{\ell_n(x)}= \frac29 m_n(x).
\end{equation}

\smallskip {(a2)} If $u_{n+1}(x)=1$, then $(a_{n+2}(x),b_{n+2}(x)) = (a_n(x)+\frac19\ell_n(x),a_n(x) + \frac19\ell_n(x)+\frac1{18}\ell_n(x))$. Using that$\la<\frac12$, we can select a small $\eps_\la>0$ such that for sufficiently large $m_n(x)$ we have
\begin{equation*}
m_{n+1}(x) =m_n(x)\Big(1+6\la\sqrt{1 +\mfrac1{m_n^2(x)}}\Big)<m_n(x)(4-\eps_\la).
\end{equation*}
By Proposition \ref{*lemmapi}, one deduces that
\begin{align*}
\fla(x)\leq F^\la_{n+2}(b_{n+2}(x) )
& = F^\la_{n+2}(a_{n+2}(x)) + \mfrac1{18}\ell_n(x)m_{n+1}(x) \nonumber \\
& = \flan(a_n(x)) + \mfrac19\ell_n(x)m_n(x) +  \mfrac1{18}\ell_n(x)m_{n+1}(x) \nonumber \\
& < \flan(a_n(x))+ \ell_n(x)m_n(x)\Big(\mfrac19+\mfrac{4-\eps_\la}{18}\Big).
\end{align*}
Now using \eqref{*becsl}, we obtain for the difference quotient that
\begin{equation}\label{*eqdiffqc}
\frac{ \fla(y) - \fla(x)}{y-x}\geq \frac{\ell_n(x) m_n(x) (\frac13 - \frac19- \frac{4-\eps_\la}{18})}{\ell_n(x)}= \frac{\eps_\la}{18}m_n(x).
\end{equation}

\medskip

(b) If $u_n(x)=1$, then $(a_{n+1}(x), b_{n+1}(x) )=(a_n(x) + \frac13\ell_n(x),a_n(x) + \frac13\ell_n(x)+\frac16\ell_n(x))$. One necessarily has $u_{n+1}(x)=0$ and again by Proposition \ref{*lemmapi},
\begin{align}\label{*alfdiff}
\fla(x) & \leq F^\la_{n+2}(b_{n+2}(x)) \leq F^\la_{n+1}(b_{n+1}(x))=
\flan(a_n(x))+\mfrac13\ell_n(x)m_n(x)+\mfrac1{18}\ell_n(x)m_{n+1}(x) \nonumber \\
& < \flan(a_n(x))+\ell_n(x)m_n(x)\Big (\mfrac13+\mfrac{4-\eps_\la}{18}\Big) < \flan(a_n(x))+\mfrac59\ell_n(x)m_n(x).
\end{align}
Since $u_n(y)\neq u_n(x)$, $u_n(y)\in\{2,3\}.$ Using that $F^\la_{n+1}$ is monotone decreasing on
$(a_n(x)+\frac12\ell_n(x),a_n(x) + \frac23\ell_n(x))$ and monotone increasing on $(a_n(x) + \frac23\ell_n(x),a_n(x)+\ell_n(x))$, one has
\begin{equation*}
\fla(y)\geq F^\la_{n+1}(y) \geq F^\la_{n+1}(a_n(x))+ \mfrac23\ell_n(x)m_n(x) = \flan(a_n(x))+ \mfrac23\ell_n(x)m_n(x).
\end{equation*}
Combining this with \eqref{*alfdiff} gives
\begin{equation}\label{*eqdiffqd}
\frac{\fla(y)-\fla(x)}{y-x}\geq \frac19 m_n(x).
\end{equation}
Since $m_n(x)\to +\infty$, the estimates \eqref{*eqdiffqa}, \eqref{*eqdiffqb}, \eqref{*eqdiffqc} and \eqref{*eqdiffqd} show  that $(\fla)'(x)=+\infty$.

\medskip

Finally, observe that $\cV_+$ can be regarded to be the attractor of an IFS generated by two similarities which map $[0,1]$ onto two disjoint subintervals of length $\frac{1}{6 \cdot 3 \cdot 3}=\frac1{54}$ and hence its similarity and Hausdorff dimension both equal $\frac{\log2}{\log 54}$.
\end{proof}

\subsection{Proof of Theorem \ref{theo-derivinfini} $(ii)$}\label{sec:proofTheoDerinf2}

\noindent

\begin{figure}[ht!]
\centering
\includegraphics[width=0.5\textwidth]{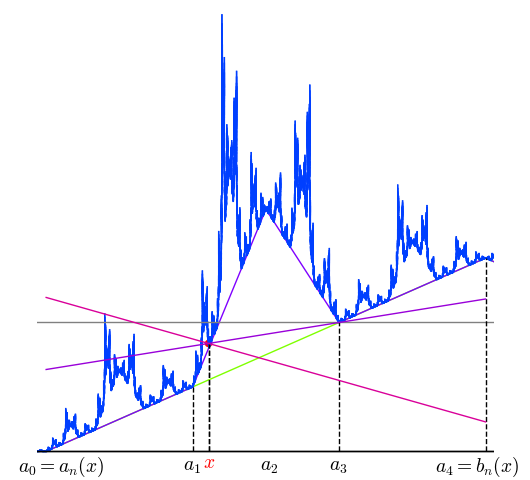}
\captionsetup{width=1\textwidth}
\caption{Non existence of infinite derivative of $\fla$ at point $x$ due to high oscillations when $\la\in(\frac23,\frac56)$.}
\label{fig:noinfderiv2}
\end{figure}

Assuming that $\la\in(\frac23,\frac56)$ is fixed, we will prove in this subsection that there are no points $x\in[0,1]$ such that $(\fla)'(x)=+\infty$ or $(\fla)'(x)=-\infty$. The main difference with the previous section, and the idea of the proof, is that   $\la$ is so large that on every $I_n(x)$, $\fla$ oscillates too much to allow the existence of a limit (even infinite) derivative, see Figure \ref{fig:noinfderiv2}.

We again only deal with the case $(\fla)'(x)=+\infty$, the other case is similar and is left to the reader. Proceeding towards a contradiction, assume that $(\fla)'(x)=+\infty$ for a certain $x\in[0,1]$. By the proof of Proposition \ref{propexphmax}, $x\not\in\cEt$, and by Lemma \ref{lemslopbound}, $|m_n(x)|\to+\infty$.

\medskip

Let us fix $n\geq0$ large enough such that $m_n(x)>0$ and
\begin{equation}\label{*deriv}
\frac{\fla(y)-\fla(x)}{y-x}>1 \text{ for all $y\in [a_n(x),b_n(x)]$ and $y\neq x$.}
\end{equation}
Consider the point $z\in I_n(x)=(a_n(x),b_n(x))$ for which $u_i(z)= u_i(x)$ for  all $i<n$ and $u_i(z)=1$ for $i\geq n$. In particular,
\begin{equation}\label{debutz}
a_n(z)=a_n(x),\quad\ell_n(z)=\ell_n(x),\quad m_n(z)=m_n(x)
\end{equation}
and
\begin{equation}
a_n(x)< z < \frac{a_n(x) + b_n(x)}{2} = a_n(x) + \mfrac12\ell_n(x).
\end{equation}
For all $k\geq0$, we have
$\ell_{n+k+1}(z)= \frac16 \ell_{n+k}(z)$, thus
\begin{equation}\label{lnk}
\ell_{n+k}(z)= \frac1{6^k} \ell_n(z).
\end{equation}
Moreover, \eqref{itermn} implies
\begin{equation*}\label{}
m_{n+k+1}(z)= m_{n+k}(z)\Big(1+6\la\sqrt{1+\mfrac1{m_{n+k}^2(z) }}\Big)>(1+6\la)m_{n+k}(z),
\end{equation*}
yielding
\begin{equation}\label{minormnk}
m_{n+k}(z)\geq (1+6\la)^k m_n(z), \text{ the inequality being strict for $k\geq1$.}
\end{equation}

\medskip

Observe now that
\begin{equation*}\label{}
F^\la_{n+k+1}(a_{n+k+1}(z)) = F^\la_{n+k}(a_{n+k}(z)+\tfrac13\ell_{n+k}(z)) = F^\la_{n+k}(a_{n+k}(z))+\mfrac13\ell_{n+k}(z)m_{n+k}(z).
\end{equation*}
It follows by induction that
\begin{equation*}\label{}
F^\la_{n+k}(a_{n+k}(z)) = \flan(a_n(z)) + \frac13\sum_{i=0}^{k-1} \ell_{n+i}(z)m_{n+i}(z).
\end{equation*}
Hence, using \eqref{debutz}, \eqref{lnk} and \eqref{minormnk},
\begin{align}\label{*maxestn}
F^\la_{n+k}(a_{n+k}(z)) & \geq \flan(a_n(x)) + \frac13\ell_n(x)m_n(x)\sum_{i=0}^{k-1} \left(\frac{1+6\la}{6}\right)^i \nonumber \\
& = \flan(a_n(x)) + \ell_n(x)m_n(x)\left(\frac{2}{5-6\la} \right)\Big(1-\left(\frac{1+6\la}{6}\right)^k\Big).
\end{align}
Since $\fla(z)>F^\la_{n+k}(a_{n+k}(z))$ and $\la\in(\frac23,\frac56)$, taking the limit $k\to\infty$ in \eqref{*maxestn} yields
\begin{align}\label{*maxest}
\fla(z)> \flan(a_n(x)) + 2\ell_n(x)m_n(x) > \flan(a_n(x)) + \ell_n(x)m_n(x) = \flan(b_n(x)).
\end{align}
This means that the value of $\fla$ at $z\in I_n(x)$ is larger than $\fla(b_n(x))$. By Proposition \ref{*lemmapi} this is clearly impossible when $\la<\frac12$ and $|m_n(x)|$ is sufficiently large. This is a drastic change in the geometry of our function. It is illustrated by Figure \ref{fig:noinfderiv2}, where one sees that images of some points in the first subinterval of $I_n(x)$ have a larger value than $\fla(a_n(x)+\frac13\ell_n(x))$.

\medskip

Next, applying \eqref{*deriv} twice, \eqref{*maxest} gives
\begin{equation}\label{comparzbx}
\fla(z)>\flan(b_n(x))=\fla(b_n(x))>\fla(x),
\end{equation}
so that necessarily $x$ must be smaller than $z$. Recalling the definition of $z$, i.e. $u_i(z)= u_i(x)$ for all $i<n$ and $u_i(z)=1$ for $i\geq n$, this implies that $x<a_n(x)+\frac12\ell_n(x)$, and $u_n(x)\in\{0,1\}$. Also, since $m_n(x)\to +\infty$, one must have $u_n(x)=1$ infinitely often. Hence it is enough to deal with the case $u_n(x)=1$ since this occurs for infinitely many $n$s.

Assume then that $u_n(x)=1$. Therefore, using that $a_n(x)+\mfrac23\ell_n(x)$ is a breakpoint of $F^\la_{n+1}$ and $x<a_n(x)+\mfrac23\ell_n(x)$, \eqref{*deriv} gives
\begin{align}\label{*fllln}
\flan(x)\leq \fla(x)\leq \fla(a_n(x)+\mfrac23\ell_n(x))
& = F^\la_{n+1}(a_n(x)+\mfrac23\ell_n(x)) \nonumber \\
& = \flan(a_n(x))+\mfrac23\ell_n(x)m_n(x).
\end{align}

On the other hand, consider an arbitrary point $x'\in (a_n(x),a_n(x)+\frac13\ell_n(x))= I_{n+1}(x')$. Since $u_n(x)=1$, $x'<x$, $m_{n+1}(x')=m_n(x)$ and $\ell_{n+1}(x')=\frac13\ell_n(x)$. By the argument above used to get \eqref{*maxest}, there exists $z'\in I_{n+1}(x')$ for which $x'<z'< a_n(x)+\frac13\ell_n(x) <x$ and
\begin{equation*}\label{*maxestbis}
\fla(z') > F^\la_{n+1}(a_{n+1}(x'))+ 2 \ell_{n+1}(x') m_{n+1}(x') =
\flan(a_n(x))+ \mfrac23 \ell_n(x) m_n(x).
\end{equation*}
This inequality combined with \eqref{*fllln} gives that $\fla(z')>\fla(x)$, which by using \eqref{*deriv} implies $z'>x$, hence a contradiction.

\addtocontents{toc}{\vspace{0.2cm}}%
\section{The measure $\mula$ and its multifractal spectrum}\label{sec:spectrum}

In this section, we investigate the multifractal spectrum of $\mula$, which paves the way for the proof of Theorem \ref{spectruma}.

\subsection{Recalls on the multifractal properties of $\mula$}\label{sec:mula}

\noindent

Let us begin by recalling key notions about multifractal analysis of measures. The lower and upper local dimensions of a probability measure $\mu$ defined on $\R$ are
\begin{equation}\label{defdimloc}
\ldimloc(\mu,x)=\liminf_{r\to 0}\frac{\log \mu(B(x,r))}{\log r} \text{ and } \udimloc(\mu,x)=\limsup_{r\to 0}\frac{\log \mu(B(x,r))}{\log r},
\end{equation}
and, when these quantities agree, their common value $\dimloc(\mu,x)$ is called the local dimension of $\mu$ at the point $x$. In addition, the level sets associated with these local dimensions are
\begin{equation*}\label{defloclevelset}
\underline{E}_\mu(\alpha) = \{x\in\R: \ldimloc(\mu,x)=\alpha\},\,\overline{E}_\mu(\alpha) = \{x \in\R: \udimloc(\mu,x)=\alpha\}
\end{equation*}
and
\begin{equation*}\label{defloclevelsetbis}
E_\mu(\alpha) = \{x\in\R: \dimloc(\mu,x)=\alpha\}.
\end{equation*}
We are especially interested in the multifractal spectrum of $\mu$, that is in the mapping $d_\mu:\R^+\to [0,1]\cup\{-\infty\}$ defined by
\begin{equation}\label{defspecmu}
d_\mu(\alpha)= \dimh(\underline{E}_\mu(\alpha)) = \dimh( \{x\in\R: \ldimloc(\mu,x)=\alpha\}).
\end{equation}
Determination of $d_\mu(\alpha)$ requires to consider the following set
\begin{equation}\label{def-Eleq}
\underline{E}_\mu^\geq (\alpha) = \{x\in\R: \ldimloc(\mu,x)\geq \alpha\} \supset \underline{E}_\mu(\alpha).
\end{equation}

\medskip

Let us return to the measure $\mula$ from Definition \ref{def:mula}. Since it is associated with a hyperbolic IFS satisfying the Open Set Condition, $\mula$ is known to enjoy many multifractal properties that will be of great importance in the following and are listed at the end of the section. Before, let us observe that this measure can be alternatively constructed
by using the Mass Distribution Principle and considering the successive intervals of generation $n$ from $\cI_n$. Precisely, starting with $\mula((0,1))=1$ then set $\mula(S_i((0,1)))= p_{\la,i}$ for $i=0,\ldots,3$. In general, assuming that $\mula(I_n(x_0))$ is already defined for an interval $I_n(x_0)=(a_n(x_0),a_n(x_0)+\eps_n(x_0)\ell_n(x_0))$ with $x_0\in[0,1]\setminus\cE$, then set
\begin{equation}\label{eq:itermu}
\begin{cases}
& \mula\big((a_n(x_0),a_n(x_0)+\tfrac13\eps_n(x_0)\ell_n(x_0))\big) = p_{\la,0} \mula(I_n(x_0)), \\
& \mula\big((a_n(x_0)+\tfrac13\eps_n(x_0)\ell_n(x_0),a_n(x_0)+\tfrac12\eps_n(x_0)\ell_n(x_0))\big) = p_{\la,1} \mula(I_n(x_0)), \\
& \mula\big((a_n(x_0)+\tfrac12\eps_n(x_0)\ell_n(x_0),a_n(x_0)+\tfrac23\eps_n(x_0)\ell_n(x_0))\big) = p_{\la,2} \mula(I_n(x_0)), \\
& \mula\big((a_n(x_0)+\tfrac23\eps_n(x_0)\ell_n(x_0),a_n(x_0)+\eps_n(x_0)\ell_n(x_0))\big) = p_{\la,3} \mula(I_n(x_0)).
\end{cases}
\end{equation}
Since the family $\{\cI_n : n\geq0\}$ generates the Borel $\sigma$-algebra in $[0,1]$, the iterative construction above defines a Borel measure supported in $[0,1]$ satisfying \eqref{defmula}. In particular, it follows from this construction that, for every $x\in[0,1]\setminus\cE$,
\begin{equation}\label{eq:muIn}
\mula(I_n(x)) = \prod_{i=0}^3 p_{\la,i}^{\beta_i(x,n)}.
\end{equation}
Therefore, it will be convenient to deal with quantities $\mula(I_n(x))$ instead of $\mula(B(x,r))$ to compute the local dimensions in \eqref{defdimloc}. The following property of $\mula$ makes this possible.

\begin{lemma}\label{lem-doubling}
Let $\la\in(\frac16,+\infty)$. Then $\mula$ is a doubling measure, i.e. there exists a constant $C_\la\geq 1$ such that, for every $x\in[0,1]$ and every $r>0$, $\mula(B(x,2r))\leq C_\la\mula(B(x,r))$.
\end{lemma}

\begin{proof}
This will follow from $p_{\la,0}=p_{\la,3}$ and from Lemma \ref{lem-neighbors}.
The cases $x=0$, or $x=1$ can be obtained from the general case $x\in(0,1)$ and we leave details to the reader. Hence we suppose that $x\in (0,1)$. Without limiting generality suppose that $0<2r<\frac16$ and select $n\geq1$ such that
\begin{equation}\label{*eqvaln}
\ell_n(x)\leq 3r < \ell_{n-1}(x).
\end{equation}
We divide the interval $I_n(x)= (a_n(x),b_n(x))$ into three subintervals $I$, $I'$ and $I''$ of equal length $\frac13|I_n(x)|= \frac13\ell_n(x)$. From the definition of $\mula$,
\begin{equation*}
\min(\mula(I),\mula(I'),\mula(I''))\geq \Big(\min_{0\leq i \leq 3} p_{\la,i}\Big)\mula(I_n(x)).
\end{equation*}
Since $B(x,r)$ necessarily contains one of the three subintervals $I$, $I'$ or $I''$, we deduce that $\mula(B(x,r))\geq C \mula(I_n(x))$ with some positive constant $C>0$.

Suppose that $x\leq \frac12(a_{n-1}(x)+b_{n-1}(x))$, the other case is similar. Choose $y$ such that $I_{n-1}(y)$ and $I_{n-1}(x)$ are contiguous and $I_{n-1}(y)$ is to the left of $I_{n-1}(x)$. By Lemma \ref{lem-neighbors} and \eqref{*eqvaln}, we have $B(x,\frac32r)\subset I_{n-1}(x)\cup I_{n-1}(y)$. From the definition of $\mula$ it follows that $\mula(I_{n-1}(x)) \leq C'\mula(I_n(x))$, where $C'=\max_{0\leq i \leq3} p_{\la,i}^{-1}>0$.

Let $m$ be the largest integer such that $I_{n-1}(x)$ and $I_{n-1}(y)$ are included in the same interval of generation $m\geq 0$. Observe that $m\leq n-2$. This implies that $u_p(x)=u_p(y)$ for every $0\leq p\leq m-1$. The intervals $I_{n-1}(x)$ and $I_{n-1}(y)$ being contiguous, the only possibility for the values of $u_p(x)$ and $u_p(y)$ are $u_m(x)=u_m(y)+1$, and $u_p(x)=0$ and $u_p(y)=3$ for every $p\in \{m+1,\ldots,n-2\}$. This shows that $\mula(I_{n-1}(x))=\mula(I_m(x))p_{\la,u_m(x)} p_{\la,0}^{n-m-2}$ and $\mula(I_{n-1}(y))=\mula(I_m(y))p_{\la,u_m(y)} p_{\la,3}^{n-m-2}$. Since $p_{\la,0}=p_{\la,3}$ and $I_m(x)=I_m(y)$ we obtain $\mula(I_{n-1}(y))\leq C'' \mula (I_{n-1}(x))$ with a suitable constant $C''>0$.
Therefore,
\begin{align*}
\mula(B(x,\tfrac32r)) \leq \mula(I_{n-1}(x))+\mula(I_{n-1}(y))
& \leq (1+C'')\mula(I_{n-1}(x)) \nonumber \\
& \leq C'(1+C'')\mula(I_n(x))  \\
& \leq C'(1+C'') C^{-1}\mula(B(x,r)). \nonumber
\end{align*}
Since this holds for any $r>0$,
\begin{align}\label{eqdoubll}
\mula(B(x,2r))\leq \mula(B(x,(\tfrac32)^2r))
& \leq C'(1+C'') C^{-1}\mula(B(x,\tfrac32r)) \nonumber \\
& \leq (C'(1+C'') C^{-1})^2\mula(B(x,r))
\end{align}
which proves the claim by taking $C_\la = (C'(1+C'') C^{-1})^2$.
\end{proof}

From Lemma \ref{lem-doubling}, especially from \eqref{eqdoubll}, it follows that for every $x\in [0,1]$, in order to compute the local dimensions, one can replace the centered ball $B(x,r)$ by $I_n(x)$, i.e.
\begin{equation}\label{doubling}
\ldimloc(\mula,x)=\liminf_{n\to\infty }\frac{\log \mula(I_n(x))}{\log |I_n(x)|},
\end{equation}
and similarly for $\udimloc(\mula,x)$ and $\dimloc(\mula,x)$. This will simplify many computations in the following.

\begin{proposition}\label{propvaluedimloc}
Let $\la\in(\frac16,\frac56)$. For every $x\in[0,1]$, we have
\begin{equation}\label{value-dimloc}
\ldimloc(\mula,x) = \ga -\limsup_{n\to \infty}\frac{\beta_1(x,n)\log(6\la+1)+\beta_2(x,n)\log(6\la-1)}{\beta_{0,3}(x,n)\log3 +\beta_{1,2}(x,n)\log 6}.
\end{equation}
\end{proposition}

\begin{proof}
Using \eqref{xint} and \eqref{lnxdef}, we obtain $|I_n(x)|=\ell_n(x) = 3^{-\beta_{0,3}(x,n)}6^{-\beta_{1,2}(x,n)}$. In addition
by \eqref{def-pila} and  \eqref{eq:muIn}, we obtain
\begin{align*}
& \hspace{-0.5cm} \frac{\log \mula(I_n(x))}{\log |I_n(x)|}
= \frac{\beta_0(x,n) \log p_{\la,0}+\beta_1(x,n) \log p_{\la,1}+\beta_2(x,n) \log p_{\la,2}+\beta_3(x,n) \log p_{\la,3}}{ -\beta_{0,3}(x,n)\log 3- \beta_{1,2}(x,n) \log 6} \\
& =  \frac{-\beta_{0,3}(x,n)\ga\log 3 +\beta_1(x,n) \log (6\la+1)  + \beta_2(x,n) \log (6\la-1) -\beta_{1,2}(x,n)\ga\log 6}{-\beta_{0,3}(x,n)\log 3-  \beta_{1,2}(x,n)\log 6} \\
& = \ga - \frac{\beta_1(x,n)\log(6\la+1)+\beta_2(x,n)\log(6\la-1)}{\beta_{0,3}(x,n)\log3 +\beta_{1,2}(x,n)\log 6} .
\end{align*}
The result is then a consequence of \eqref{doubling}.
\end{proof}

Looking at \eqref{value-dimloc}, we introduce the notion of digit frequencies of $x$.

\begin{definition}\label{*defregp}
A point $x\in [0,1]$ is called {\em regular} if, for every $i=0,\ldots,3$, the frequency $\beta_i(x)$ of the digit $i$ exists, i.e.
\begin{equation}\label{*F1*a}
\lim_{n\to\infty}\frac{\beta_i(x,n)}{n}=\beta_i(x).
\end{equation}
We denote by $\cR$ the set of regular points. In particular, if $x\in\cR$, then necessarily $(\beta_i(x))_{i=0,\ldots,3}$ is a probability vector. Observe that any $x\in \cE$ is regular.
\end{definition}

Let us now turn to the $L^q$-spectrum $\tau_\mula$ of the measure $\mula$. We have chosen with Definition~\ref{def:Lqspectrum} an implicit but tractable definition of $\tau_\mula(q)$. Let us mention here that a ``direct'' definition of $\tau_\mula$ in terms of the $q$-moments of $\mula$ is the following, see e.g. \cite{BRMICHPEY}:
\begin{equation}\label{deftaubis}
\tau_\mu: q\in \R\mapsto \liminf_{j\to\infty} \frac{-1}{\log 3^j} \log\sum_{k=0}^{3^j-1} \big(\mu([k3^{-j},(k+1)3^{-j}])\big)^q \in \R\cup\{-\infty,+\infty\},
\end{equation}
where the summation runs only through indices such that $\mu([k3^{-j},(k+1)3^{-j}])\neq0$. The equality between \eqref{deftaubis} and the solution to \eqref{eq:Lqspectrum} in Definition~\ref{def:Lqspectrum} is proved for self-similar measures that satisfy the Open Set Condition in Corollary C of \cite{LauNgai99}.

\medskip

The properties of $\mula$ are gathered in the two following propositions. They follow from Theorem $1$ of \cite{BRMICHPEY} for the upper bounds for all dimensions, Theorems $2$ and $3$ of \cite{BRMICHPEY} for the lower bounds since $\mula$ is doubling, or Theorems A and B, Corollary C of \cite{LauNgai99} which treats more general cases. The reader may also check \cite{Olsen} for other notions of dimensions and formalisms, and \cite{LauFeng,Shmerkin2,BaFe21,Shmerkin1,Fen23} for recent spectacular advances on the multifractal formalism for self-similar measures with overlaps (our contractions do not overlap).

\medskip

Proposition \ref{recalls-mula1} below provides first results on  $\tau_\mula$ and $\tau^*_\mula$. We refer to \cite{CaMau92,falconer2003} and Sections $2$ and $3$ of \cite{LauNgai99} for the general expressions for self-similar measures. The explicit values for the measure $\mula$ in consideration follow directly from basic computations based on the implicit equation \eqref{eq:Lqspectrum} (by letting $q$ tend to infinity), see Section $1$ of \cite{CaMau92}.

\begin{proposition}\label{recalls-mula1}
Let $\la\in(\frac16,\frac56)$ and consider $\mula$, $\tau_\mula$ and $\tau^*_\mula$. Denote by $[\alphat_{\la,\min},\alphat_{\la,\max}]$ the support of $\tau_\mula^*$. Then, the mapping $\alpha\mapsto \tau^*_\mula(\alpha)$ is real analytic concave on the interior of its support, and we have
\begin{align*}
\alphat_{\la,\min} & = \lim_{q\to+\infty} \tau_\mula'(q) =\min _{0\leq i \leq 3}\mfrac{\log p_{\la,i}}{\log r_i} = \ga-\mfrac{\log (6\la+1)}{\log6} \\
\alphat_{\la,\max} & = \lim_{q\to-\infty} \tau_\mula'(q) = \max_{0\leq i \leq 3}\mfrac{\log p_{\la,i}}{\log r_i} =\max \Big(\ga,\ga-\mfrac{\log (6\la-1)}{\log6}\Big).
\end{align*}
Moreover, $\tau^*_\mula$ reaches its maximum at point $\alphat_{\la,\cL}$, with $\tau^*_\mula(\alphat_{\la,\cL})=1$, where
\begin{equation}\label{*lebas}
\alphat_{\la,\cL} =\tau_\mula'(0) = \ga -\frac{\log(36\la^2-1)}{4 \log 3+2\log 6}.
\end{equation}
\end{proposition}

\medskip

Observe that as long as $\la >\frac13$, then $\alphat_{\la,\max}=\ga$, and hence it depends only on $\la$ through the constant $\ga$. When $\la\leq \frac13$, $\alphat_{\la,\max}=\ga-\mfrac{\log(6\la-1)}{\log6}$ so that the dependence on $\la$ is more intricate. We will come back to this fact in Section \ref{*secdmula} below where in point $(7)$ we also comment on \eqref{*lebas}.

\begin{proposition}\label{recalls-mula2}
Let $\la\in(\frac16,\frac56)$ and consider $\mula$, $\tau_\mula$ and $\tau^*_\mula$.
\begin{itemize}[parsep=-0.15cm,itemsep=0.25cm,topsep=0.2cm,wide=0.175cm,leftmargin=0.85cm]
\item[$(i)$] For every $\alpha \in [\alphat_{\la,\min},\alphat_{\la,\max}]$,
\begin{equation}\label{form-mula}
d_\mula(\alpha)=\tau^*_\mula(\alpha) = \dimh(E_\mula(\alpha)).
\end{equation}
In particular, $\mula$ satisfies the multifractal formalism.
\item[$(ii)$] For every $\alpha\in[\alphat_{\la,\cL},\alphat_{\la,\max}]$, i.e. in the decreasing part of the spectrum,
\begin{equation*}
d_\mula(\alpha)=\tau^*_\mula(\alpha) = \dimh(\underline{E}_\mula^\geq(\alpha)).
\end{equation*}
\end{itemize}
\end{proposition}

\medskip

Let us make some remarks about these results.

\medskip

\noindent{\bf Comments on Proposition \ref{recalls-mula2}.}
\begin{itemize}[parsep=-0.15cm,itemsep=0.25cm,topsep=0.2cm,wide=0.175cm,leftmargin=0.85cm]
\item[$(a)$] The main reason for the equality between $d\mula(\alpha) = \dim \underline{E}_\mula(\alpha)$ and the Hausdorff dimension of $E_{\mula}(\alpha)$ (last equality of item $(i)$) is the application of the law of large numbers and Birkhoff's ergodic theorem.
\item[$(b)$] Given $\alpha\geq\alphat_{\la,\cL}$, item $(ii)$ follows from the large deviations theory (and in particular, Theorem 1 of \cite{BRMICHPEY}), and states that the majority of the points that have a lower local dimension greater than $\alpha $ are actually those with the limit local dimension equal to $\alpha$. This is not true when $\alpha < \alphat_{\la,\cL}$, in this case one has the opposite property that $\dimh(\{x\in[0,1]: h_\mula(x)\leq\alpha\}) = \tau^*_\mula(\alpha)$, but we will not use this later on.
\item[$(c)$] Writing that $\{x\in[0,1] : \ldimloc(\mula,x)\neq \alphat_{\la,\cL}\} = \bigcup _{n\geq 1} \{x\in[0,1] : \ldimloc(\mula,x)\leq \alphat_{\la,\cL}-1/n\}\cup \{x\in[0,1] : \ldimloc(\mula,x)\geq \alphat_{\la,\cL}+1/n\} $, comment $(b)$ combined with Proposition \ref{recalls-mula1} shows that every set appearing on the right hand-side has a Hausdorff dimension strictly less than 1, hence a Lebesgue measure 0. Hence, the set $\{x\in[0,1] : \ldimloc(\mula,x)=\alphat_{\la,\cL}\}$ has full Lebesgue measure in $[0,1]$.
\end{itemize}

The next result asserts that one can actually specify that those most important points in a given $\underline{E}_\mula(\alpha)$ are regular, i.e. they have limit digit frequencies in their decomposition \eqref{decomp-x} and \eqref{def-xn}.

\begin{lemma}\label{lem:Fbeta}
Let $\la\in(\frac16,\frac56)$ and consider $\mula$, $\tau_\mula$ and $\tau^*_\mula$. For any $(\beta_0,\beta_1,\beta_2,\beta_3)\in[0,1]^4$ set
\begin{equation}\label{def:Fbeta}
F_{\beta_0,\beta_1,\beta_2,\beta_3}=\{x\in [0,1]: x\text{ is regular and } \beta_i(x)=\beta_i \text{ for $i=0,\ldots,3$} \}.
\end{equation}
Then, for every $\alpha \in[\alphat_{\la,\min},\alphat_{\la,\max}]$, there exists a probability vector $(\beta_0,\beta_1,\beta_2,\beta_3)\in [0,1]^4$ such that $F_{\beta_0,\beta_1,\beta_2,\beta_3}\subset E_\mula(\alpha)$ and
\begin{equation*}
\dimh(F_{\beta_0,\beta_1,\beta_2,\beta_3}) = d_\mula(\alpha).
\end{equation*}
\end{lemma}

\begin{proof}
In Lemmas $2.4$ and $5.1$ of \cite{Gatzouras} it is proved that
given a probability vector $(\beta_i)_{i=0,\ldots,3}$, the set  $F_{\beta_0,\beta_1,\beta_2,\beta_3}$ has Hausdorff dimension equal to
\begin{equation}\label{dimFbeta}
\dimh(F_{\beta_0,\beta_1,\beta_2,\beta_3}) =\frac{\sum_{i=0}^3 \beta_i\log \beta_i}{\sum_{i=0}^3 \beta_i\log r_i}.
\end{equation}
Next, given a possible local dimension $\alpha$ for $\mula$, we call $q_\alpha\in \R$ the unique real number such that $\alpha=\tau'_{\mula}(q_\alpha)$ (which exists since $\tau_\mula$ is analytic concave). An immediate  computation using the implicit formula \eqref{eq:Lqspectrum} gives
\begin{equation*}
\alpha = \frac{\sum_{i=0}^3 p_{\la,i}^{q_\alpha} r_i ^{-\tau_\mula(q_\alpha)}\log p_{\la,i}}{\sum_{i=0}^3 p_{\la,i}^{q_\alpha} r_i ^{-\tau_\mula(q_\alpha)} \log r_i}.
\end{equation*}
The multifractal formalism ensures that $\tau_\mula^*(\alpha) = \alpha q_\alpha -\tau_\mula(q_\alpha)$. Then, setting $\beta_i=  p_{\la,i}^{q_\alpha} r_i ^{-\tau_\mula(q_\alpha)}$,    \eqref{eq:Lqspectrum} shows that $(\beta_i)_{i=0,\ldots,3}$ is a probability vector, and then \eqref{dimFbeta} yields
\begin{eqnarray*}
\dimh(F_{\beta_0,\beta_1,\beta_2,\beta_3})
& = & \frac{\sum_{i=0}^3 \beta_i\log \beta_i}{\sum_{i=0}^3 \beta_i\log r_i} = \frac{\sum_{i=0}^3 \beta_i \big( q_\alpha  \log p_{\la,i}-\tau_\mula(q_\alpha)\log r_i \big)}{\sum_{i=0}^3 \beta_i\log r_i}\\
& = & \alpha q_\alpha -\tau_\mula(q_\alpha) = \tau_\mula^*(\alpha),
\end{eqnarray*}
hence the result by \eqref{form-mula}.
\end{proof}

\subsection{Behavior of $d_\mula$ in terms of $\la$}\label{*secdmula}

\noindent

This gathers additional remarks, complementary to the previous ones but enlightening the dependence of $d_\mula$ on $\la$.

\medskip

\begin{figure}[hb!]
\centering
\captionsetup{width=1\textwidth}
\includegraphics[scale=0.45]{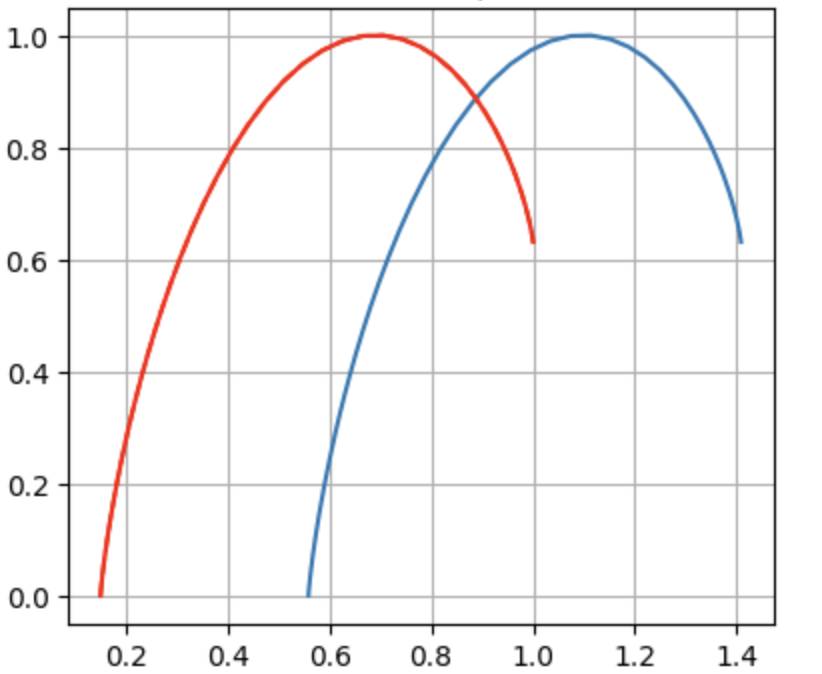}
\includegraphics[scale=0.45]{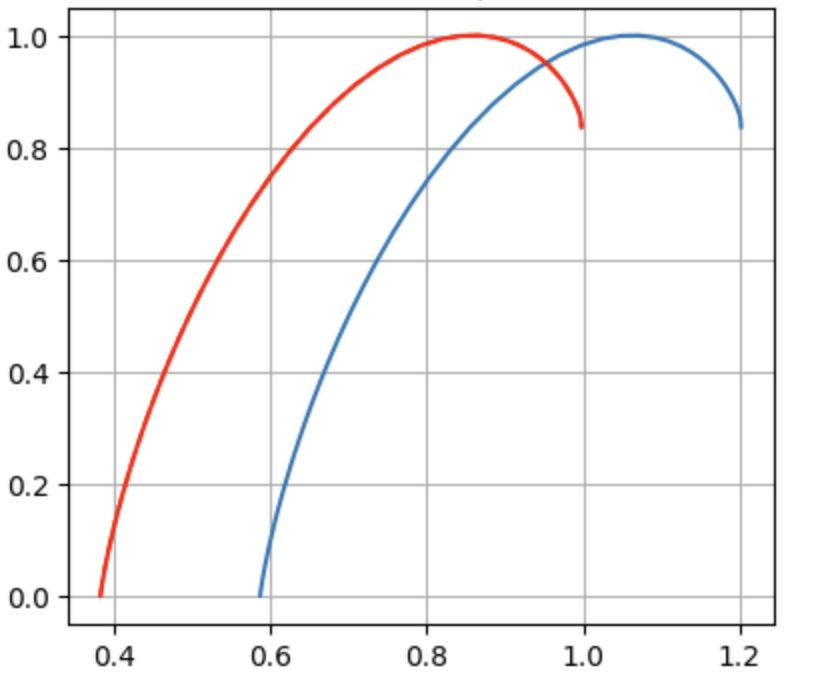}
\includegraphics[scale=0.45]{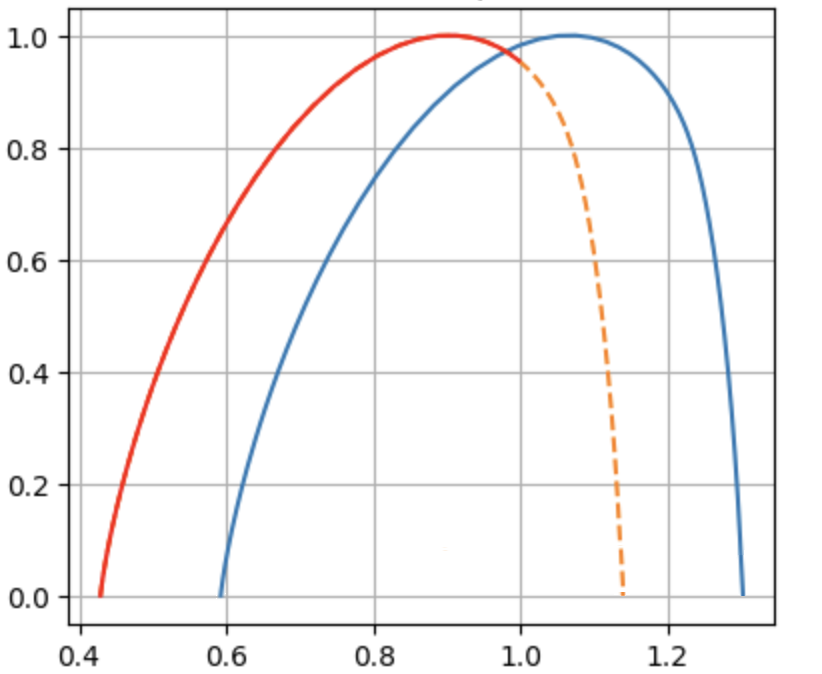}
\includegraphics[scale=0.45]{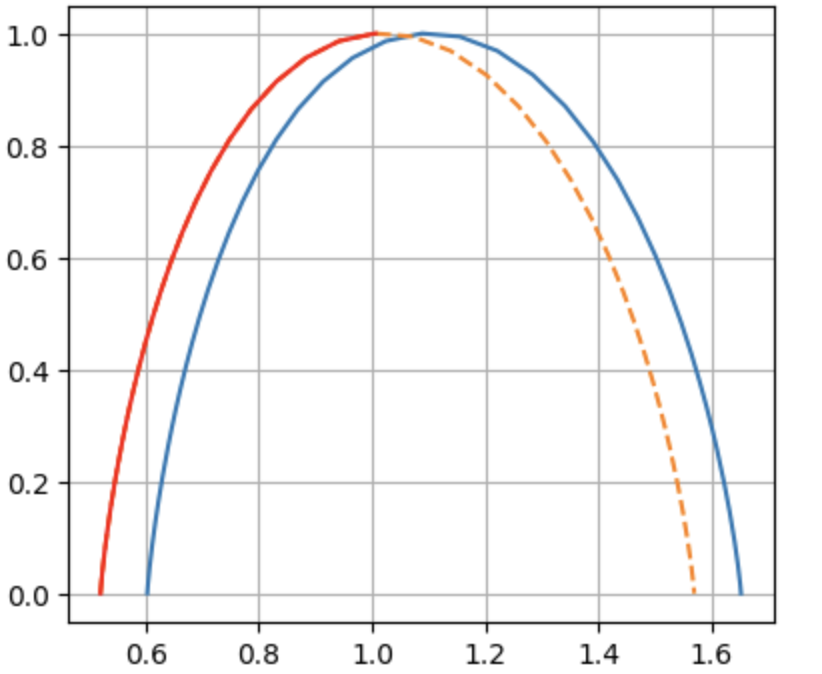}
\captionsetup{width=1\textwidth}
\caption{The multifractal spectrum $d_\mula$ of $\mula$ in blue, of $\mula$ translated by $1-\ga$ in dotted red, of $\fla$ in red, for $\la=0.6$, $0.4$, $0.3$ and $\frac{\sqrt2}{6}+0.01$ from top to bottom and left to right. The dotted red and red graphs coincide as long as $\la>\frac13$.}
\label{fig-spectra}
\end{figure}

\begin{itemize}[parsep=-0.15cm,itemsep=0.25cm,topsep=0.2cm,wide=0.175cm,leftmargin=0.85cm]
\item[$(1)$] The minimal and maximal exponents $\alphat_{\la,\min}$ and $\alphat_{\la,\max}$ are analytic functions of $\la$, since the dependence on $\ga$ in $\la$ given by \eqref{def-gamma} is analytic. Also, in the range of parameters $\la$ of interest, $d_\mula(\alphat_{\la,\min})=\tau_{\mula}^*(\tau_\mula'(+\infty))=0$.

\item[$(2)$] When $\la\in (\frac{\sqrt2}6,\frac13)$, we have
\begin{equation*}
\hspace{0.75cm}\alphat_{\la,\min}= \ga-\frac{\log (6\la+1)}{\log6} =\frac{\log p_{\la,1}}{\log r_1}  < \frac{\log p_{\la,0}}{\log r_0} < \frac{\log p_{\la,2}}{\log r_2}= \ga-\frac{\log (6\la-1)}{\log6}=\alphat_{\la,\max}.
\end{equation*}
The level set $\underline{E}_\mula(\alphat_{\la,\max})$ consists of
the points for which $\lim_{n\to\infty}\frac{\beta_2(x,n)}n=1$. Indeed, recalling \eqref{eq:itermu} which explains how $\mula$ spreads the mass of an interval of generation $n$ on its four subintervals, the only possibility for a point $x$ to have a local dimension for $\mula$ equal to $\frac{\log p_{\la,2}}{\log r_2}$ is that most of the digits of $x$ equal $2$, i.e. $\lim_{n\to\infty}\frac{\beta_2(x,n)}n = 1$. Then, applying again the multifractal formalism satisfied by $\mula$, this set has a Hausdorff dimension equal to $d_\mula(\alphat_{\la,\max})=\tau^*_\mula(-\infty)=0$.

{\item[$(3)$] The mapping $(\la,\alpha)\mapsto \tau^*_\mula(\alpha)$ is continuous on $\{(\la,\alpha): \la\in(\frac{\sqrt{2}}{6},\frac{1}{3}),\alpha\in [\alphat_{\la,\min},1]\}$, by properties of the Legendre transform. Hence, applying \eqref{form-mula}, the mapping
\begin{equation}\label{dmu-map}
\la\in \big(\mfrac{\sqrt2}{6},\mfrac13\big) \mapsto \mathrm{Gr}(d_\mula) \in\cK(\R^2)
\end{equation}
is continuous, where we recall that $\mathrm{Gr}(d_\mula)$ denotes the (compact) graph of $d_\mula$ and $\cK(\R^2)$ the space of compact subsets in $\R^2$ equipped with the Hausdorff metric. Propositions \ref{prop_lower}, \ref{*sppropa} and \ref{*sppropb} will imply that for these parameter values the
spectrum of $\fla$ is a translation of that of $\mula$, but truncated. Figure \ref{fig:HVKsection6} displays $d_\fla(\alpha)$ and $d_\mula(\alpha+\ga-1)$ for the classical von Koch function, i.e. when $\la=\frac{\sqrt 3}6.$ On $[\alpha_{\la,\min},1]$ these two curves coincide, but the domain of $d_\mula(\alpha+\ga-1)$ is larger and as Remark $(2)$ implies it reaches the $x$ axis at $\alphat_{\la,\max}-\ga+1$.} See also Figure \ref{fig-spectra} for examples with different parameter values.

\item[$(4)$] Precisely at $\la=\frac13$,
\begin{equation}\label{coincidence}
\frac{\log p_{\la,2}}{\log r_2}  = \frac{\log p_{\la,0}}{\log r_0},
\end{equation}
and now the level set $\underline{E}_\mula(\alphat_{\la,\max})$ contains the points for which $\lim_{n\to\infty} \frac{\beta_1(x,n)}n=0$. In particular, it contains the set $\{x\in[0,1] : \beta_1(x,n)=0 \mbox{ for every integer $n$}\}$, which is the inhomogeneous Cantor set described in the introduction, with Hausdorff and box dimensions both equal to the solution to \eqref{def-dim-s}. It is not difficult to see that $d_{\mu_\frac13}(\alphat_{\frac13,\max})$ also equals $s$, again using the multifractal formalism and the coincidence \eqref{coincidence}. So, even if the multifractal spectrum $\alpha\mapsto d_{\mu_\frac13}(\alpha)$ is still continuous, there is a jump discontinuity at $\la= \frac13$ from both sides since $d_{\mu_\frac13}(\alphat_{\frac13,\max})> \frac{\log 2}{\log3}=d_\mula(\alphat_{\la,\max})$ for $\la\in(\frac13,\frac56)$, and $d_\mula(\alphat_{\la,\max})=0$ for $\la\in(\frac16,\frac13)$, as we observed in $(2)$.

\item[$(5)$] When $\la\in (\frac13,\frac56)$, we have
\begin{equation*}
\alphat_{\la,\min}= \ga-\frac{\log(6\la+1)}{\log6}= \frac{\log p_{\la,1}}{\log r_1}  < \frac{\log p_{\la,2}}{\log r_2} < \frac{\log p_{\la,0}}{\log r_0} =\ga =\alphat_{\la,\max}.
\end{equation*}
In particular, the level set $\underline{E}_\mula(\alphat_{\la,\max})$ contains the points for which $\lim_{n\to \infty}\frac{\beta_{0,3}(x,n)}n=1$. This set has Hausdorff dimension $\tau^*_\mula(\alpha_{\la,\max})=\lim_{q\to-\infty}\tau_\mula'(q) = \frac{\log2}{\log3}$, and is independent of $\la$ in this range (it contains and has the same dimension as the triadic Cantor set).

\item[$(6)$] When $\la\in(\frac13,\frac56]$, the mapping \eqref{dmu-map} is again continuous, still by \eqref{form-mula}.

\item[$(7)$] When $x$ is chosen randomly Lebesgue uniformly in $[0,1]$, by the ergodic theorem applied to $T$, for Lebesgue almost every $x$ we have
\begin{equation}\label{randomLebesgue}
\frac{\beta_{0,3}(x,n)}{n}\to \frac23, \,
\frac{\beta_{1,2}(x,n)}{n} \to \frac26, \,
\frac{\beta_1(x,n)}{n} \to \frac16 \text{ and }
\frac{\beta_2(x,n)}{n} \to \frac16.
\end{equation}
Thus, we get heuristically $|I_n(x)| \approx 3^{-2n/3}6^{-n/3}$ and, using \eqref{eq:muIn},
\begin{equation*}
\mula(I_n(x)) \approx p_{\la,0}^{n/3}p_{\la,1}^{n/6}p_{\la,2}^{n/6}p_{\la,3}^{n/3}.
\end{equation*}
Considering $\frac{\log \mula(I_n(x))}{\log |I_n(x)|}$, simplifying and taking its limit when $n$ tends to infinity explains the Lebesgue-almost sure value $\alphat_{\la,\cL}$ given in \eqref{*lebas}.
\end{itemize}

\addtocontents{toc}{\vspace{0.2cm}}%
\section{The lower bound for $d_\fla$ when $\la\in(\frac16,\frac56)$ and the decreasing part of $d_\fla$}\label{sec:spectrumirr}

The lower bound we are going to prove for $d_\fla$ holds for every $\la\in (\frac16,\frac56)$ but it turns out that the analysis in the previous section allows us to deal with the subcase when $\la\in(\frac13,\frac56)$ in an easier way. Indeed, on one hand, Proposition \ref{propvaluedimloc} gives the explicit value of $\ldimloc(\mula,x)$, and on the other hand, Proposition \ref{propholdabis} provides an upper bound for the pointwise exponent of $\fla$. Recalling Definition \ref{*defI} of the set $\cI$, the combination of these facts yields the following key result.

\begin{corollary}\label{coro-compare}
For every $x\in [0,1]$, $h_\fla(x) \leq 1-\ga +\ldimloc(\mula,x) $, with equality when $x\in\cI$ and the right hand-side belongs to $[0,1]$.
\end{corollary}

Observe also that $1-\ga + \alpha_{\la,\min} \to  0$ when $\la\nearrow \frac56$, which confirms that our analysis is not applicable when $\la>\frac56$. Proposition \ref{prop-2} below gives the multifractal spectrum of $\fla$, i.e. Theorem \ref{spectruma}, when $\la\in (\frac13,\frac56)$.

\begin{proposition}\label{prop-2}
Let $\la\in (\frac13,\frac56)$. Then $d_\fla(\alpha)= \tau^*_\mula(\alpha+\ga-1)$  for every $\alpha\in [\alpha_{\la,\min},1]$.
\end{proposition}

\begin{proof}
In this range of parameters, by Definition \ref{*defI} and Lemma \ref{lemslopbound} we have $\cI= [0,1]\setminus\cEt$. When $x\in \cEt$, $h_\fla(x)=1$, thus
\begin{equation*}
d_\fla(1)\geq \dimh(\cEt)= \frac{\log2}{\log3} = d_\mula(\ga) =d_\mula(1+\ga-1).
\end{equation*}
When $x\notin\cEt$, then by Lemma \ref{lemslopbound}
$x\in\cI$, and by Corollary \ref{coro-compare} $h_\fla(x) = \ldimloc(\mula,x) +1-\ga$. This yields the result, since for every $\alpha\geq 0$, $E_\fla(\alpha)=\underline{E}_\mula(\alpha+\ga-1)$.
\end{proof}

Now we look for an extension of Proposition \ref{prop-2} to the subcase when $\la\in(\frac16,\frac13]$.

\begin{lemma}\label{lemregulinf}
Let $\la\in(\frac16,\frac56)$. Let $(\beta_i)_{i=0,\ldots,3}$ be a probability vector and $x\in F_{\beta_0,\beta_1,\beta_2,\beta_3}$ (recall \eqref{def:Fbeta}).
\begin{enumerate}[parsep=-0.15cm,itemsep=0.25cm,topsep=0.2cm,wide=0.175cm,leftmargin=0.5cm]
\item[$(i)$] If $\beta_1 \log(6\lambda+1)+\beta_2\log(6\lambda -1)>0$, then $x\in\cI$. 
\item[$(ii)$] If $\beta_1 \log(6\lambda+1)+\beta_2 \log(6\lambda -1) \leq 0$, then $h_\fla(x)=1$.
\end{enumerate}
\end{lemma}

\begin{proof}
If $\lll \in (\frac13,\frac56]$ then only case $(i)$ is possible and as we observed above
when $x\notin\cEt$, then by Lemma \ref{lemslopbound}
$x\in\cI$. If $x\in\cEt$ then $x\notin F_{\beta_0,\beta_1,\beta_2,\beta_3}$ for a probability vector $(\beta_i)_{i=0,\ldots,3}$ satisfying the hypothesis of $(i)$.

\medskip

Assume that $\la\in(\frac16,\frac13)$ so that $0<6\la-1<1$. Since $x\in F_{\beta_0,\beta_1,\beta_2,\beta_3}$, for every $\eps>0$ there exists $N_\eps \geq1$ such that
\begin{equation}\label{*F3*a}
\Big|\frac{\beta_i(x,n)}{n}-\beta_i\Big |<\eps \text{ for all $i=0,1,2,3$ and $n\geq N_\eps$}.
\end{equation}

Let us prove $(i)$. From the assumption in $(i)$, $x\in F_{\beta_0,\beta_1,\beta_2,\beta_3}$, \eqref{defE} and \eqref{defEtilde} it also follows that $x\notin\cEt$. Therefore, going back to Lemma \ref{lemsloprec}, for all $n\geq0$,
\begin{equation}\label{slopabsrec}
|m_{n+1}(x)|
= \left\{
\begin{array}{lll}
|m_n(x)| & & \text{if $u_n(x)\in\{0,3\}$,} \\
6\la\sqrt{1+m^2_n(x)} + |m_n(x)| & \geq (6\la+1)|m_n(x)|& \text{if $u_n(x)=1$,} \\
6\la\sqrt{1+m^2_n(x)} - |m_n(x)| & \geq (6\la-1) |m_n(x)|& \text{if $u_n(x)=2$.}
\end{array}\right.
\end{equation}
We may assume without loss of generality that $m_1(x)=6\la$, hence
\begin{equation}\label{slopabsmin}
|m_n(x)| \geq 6\la(6\la+1)^{\beta_1(x,n)-1}(6\la-1)^{\beta_2(x,n)}.
\end{equation}
By assumption we can select $\eps>0$ small enough such that
\begin{equation}\label{*F2*e*}
(\beta_1-\eps)\log(6\lambda+1)+(\beta_2+\eps)\log(6\lambda -1)>0.
\end{equation}
Recalling that $\log(6\la-1)<0$, it follows that, for all $n\geq N_\eps$,
\begin{align}
\nonumber\log |m_n(x)|
& \geq \log(6\la) + (\beta_1(x,n)-1)\log(6\la+1) +\beta_2(x,n)\log(6\la-1) \nonumber \\
\label{slopabsminbis}
& \geq \log(6\la) + n\Big((\beta_1-\eps)\log(6\la+1) +(\beta_2+\eps)\log(6\la-1)\Big)-\log(6\la+1).
\end{align}
Hence, $\log |m_n(x)|\to +\infty$, i.e. $x\in\cI$.

\medskip

Let us prove $(ii)$. Suppose that $x\in F_{\beta_0,\beta_1,\beta_2,\beta_3}$ is fixed and assume without loss of generality that $m_1(x)>0$. By Proposition \ref{propholda} we only have to prove that $h_\fla(x)\geq1$. Given $\eps_1>0$, choose $\eps\in(0,\eps_1)$ small enough such that
\begin{equation}\label{*F2*f}
\max\big(\beta_1 \log(6\la+1+\eps) +\beta_2 \log(6\la -1+\eps) , 2\eps |\log(6\la+1+\eps)| , 2\eps|\log(6\la-1+\eps)|\big)  \leq \frac{\eps_1}3
\end{equation}
and $0<6\la-1+\eps<1$. Fix $N_\eps\geq 1$ such that \eqref{*F3*a} holds for $i=0,\ldots,3$, and $M_\la>0$ so large that for all $m\geq M_\la$,
\begin{equation}\label{*F3*c}
\Big |1+6\la\sqrt{1+\mfrac1{m^2}} - (6\la+1) \Big |< \eps  \ \text{ and } \
\Big |-1+6\la\sqrt{1+\mfrac1{m^2}} - (6\la -1) \Big |< \eps.
\end{equation}

\medskip

\noindent{\it Case 1:} The sequence $(|m_n(x)|)_{n\geq 1}$ is bounded by some $M>0$.

Then, by Lemma \ref{lemslopbound}, $0< m_n^-= \sqrt{36\la^2-1} \leq |m_n(x)|\leq M$, and since $\ell_n(x) \leq 3^{-n}$, one has $\lim_{n\to\infty} \frac{\log |m_n(x)|}{|\log \ell_n(x)|}=0$. Then, Proposition \ref{propholda} yields $h_\fla(x)=1$.

\medskip

\noindent{\it Case 2:} $\lim_{n\to\infty} |m_n(x)|=+\infty$.

Let $N\geq N_\eps$ be so large that for every $n\geq N$, $|m_n(x)|\geq M_\la$. By construction, it follows from \eqref{*F3*c} and Lemma \ref{lemsloprec} that, for all $n\geq N$,
\begin{equation}\label{*67a}
|m_{n+1}(x)| \leq \left\{
\begin{array}{ll}
|m_n(x)| & \text{if $u_n(x)\in\{0,3\}$,} \\
(6\la+1+\eps)|m_n(x)| & \text{if $u_n(x)=1$,} \\
(6\la-1+\eps)|m_n(x)| & \text{if $u_n(x)=2$.}
\end{array}\right.
\end{equation}
This implies that for every $n\geq N$,
\begin{equation*}
|m_n(x)|\leq |m_N(x)| (6\la+1+\eps)^{\beta_1(x,n)-\beta_1(x,N)}(6\la -1+\eps)^{\beta_2(x,n)-\beta_2(x,N)}.
\end{equation*}
By the choice of $N_\eps$ we can use \eqref{*F3*a} to obtain, for $n\ge 3N$,
\begin{align}\label{eq-la13-2}
|m_n(x)| & \leq |m_N(x)|(6\la+1+\eps)^{(\beta_1+2\eps)(n-N)}
(6\la-1+\eps)^{(\beta_2-2\eps)(n-N)} \nonumber \\
& \leq  C (6\la+1+\eps)^{(\beta_1+2\eps)n}(6\la-1+\eps)^{(\beta_2-2\eps) n},
\end{align}
where $C>0$ depends on $x$ and $N$ only. Taking logarithm and using \eqref{*F2*f}, we infer
\begin{equation}\label{*logmn}
\log (\sqrt{36\la^2-1})\leq \log|m_n(x)|\leq \log  C  + n\eps_1.
\end{equation}
The parameter $\eps_1$ being arbitrarily small, one deduces that $\lim_{n\to\infty} \frac{\log |m_n(x)|}{|\log\ell_n(x)|}=0$, and again Proposition \ref{propholda} gives $h_\fla(x)=1$.

\medskip

\noindent{\it Case 3:} $\liminf_{n\to\infty} |m_n(x)|< \limsup_{n\to \infty}|m_n(x)|=+\infty$.

In this case the main difficulty is that \eqref{*F3*c} works only for $|m_n(x)|\geq M_\la$, but the idea is simple since for the other cases $m_n(x) $ is bounded.

There exist $M\geq M_\la$ and a unique strictly increasing sequence $(n_k)_{k\geq 1}$ such that $n_1\geq N_\eps$, $|m_{n_1}(x)|\leq M$ and, for every $k\geq 1$,
\begin{equation*}
\left\{
\begin{array}{lll}
& |m_n(x) | > M & \mbox{ when } n_{2k-1} < n \leq n_{2k},\\
\sqrt{36\la^2-1}\leq & |m_n(x)| \leq M & \mbox{ when } n_{2k} < n \leq n_{2k+1}.
\end{array}\right.
\end{equation*}
Let $k\geq 1$. When $n\geq n_{2k-1} \geq N_{\eps}$, one has by \eqref{*F3*a}
\begin{equation*}\label{eq-beta}
|\beta_i(x,n)-\beta_i(x,n_{2k-1})| \leq (n-n_{2k-1})\beta_i + 2n \eps, \mbox{ for $i=0,\ldots,3$.}
\end{equation*}
Recall that $6\la-1+\eps<1$, and hence $\eps<1$ as well. Using \eqref{*F2*f} to obtain the last line we deduce
\begin{align}\label{eq-la13}
|m_n(x)|
& \leq |m_{n_{2k-1}}(x)|(6\la+1+\eps)^{\beta_1(x,n)-\beta_1(x,n_{2k-1})}(6\la-1+\eps)^{\beta_2(x,n)-\beta_2(x,n_{2k-1})} \nonumber \\
& \leq M(6\la+1+\eps)^{(n-n_{2k-1})\beta_1 + 2n \eps}(6\la-1+\eps)^{ (n-n_{2k-1})\beta_2-2n\eps} \nonumber \\
& \leq M(6\la+1+\eps)^{(n-n_{2k-1})(\beta_1 + 2 \eps)+2n_{2k-1}\eps}(6\la-1+\eps)^{ (n-n_{2k-1})(\beta_2-2\eps)-2n_{2k-1}\eps}
\nonumber \\
& \leq M e^{(n-n_{2k-1})\eps_1}(6\la+1+\eps)^{2n_{2k-1}\eps}(6\la-1+\eps)^{-2n_{2k-1}\eps}\leq M\cdot 4^{3n\eps_1},
\end{align}
where we used $\eps<\eps_1$, $\la<\frac13$, $e<4$. Therefore,
\begin{equation*}
\log |m_n(x)|\leq \log M + 4n\eps_1\log 10
\end{equation*}
holds for $n_{2k-1}< n\leq n_{2k}$, where $n\leq n_{2k}$ was used to ensure that $|m_j(x)|>M_\la$ for all $n_{2k-1}\leq j < n$ in order to apply \eqref{*67a} iteratively. When $n_{2k}< n\leq n_{2k+1}$, $|m_n(x)|\leq M$, so \eqref{eq-la13} also holds in this case. Since $\eps_1$ is arbitrarily small, we obtained that $\lim_{n\to \infty} \frac{\log |m_n(x)|}{|\log\ell_n(x)|}=0$, and Proposition~\ref{propholda} gives $h_\fla(x)=1$.

\medskip

Finally, the case $\la=\frac13$ is treated similarly.

For $(i)$ one can use that $\log (6\la-1)=0$ and hence by the assumption $\beta_1>0$. Moreover \eqref{slopabsrec} implies that $|m_{n+1}(x)|\geq |m_n(x)|$ even in case of $u_n(x) =2$.

For $(ii)$ we can argue similarly as above except that in \eqref{eq-la13-2} and \eqref{eq-la13}, $-2n\ep$ is replaced by $+2n\eps$, since $6\la-1=1$ and this time we cannot assume that $6\la-1+\eps<1$.
\end{proof}

\begin{proposition}\label{prop_lower}
Let $\la\in(\frac16,\frac56)$. Then, for every $\alpha\in[\alpha_{\la,\min},1]$, $d_\fla(\alpha)\geq \tau^*_\mula(\alpha+\ga-1)$.
\end{proposition}

\begin{proof}
Let $\alpha\in[\alpha_{\la,\min},1]$ and set $\alphat=\alpha+\ga-1$. By Lemma \ref{lem:Fbeta}, let us choose parameters $(\beta_i)_{i=0,\ldots,3}$ such that $F_{\beta_0,\beta_1,\beta_2,\beta_3}\subset E_\mula(\alphat)$ and
\begin{equation*}
d_\mula(\alphat)=\tau^*_\mula(\alphat)=\dimh(F_{\beta_0,\beta_1,\beta_2,\beta_3}).
\end{equation*}
Next, we show that $F_{\beta_0,\beta_1,\beta_2,\beta_3}\subset E_\fla(\alpha)$. Suppose $x\in F_{\beta_0,\beta_1,\beta_2,\beta_3}$. Then its digit frequencies \eqref{*F1*a} exist and equal $\beta_i(x)=\beta_i$ for $i=0,\ldots,3$, and then the local dimension of $\mula$ in formula \eqref{value-dimloc} is given by
\begin{align*}
\alphat=\ldimloc(\mula,x) = \ga - \frac{ \beta_1\log(6\la+1)+\beta_2\log(6\la-1)}{(\beta_0+\beta_3)\log 3 + (\beta_1+\beta_2)\log 6}.
\end{align*}

If $\alpha<1$ then $\alphat<\ga$, which implies that $\beta_1 \log(6\la+1)+\beta_2\log(6\la -1)>0$. Hence, by $(i)$ of Lemma \ref{lemregulinf}, we have $x\in\cI$. Thus, by Corollary \ref{coro-compare}, we have $h_\fla(x) = 1-\ga+\alphat = \alpha$ and hence $F_{\beta_0,\beta_1,\beta_2,\beta_3}\subset E_\fla(\alpha)$. Therefore, $d_\fla(\alpha) \geq \dimh (F_{\beta_0,\beta_1,\beta_2,\beta_3}) = \tau^*_\mula(\alphat) = \tau^*_\mula(\alpha+\ga-1)$.

If $\alpha=1$ then $\widetilde \alpha=\ga$, which implies that $\beta_1\log(6\la+1)+\beta_2\log(6\la-1)=0$. Hence, by $(ii)$ of Lemma \ref{lemregulinf}, we have $h_\fla(x)=\alpha=1$ and hence $F_{\beta_0,\beta_1,\beta_2,\beta_3}\subset E_\fla(1)$. Therefore, $d_\fla(\alpha) = d_\fla(1) \geq \dimh(F_{\beta_0,\beta_1,\beta_2,\beta_3}) = \tau^*_\mula(\alphat) = \tau^*_\mula(\ga)= \tau^*_\mula(\alpha+\ga-1)$.
\end{proof}

Finally, we conclude on the decreasing part of the multifractal spectrum $d_\fla$ when $\la\in(\frac{\sqrt2}6,\frac13]$.

\begin{proposition}\label{*sppropa}
Let $\la\in(\frac{\sqrt2}6,\frac13]$. Then, for every $\alpha\in [\alpha_{\la,\cL},1]$, $d_\fla(\alpha)=\tau^*_\mula(\alpha+\ga-1)$.
\end{proposition}

\begin{proof}
Let $\alpha\in[\alpha_{\la,\cL},1]$. By Proposition \ref{prop_lower}, it is enough to prove the upper bound $d_\fla(\alpha)\leq \tau^*_\mula(\alpha+\ga-1)$. By Propositions \ref{propholdabis} and \ref{propvaluedimloc}, if $x\in E_\fla(\alpha)$, then necessarily $x\in \underline{E}_\mula^\geq(\alpha+\ga-1)$. Thus item $(ii)$ of Proposition \ref{recalls-mula2} implies that
$d_\fla(\alpha)\leq \dimh\big(\underline{E}_\mula^\geq(\alpha+\ga-1)\big) = \tau^*_\mula(\alpha+\ga-1)$,
hence the result.
\end{proof}

\addtocontents{toc}{\vspace{0.2cm}}%
\section{The upper bound for the increasing part of $d_\fla$ when $\la\in [\frac{\sqrt2}6,\frac13)$}\label{sec:increasespec}

It is interesting that the above arguments are not sufficient to prove that $d_\fla(\alpha)=\tau_{\mula}^*(\alpha+1-\ga)$ in the increasing part of the spectrum, since item $(ii)$ of Proposition \ref{recalls-mula2} does not apply any more. We have to develop a much more involved argument to prove the upper bound. Proposition \ref{prop_lower} and Proposition \ref{*sppropb} below will imply \eqref{res-spectrum}, i.e. item $(ii)$ of Theorem \ref{spectruma}.

\begin{proposition}\label{*sppropb}
Let $\la\in [\frac{\sqrt2}6,\frac13]$. Then, for every $\alpha\in [\alpha_{\la,\min},\alpha_{\la,\cL}]$, $d_\fla(\alpha)\leq \tau^*_\mula(\alpha+\ga-1)$.
\end{proposition}

\medskip

The rest of this section is devoted to the proof of Proposition \ref{*sppropb}.

\medskip

\noindent{\it Proof ideas.} To make the reading of the proof easier, we present here the strategy and the different steps of the proof.
\begin{itemize}[parsep=-0.15cm,itemsep=0.25cm,topsep=0.2cm,wide=0.175cm,leftmargin=0.5cm]
\item[$\bullet$] Observe that the case $\la=\frac{\sqrt 2}6$ and $\alpha_0=1$ is obvious since in this case $d_\fla(\alpha_0)=1=\tau^*_\mula(\alpha_0+\ga-1)$. So in the rest of our argument for the case $\la=\frac{\sqrt 2}6$ we suppose $\alpha_0<\alpha_{\la,\cL}= 1$.
\item[$\bullet$] Suppose that for some $\alpha_0 \in [\alpha_{\la,\min},\alpha_{\la,\cL})$, points not in $\cI$ can modify the spectrum, that is, $\dimh(\{x\notin \cI: h_\fla(x)=\alpha_0\})> \dimh(\{x\in \cI: h_\fla(x)=\alpha_0\})$. From Corollary \ref{*F3*a}, if $\alpha_0<1$ and $x\in\cI$, then $h_\fla(x)=1-\ga + \underline{\dim}_{\text{loc}}(\mula,x)$. Therefore $\{x\in \cI: h_\fla(x)=\alpha_0\}\subset \{x : 1-\ga + \underline{\dim}_{\text{loc}}(\mula,x) = \alpha_0\} = \underline{E}_\mula(\alpha_0+\ga-1)$. And hence $\dimh(\{x\in \cI: h_\fla(x)=\alpha_0\}) \leq \tau^*_\mula(\alpha_0+\ga-1)$. When $x\notin\cI$,
\begin{equation}\label{*notinI}
\liminf_{n\to\infty} |m_n(x)|<+\infty.
\end{equation}
Since $|m_n(x)|$ is small infinitely often, for a sufficiently large $M$ select a sequence $n_j(x)\to\infty$ such that for every $j$,
\begin{equation}\label{*mnx26a}
\frac{M}{2(6\la+1)} < |m_{n_j}(x)|< M.
\end{equation}
Our goal is to build a large set of regular points $x$ where $h_\fla(x)$ is small. In order to obtain such a small $h_\fla(x)$, it will be necessary that sometimes $|m_n(x)|$ is ``rapidly increasing'' for $n_j(x)\leq n\leq n'_j(x)$. Using these rapidly increasing stretches, one can build a self-similar set $\bbF$, attractor of an IFS $\{\widetilde\Phi_k\}_{k=1,\ldots,K}$ satisfying the Open Set Condition, which contains a subset of points in $\cI$ with local $\mula$-dimension equal to some $\alpha_\bbF$ with Hausdorff dimension larger than $\tau^*_\mula(\alpha_\bbF)$, leading to a contradiction.
\item[$\bullet$] Therefore, proceeding towards a contradiction, suppose that for $\alpha_0\in (\alpha_{\la,\min},\alpha_{\la,\cL})$
\begin{equation}\label{ineg-dim}
1=\dimh(E_\fla(\alpha_{\la,\cL}))>\dimh(E_\fla(\alpha_0))=:d_0> \tau^*_\mula(\alpha_0+\ga-1) =: \dos.
\end{equation}
We choose $d_1,\ldots,d_4$ such that
\begin{equation}\label{eq-di}
\dos < d_2 < d_3 < d_1 < d_0 < d_4 < 1,
\end{equation}
and by continuity and monotonicity of $\tau^*_\mula$, we select $\eps_0>0$ such that $\alpha_{\la,\min}<\alpha_0+\eps_0<\alpha_{\la,\cL}$ and
\begin{equation}\label{*epso}
d_2>\tau^*_\mula(\alpha+\ga-1)
\text{ for } \alpha\in(\alpha_{\la,\min},\alpha_0+\eps_0),
\end{equation}
see Figure \ref{fig:HVKsection6} below. As later arguments show, $d_4$ can be chosen as close as we wish to $1$, therefore we select $d_2$ and $d_3$ before fixing $d_4$, this way we can ensure
\begin{equation}\label{*d2d1d5}
d_3+d_4-d_2>1.
\end{equation}
\end{itemize}

\begin{figure}[h!]
\includegraphics[width=0.75\textwidth]{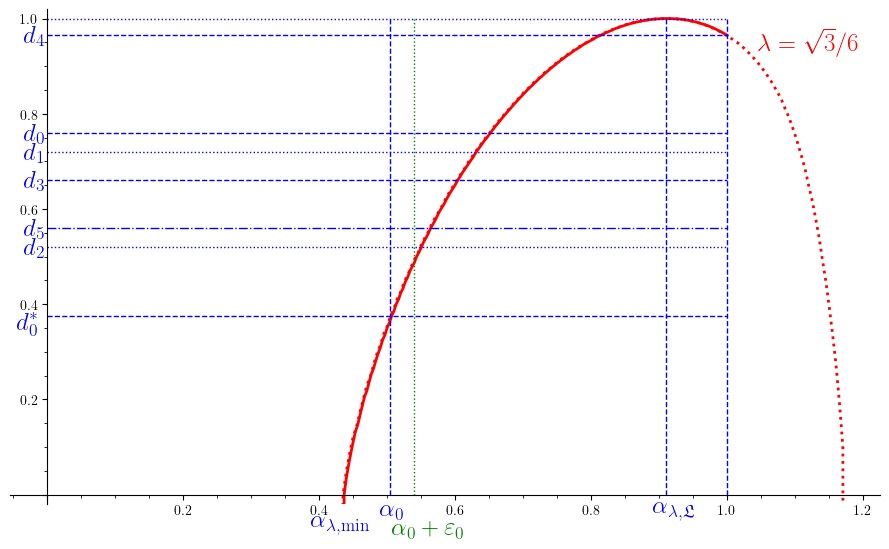}
\captionsetup{width=1\textwidth}
\caption{Illustration to the proof of Proposition \ref{*sppropb}. The thick red continuous curve is $d_\fla(\alpha)$ when $\la=\frac{\sqrt3}6$, while the thick dotted line shows $d_\mula(\alpha+\ga-1)$.}
\label{fig:HVKsection6}
\end{figure}

\subsection{Construction of the IFS $\{\widetilde\Phi_k\}_{k=1,\ldots,K}$ and its attractor $\bbF$}

\noindent

In this section, we build an IFS $\{\widetilde\Phi_k\}_{k=1,\ldots,K}$ satisfying the Open Set Condition. Each map $\widetilde\Phi_k$ will be a bijective affine map from $[0,1]$ onto some interval $\bbJ_k$ (the monotonicity of $\widetilde\Phi_k$ will be explicit). All intervals $\bbJ_k$ are of the form $\overline{I}_n(x)$ for some $x$ and $n$ (recall our notation $\overline{I}_n(x)$ for the closed interval $\overline{I_n(x)}$, the closure of the open interval $I_n(x)$). The delicate part is to identify these intervals $\bbJ_k$, since we want to control very precisely the behavior of $\mula$ and $\fla$ on them. To do so, we will proceed in two steps, first building some intervals on which we can ensure that $|m_n(x)|$ is quite large, and then inside each of them we will find smaller intervals where the slopes $m_n(x)$ are precisely controlled.

\subsubsection{The family of the intervals $\bbJ_j'$}\label{*ssecJj}

\noindent

Recall \eqref{ineg-dim}. Since $\dimh(E_\fla(\alpha_0))=d_0$ and $d_1<d_0$, there exists a compact set $E\subset E_\fla(\alpha_0)\setminus (\cE\cup\cI)$ such that $\cH^{d_1}(E)=+\infty$. Since $E\cap\cI=\emptyset$, one has
\begin{equation*}
E = \bigcup_{M\geq1} \{x\in E: |m_n(x)|<M \mbox{ for infinitely many integers } n\}.
\end{equation*}

Recall that by Corollary \ref{coro-compare} and \eqref{ineg-dim}, $\dimh(\{ x\in\cI:h_\fla(x)=\alpha_0\})=\tau^{*}_\mula(\alpha_0+\ga-1)=\dos<d_1$. Since $\cH^{d_1}(E)=+\infty$, we can suppose, up to a small modification of $d_1$, that the compact set $E$ was chosen so that there exists $M>100$ such that for all $x\in E$, $|m_{\overline n_j(x)}(x)|<M$ for infinitely many integers $\overline n_j(x)$, and $\cH^{d_1}(E)=+\infty$, in a heuristic way one should think that $|m_{\overline n_j(x)}(x)|$ is ``very small''.

After fixing in \eqref{*epso} $\eps_0$ next we will choose sufficiently small $\eps_1$ and $\eps_2$. These will have to satisfy \eqref{*epschoose}, but otherwise their choice will be independent from $\eps_0$.

On the other hand, let $\eps_1>0$ be such that $\alpha_0+\eps_1<1$. Suppose that $x\in E$ is fixed.

The fact that $h_\fla(x)=\alpha_0$ implies, for instance by \eqref{expholda}, that for infinitely many integers $\overline n'_j(x)$ one has $|m_{\overline n'_j(x)}(x)||I_{\overline n'_j(x)}(x)|\geq |I_{\overline n'_j(x)}(x)|^{\alpha_0+\eps_1}$, i.e.
\begin{equation}\label{*mnx38a}
\log |m_{\overline n'_j(x)}(x)| \geq (1-\alpha_0-\ep_1)|\log \ell_{\overline n'_j(x)}(x)|,
\end{equation}
in a heuristic way one should think that $|m_{\overline n'_j(x)}(x)|$ is ``large''. Now we want to interleave small and large values of $|m_{n}(x)|$.

We choose a suitably small $\eps_2>0$ keeping in mind that we  will need \eqref{*epschoose} in the end.

We construct iteratively two sequences $(n_j(x))_{j\geq0}$ and $(n'_j(x))_{j\geq0}$ as follows:
\begin{itemize}[parsep=-0.15cm,itemsep=0.25cm,topsep=0.2cm,wide=0.175cm,leftmargin=0.85cm]
\item[(a)] $n_0(x) =0$ and $n_0'(x)$ is so large that
\begin{equation}\label{*etea}
\text{$\eps_2(1-\alpha_0-\ep_1)|\log \ell_{n'_0(x)} (x)| > \log M$.}
\end{equation}
\item[(b)] For $j\geq 1$, let $\widetilde n_j(x)$ be the smallest integer strictly larger than $n'_{j-1}(x)$ such that
\begin{equation}\label{*mntj}
|m_{\widetilde n_j(x)}(x)|<M,
\end{equation}
that is $|m_{\widetilde n_j(x)}(x)|$ is very small after a large value given by
$|m_{n'_{j-1}(x)}(x)|$.
\item[(c)] Then let $n_j'(x)$ be the smallest integer strictly larger than $\widetilde n_j(x)$ satisfying \eqref{*mnx38a} with $n_j'(x)$ instead of $\overline n_j'(x)$.
That is, after $\widetilde n_j(x)$ we look for the first time when $|m_n(x)|$ is large.
Since $n_j'(x)>n'_0(x)$ we obtain $|\log\ell_{n'_j(x)}(x)|> |\log\ell_{n'_0(x)}(x)|$. Therefore,
by \eqref{*etea} and \eqref{*mntj}
\begin{equation}\label{*mnx39a}
\log |m_{\widetilde n_j(x)}(x)|<\eps_2(1-\alpha_0-\eps_1)|\log\ell_{n'_j(x)}(x)|,
\end{equation}
which in our heuristic interpretation means that $|m_{\widetilde n_j(x)}(x)|$
satisfies a property which we call in our heuristic interpretation ``medium smallness''.
\item [(d)]
Finally, choose the largest integer $n$ such that $\widetilde n_j(x)\leq n < n'_j(x) $ and \eqref{*mnx39a} holds with $\widetilde n_j(x)$ replaced by $n$. Denote this integer by $n_j(x)$. This means that $|m_{n_j(x)}(x)|$ is medium small, moreover
$|m_n(x)|$ is larger than medium small, but smaller than large for $n_j(x)<n<n'_j(x)$. This way $|m_n(x)|$ remains controlled between the two integers $n_j(x)<n<n'_j(x)$.
\end{itemize}
Items (b) and (c) are made possible by the remarks preceding the construction. Thus we have:
\begin{itemize}[parsep=-0.15cm,itemsep=0.25cm,topsep=0.2cm,wide=0.175cm,leftmargin=0.5cm]
\item[$\bullet$] For every $j\geq 1$, $n_j(x)<n'_j(x)<n_{j+1}(x)$, and $\lim_{j\to\infty} n_j(x)=\lim_{j\to\infty} n'_j(x) =\infty$.
\item
For $n_j(x)<n<n'_j(x)$, $|\oI_n(x)|>\ell_{n'_j(x)}(x) = |\oI_{n'_j(x)}(x)|$, hence
\begin{equation} \label{*mnx38aa}
\log|m_n(x)| <  (1-\alpha_0-\eps_1) |\log |\oI_n(x)||< (1-\alpha_0-\eps_1) |\log\ell_{n'_j(x)}(x)|
\end{equation}
and by (d) we know that \eqref{*mnx39a} does not hold for $n\in(n_j(x),n'_j(x))$ and hence
\begin{equation}\label{*31A*a}
\log |m_n(x)| \geq \eps_2(1-\alpha_0-\eps_1)|\log \ell_{n'_j(x)}(x)|\geq \eps_2(1-\alpha_0-\eps_1)n'_j(x)\log 3.
\end{equation}
\end{itemize}


For all integers $j,k\geq0$, set
\begin{equation}\label{*mnx40a}
X_{j,k}=\big\{x\in E:\ 6^{-(k+1)}\leq |\oI_{n_j(x)}(x)|<6^{-k} \big\},
\end{equation}
and put $K_j=\min \{k: X_{j,k}\not=\emptyset \}$. It is easy to see from its definition that $K_j\to\infty$ as $j\to\infty$. Using the hierarchical structure of the intervals from $\cI_n$, one can extract from the family $\{\oI_{n_j(x)}(x)\}_{x\in X_{j,k}}$ a finite number say $L_{j,k}$ many non-overlapping intervals $\{\bbJ_{j,k,l}:= \oI_{n_j(x_l)}(x_l)\}_{l\in \{1,\ldots,L_{j,k}\}}$, such that
\begin{equation}\label{**mnx40}
\bigcup_{x\in X_{j,k}}\oI _{n_j(x)}(x)=\bigcup_{l=1}^{L_{j,k}}\bbJ_{j,k,l}.
\end{equation}
By non-overlapping intervals, we mean here and in the rest of the paper that the interiors of the intervals have empty intersection. Observe that by \eqref{*mnx40a}, $L_{j,k}<6\cdot 6^k$ and the intervals $\bbJ_{j,k,l}$ $l=1,...,L_{j,k}$ are approximately of the same length.

Now one can extract in a same way from the family $\{\oI_{n'_j(x)}(x)\}_{x\in X_{j,k}}$ a finite number, say $L'_{j,k}$ many intervals $\{\oI_{n'_j(x'_{l'})}(x'_{l'})\}_{l'\in \{1,\ldots,L'_{j,k}\}}$ such that
\begin{equation}\label{defE'jk}
\bigcup_{x\in X_{j,k}}\oI_{n'_j(x)}(x)
=\bigcup_{l'=1}^{L'_{j,k}} \oI_{n'_j(x'_{l'})}(x'_{l'}).
\end{equation}
We reorganize them by grouping those that belong to the same $\bbJ_{j,k,l}$. Let
\begin{equation*}
\cL_{j,k,l}= \big\{l'\geq 1: \oI_{n'_j(x'_{l'})}(x'_{l'})\subset \bbJ_{j,k,l}\big\},
\end{equation*}
and write $\bbJ'_{j,k,l,l'}:=\oI _{n'_j(x'_{l'})}(x'_{l'})$ for those intervals $\oI _{n'_j(x'_l)}(x'_l)$ that are included in $\bbJ_{j,k,l}$, i.e. the set in \eqref{defE'jk} is rewritten as
\begin{equation}\label{***mnx40}
\bigcup_{l=1}^{L_{j,k}}\bbJ'_{j,k,l} \,\mbox{ with }\ \bbJ'_{j,k,l}= \bigcup_{l'\in\cL_{j,k,l}}\bbJ'_{j,k,l,l'}\subset \bbJ_{j,k,l}.
\end{equation}
Observe that the definition makes sense since each $\bbJ'_{j,k,l,l'}$ is a subset of exactly one $\bbJ_{j,k,l}$, and these intervals are non-overlapping and organized in a hierarchical manner, see Figure \ref{fig:CantorIntervalJ}.

\begin{lemma}\label{lemma-lj}
Let $d_3>0$ be such that $d_2<d_3<d_1$. Then, there are infinitely many integers $j$ with the following property:
there exist integers $k_j$ and $l_j\leq L_{j,k_j}$ such that
\begin{equation}\label{*mnx41b}
\sum_{l'\in\cL_{j,k_j,l_j}}|\bbJ'_{j,k_j,l_j,l'}|^{d_3}>100\cdot 6^{-k_j}.
\end{equation}
\end{lemma}
In the above inequality the unusual constant $100$ stands for a sufficiently large constant of which we do not need an optimal value. In later arguments it is clear that this value is suitable and one can easily keep track of it in the later calculations.

\begin{proof}
Suppose that for $j>j_0$, for every $k$ and $l\leq L_{j,k_j}$,
\begin{equation}\label{*mnx41d}
\sum_{l'\in \cL_{j,k,l}}|\bbJ'_{j,k,l,l'}|^{d_3}\leq 100\cdot 6^{-k}.
\end{equation}
Consider the cover of $E$ given by $\displaystyle\bigcup_{k=K_j}^{\infty}\bigcup_{l=1}^{L_{j,k}}\bigcup_{l'\in  \cL_{j,k,l}}\bbJ'_{j,k,l,l'}$. Then, using \eqref{*mnx40a}, \eqref{**mnx40}, $L_{j,k}\leq 6\cdot 6^k$, $d_1>d_3$ and later \eqref{*mnx41d}, by also taking into consideration that from $J'_{j,k,l,l'}\subset J_{j,k,l}$ it follows $|J'_{j,k,l,l'}| \leq 6^{-k}$, therefore we obtain
\begin{align*}
\sum_{k=K_j}^{\infty}\sum_{l=1}^{L_{j,k}} \sum_{l'\in \cL_{j,k,l}}|\bbJ'_{j,k,l,l'}|^{d_1}
&  = \sum_{k=K_j}^{\infty}\sum_{l=1}^{L_{j,k}}
\sum_{l'\in\cL_{j,k,l}}|\bbJ'_{j,k,l,l'}|^{d_1-d_3}|\bbJ'_{j,k,l,l'}|^{d_3} \\
& \leq \sum_{k=K_j}^{\infty}\sum_{l=1}^{L_{j,k}}
\sum_{l'\in  \cL_{j,k,l}} 6^{-(d_1-d_3)k} |\bbJ'_{j,k,l,l'}|^{d_3} \\
& = \sum_{k=K_j}^{\infty} 6^{-(d_1-d_3)k} \sum_{l=1}^{L_{j,k}}
\sum_{l'\in \cL_{j,k,l}} |\bbJ'_{j,k,l,l'}|^{d_3} \\
& \leq \sum_{k=K_j}^{\infty}6^{-(d_1-d_3)k}L_{j,k}\cdot 100\cdot 6^{-k} \leq 600\sum_{k=K_j}^{\infty} 6^{-(d_1-d_3)k}<\infty.
\end{align*}
This would imply $d_1\geq \dimh(E)>d_1$, a contradiction.
\end{proof}

\begin{definition}\label{def-Jj}
The set $\bbJ'_j $ consisting of finitely many closed intervals is defined for those $j\geq 1$ for which $l_j$ exists in Lemma \ref{lemma-lj}. We put $\obbJ_j= \bbJ_{j,k_j,l_j}$, $\cL_j= \cL_{j,k_j,l_j}$, and $\bbJ'_{j, l'} = \bbJ'_{j,k_j,l_j,l'}$, and the family $\bbJ'_j$ is defined by
\begin{equation}\label{*defJjv}
\bbJ'_j = \bigcup_{l'\in \cL_j}\bbJ'_{j,l'}\subset\obbJ_j.
\end{equation}
\end{definition}

\begin{figure}[hb!]
\centering
\begin{tikzpicture}[scale=2.5]
 \draw[thick] (0,0) -- (4,0);
 \node at (0,-0.05) [below] {\(0\)};
 \node at (4,-0.05) [below] {\(1\)};
 \draw[thick, black] (0,0.05) -- (0,-0.05);
 \draw[thick, black] (4,0.05) -- (4,-0.05);
\foreach \x in {0.2,0.8,1.2,1.6,2.0,2.4,2.8,3.2}
 \draw[thick, blue] (\x,0.05) -- (\x,-0.05)   (\x,0+0.01) -- (\x+0.1,0+0.01)  (\x+0.1,0.05) -- (\x+0.1,-0.05) ;
 \node at (0.28,0.4) [below, blue] {\(\bbJ_{j,k_j, 1}\)};
 \node at (2.1,0.4) [below, blue] {\(\obbJ_j=\bbJ_{j,k_j, l_j}\)};
 \node at (3.35,0.4) [below, blue] {\(\bbJ_{j,k_j, L_{j,k}}\)};
 \draw[->, dotted, blue, line width=1pt] (3.5,-0.5) -- (3.25,-0.1);
 \draw[->, dotted, blue, line width=1pt] (3.5,-0.5) -- (2.85,-0.1);
 \draw[->, dotted, blue, line width=1pt] (3.5,-0.5) -- (2.45,-0.1);
 \node at (3.6,-0.75) [above, blue] {\(\text{intervals $\obbJ_j$}\)};
 \draw[thick] (1,-1.5) -- (3.5,-1.5);
\foreach \x in {1.2, 1.6, 2.0, 2.4, 2.8,3.0,3.2}
 \draw[thick, red] (\x,-1.45) -- (\x,-1.55)  (\x,-1.49) -- (\x+0.1,-1.49)  (\x+0.1,-1.45) -- (\x+0.1,-1.55);
 \node at (1.25,-1.6) [below, red] {\(\bbJ'_{j,k_j, l_j,1}\)};
 \node at (2.5,-1.6) [below, red] {\(\bbJ'_{j, l'}=\bbJ'_{j,k_j, l_j,l'}\)};
 \draw[->, dotted, red, line width=1pt] (0.5,-1)-- (1.2,-1.4);
 \draw[->, dotted, red, line width=1pt] (0.5,-1) -- (2.85,-1.4);
 \draw[->, dotted, red, line width=1pt] (0.5,-1)  -- (2.45,-1.4);
 \draw[dashed] (2,0) -- (1,-1.5)   (2.1,0) -- (3.5,-1.5) ;
 \draw[thick, black] (1,-1.55) -- (1,-1.45);
 \draw[thick, black] (3.5,-1.55) -- (3.5,-1.45);
 \node at (0.5,-.5) [below, red] {\(\text{intervals composing}\)};
 \node at (0.5,-.7) [below, red] {\(\text{the family } \bbJ_j'\)};
\end{tikzpicture}
\captionsetup{width=1\textwidth}
\caption{Representation of the hierarchical structure of the intervals $\bbJ_{j,k_j,l_j}$ and $\bbJ'_{j,k_j,l_j,l'}$.}
\label{fig:CantorIntervalJ}
\end{figure}

\subsubsection{The family of the intervals $\widetilde\bbJ_{l,\rho_j}$, and the construction of $\widetilde\bbF$}

\noindent

By the multifractal formalism and Comment (c) after Proposition  \ref{recalls-mula2}, the set $E_\fla(\alpha_{\la,\cL})$ is of full Lebesgue measure in $[0,1]$. For this exponent $\alpha_{\la,\cL}$, by the Lemma \ref{lemregulinf}, the set $E_\fla(\alpha_{\la,\cL}) \cap \cI$ still has full Lebesgue measure, and the elements $x$ of this subset satisfy $\lim_{n\to\infty} |m_n(x)|=+\infty$. Up to extraction, there exist an integer $N_\cL$ and another, smaller, subset of $E_\fla(\alpha_{\la,\cL}) \cap \cI$ with positive Lebesgue measure of points $x$ such that, for all $n\geq  N_\cL$,
\begin{equation}\label{*mnx35a}
|\widetilde{m}_n(x)|>100,
\end{equation}
where we recall that $\widetilde{m}_n(x)$ is the idealized slope over $x$, see \eqref{itermnideal} and \eqref{seqideal}.

Consider then $E_\cL\subset E_\fla(\alpha_{\la,\cL})\cap\cI$, a compact subset such that $\dimh(E_\cL)=1$ and, for every $x\in E_\cL$ and $n\geq N_\cL$, \eqref{*mnx35a} is true.

Given $\rho\in(0,1)$, for every $x\in E_\cL$, choose $n(x,\rho)\geq1$
such that $\ell_{n(x,\rho)}(x)<\rho\leq \ell_{n(x,\rho)-1}(x)$. We can extract from the family $\{\oI _{n(x,\rho)}(x)\}_{x\in E_\cL}$ a finite number say $\widetilde L_\rho$ many non-overlapping intervals $\{\oI _{n(x_l,\rho)}(x_l):=\widetilde\bbJ_{l,\rho}\}_{l=1,\ldots,\widetilde  L_\rho}$, such that
\begin{equation}\label{def-Jgras}
\bigcup_{x\in E_\cL} {\oI _{n(x,\rho)}(x)} = \bigcup_{l=1}^{\widetilde L_\rho} \widetilde\bbJ_{l,\rho}.
\end{equation}
By construction, for any $l\geq1$ one has
\begin{equation}\label{*mnx35b}
\frac{\rho}6\leq |\widetilde\bbJ_{l,\rho}|<\rho.
\end{equation}

Let us now choose $d_4<1$. We suppose that it is very close to $1$, such that \eqref{*d2d1d5} holds.
Since $\dimh(E_\cL)=1>d_4$, there is $\rho_\cL>0$ such that if $\rho<\rho_\cL$, then $\sum_{l=1}^{\widetilde L_\rho}|\widetilde\bbJ_{l,\rho}|^{d_4}>100$, which by \eqref{*mnx35b} implies
\begin{equation}\label{*mnx36a}
\widetilde L_\rho>\frac{100}{\rho^{d_4}}.
\end{equation}

\begin{definition}\label{*defwtB}
Fix large integers $j$ and $k_j> N_\cL$ such that \eqref{*mnx41b} is satisfied. Set $\rho_j=6^{-k_j}$ and $\widetilde L_j=\widetilde L_{\rho_j}$. The family of intervals $\widetilde\bbJ_{l,\rho_j}$ is used to define
\begin{equation}\label{def-FA}
\widetilde\bbF = \bigcup_{l=1}^{\widetilde L_{ j}}\widetilde\bbJ_{l,\rho_j} \mbox{ with } \widetilde L_j>100\cdot 6^{k_j {d_4}},
\end{equation}
where the lower bound for $\widetilde L_j$ can be obtained from \eqref{*mnx36a}.
\end{definition}

By \eqref{*mnx35b} used with $\rho_j$, one has
\begin{equation}\label{*Jlrhoj}
6^{-(k_j+1)}<|\widetilde\bbJ_{l,\rho_j}|\leq 6^{-k_j}.
\end{equation}
Recalling that the intervals $\widetilde\bbJ_{l,\rho_j}$ are of the form ${\oI_{n(x_l,\rho_j)}(x_l)}$ with $x_l\in E_\cL$, let us introduce the simple notation
\begin{equation}\label{*bbnAdef}
\widetilde\bbn(x_l)=n(x_l,\rho_j),\text{ that is } x_l\in\oI_{\widetilde\bbn(x_l)}(x_l)=\widetilde\bbJ_{l,\rho_j}\subset\widetilde\bbF.
\end{equation}
Later we also need the following observation
\begin{equation}\label{*nlest} 
 6^{-\widetilde\bbn(x_l)}=6^{-n(x_{l},\rho_{j})}\leq |\oI_{n(x_{l},\rho_{j})}(x_{l})|=| \widetilde\bbJ_{l,\rho_j}|\leq 6^{-k_{j}}.
\end{equation}
The set $\widetilde\bbF$ consists of intervals that satisfy \eqref{*mnx35a} and \eqref{*mnx36a} with $\rho=\rho_j$. Intuitively, the increment of the ``idealized'' slope of $\flan$ as $n$ increases on these intervals is uniformly bounded below with a reasonably large constant and by \eqref{def-FA} we have ``many'' of them.

\subsubsection{The construction of $\bbF$ and its dimension}

\noindent

The integer $j$ is now supposed to be fixed and sufficiently large  and satisfies the previous constraints so that the family $\bbJ'_j$ and $\widetilde\bbF$ are well defined. About the largeness of $j$ see later the paragraph after \eqref{*Inx*}.

To build $\bbF$, we start with the intervals $\widetilde\bbJ_{l,\rho_j}$
of $\widetilde\bbF$, where by \eqref{*mnx35a} the intermediary functions $\flan$ have uniformly sufficiently large slopes on them. Then, inside each interval $\widetilde\bbJ_{l,\rho_j}$ of $\widetilde\bbF$ we are going to build subintervals using affine images, $\Psi_l(\bbJ'_j)$ of intervals from families of type $\bbJ'_j$ constructed in subsection \ref{*ssecJj}. On the intervals $\Psi_l(\bbJ'_j)$ we control the growth of the slope of the intermediate steps $\flan$ which lead to $\fla$. Our objective is then to use the subintervals from the families $\Psi_l(\bbJ'_j)$, denoted by $\Psi_l(\bbJ'_{j,l'})$ to define affine maps from $[0,1]$ onto these intervals that will form a finite IFS. The intuition is that since the value of the slopes $|m_n(x)|$ of $\fla$ is controlled at every scale between $1$ and the length of these subintervals, by self-similarity, the growth of $|m_n(x)|$ will be controlled on the attractor of the IFS, ensuring that points belonging to this attractor will satisfy $\lim_{n\to \infty} |m_n(x)|=+\infty$ and thus will belong to $\cI$. Finally, we will prove that the attractor has a too large Hausdorff dimension, contradicting our assumptions \eqref{ineg-dim} and \eqref{*epso}.

The set $\widetilde\bbF$ being fixed, we now use the family of intervals $\bbJ'_j\subset\obbJ_j$ from Definition \ref{def-Jj}. For any $l\in\{1,\ldots,\widetilde L_j\}$,
\begin{equation}\label{*defpsil}
\text{we select an affine bijection $\PPP_l:\obbJ_j\to\widetilde\bbJ_{l,\rho_j}$.}
\end{equation}
Although Figure \ref{fig:fig79} below was intended to help to understand the proof of Lemma \ref{lem-infty}, and some notation is used on it from the proof of this lemma, it may be helpful to look at it during these earlier stages of our definitions.
This choice is meant to ensure that the construction will keep memory only of what happens between $\obbJ_j$ and the intervals $\bbJ'_{j,l'}$, where the variation of the slopes is controlled by \eqref{*mnx38aa} and \eqref{*31A*a}. There are two such bijections, and to determine the monotonicity of $\PPP_l$, we look at its domain and its range. More precisely, let $x_0\in\widetilde \bbJ_{l,\rho_j}\setminus\cE$ and $x'_0\in\overline{\bbJ}_j\setminus\cE$. By definition $\oI_{\widetilde\bbn(x_0)}(x_0)=\widetilde\bbJ_{l,\rho_j}$, and 
by \eqref{**mnx40} and Definition \ref{def-Jj},
 $\oI_{{n}_j(x'_0)}(x'_0)= \obbJ_j$. Now $\PPP_l$ is chosen monotone increasing (resp. monotone decreasing) if $\eps_{\widetilde\bbn(x_0)}(x_0)\eps_{n_j(x'_0)}(x'_0)=1$ (resp. $\eps_{\widetilde\bbn(x_0)}(x_0)\eps_{n_j(x'_0)}(x'_0)=-1$). This choice is natural since it reflects the possible change in monotonicity between the functions $F^{\la}_{\widetilde\bbn(x_0)}$ and $F^{\la}_{n_j(x'_0)}$ on the corresponding intervals. Observe that the family of intervals $(\PPP_l(\bbJ'_{j,l'}))_{l=1,\ldots,\widetilde L_j,l'\in \cL_j}$ consists of non-overlapping intervals.

Consider now the affine contractive bijections $\FFF_{l,l'}:[0,1]\to  \PPP_l(\bbJ'_{j,l'})$, $l=1,\ldots,\widetilde L_j$, $l'\in\cL_j$. As above, the monotonicity of $\FFF_{l,l'}$ must be compatible with our construction. For this, pick $x_0\in\PPP_l(\bbJ'_{j,l'})\setminus\cE$ and choose $n_{l,l'}(x_0)$ in such a way that $\oI_{n_{l,l'}(x_0)}(x_0)= \PPP_l(\bbJ'_{j,l'})$. Then $\FFF_{l,l'}$ is monotone increasing (resp. monotone decreasing) if $\eps_{n_{l,l'}(x_0)}(x_0)=1$ (resp. $\eps_{n_{l,l'}(x_0)}(x_0)=-1$).

\begin{definition}
The first generation, $\widetilde\bbG$, of the compact attractor $\bbF$
is defined by
\begin{equation}\label{def-FB}
\widetilde\bbG= \bigcup_{l=1}^{\widetilde L_j}\PPP_l(\bbJ'_j) = \bigcup_{l=1}^{\widetilde L_j}\bigcup_{l'\in\cL_j} \PPP_l(\bbJ'_{j, l'})\subset\bigcup_{l=1}^{\widetilde L_j}\widetilde \bbJ_{l,\rho_j}.
\end{equation}
\end{definition}

By \eqref{*mnx40a} and \eqref{*mnx35b},
\begin{equation}\label{*bdtart}
6^{-(k_j+1)}\leq |\widetilde \bbJ_{l,\rho_j}|< \rho_j=6^{-k_j} \text{  and } 6^{-(k_j+1)}\leq | \bbJ_j|<6^{-k_j},
\end{equation}
hence
\begin{equation}\label{*slopepsil}
\text{the absolute value of the slope of $\PPP_l$ is between $\frac16$ and $6$.}
\end{equation}
In particular, denoting by $r_{l,l'}=|\PPP_l(\bbJ'_{j,l'})|$ the contraction ratio of $\FFF_{l,l'}$, one has
\begin{equation}\label{*32D*a}
\frac16|\bbJ'_{j,l'}|\leq r_{l,l'}\leq 6 |\bbJ'_{j,l'}|.
\end{equation}
By \eqref{*mnx41b} and \eqref{*32D*a}, for a fixed $l$,
\begin{equation}\label{*32D*b}
\sum_{l'\in\cL_j}|r_{l,l'}|^{d_3}>\frac{100}{6^{d_3}}\cdot 6^{-k_j}>10\cdot 6^{-k_j}.
\end{equation}
By \eqref{*mnx36a} we have
$\widetilde L_j>100\rho_j^{-d_4}= 100\cdot 6^{k_j d_4}$. Therefore, using \eqref{*32D*a}, \eqref{*32D*b} and the fact that $d_2<d_3$, and finally using \eqref{*d2d1d5}, we obtain
\begin{align}
\label{*rllr}
\sum_{l=1}^{\widetilde L_j}\sum_{l'\in\cL_j} |r_{l,l'}|^{d_2}
& = \sum_{l=1}^{\widetilde L_j}\sum_{l'\in \cL_j} |r_{l,l'}|^{d_3}|r_{l,l'}|^{d_2-d_3} \geq\sum_{l=1}^{\widetilde L_j}\sum_{l'\in \cL_j} |r_{l,l'}|^{d_3}6^{-k_j(d_2-d_3)} \\
& \geq \sum_{l=1}^{\widetilde L_j} 10\cdot 6^{-k_j} 6^{-k_j(d_2-d_3)} \ \geq 100\cdot 6^{k_j(d_4+d_3-d_2)}\cdot 10\cdot 6^{-k_j}>1000. \nonumber
\end{align}

We can now give the precise definition of the set $\bbF$.

\begin{definition}\label{def-F}
We call $\bbF$ the attractor of the IFS $(\FFF_{l,l'})_{l\in\{1,\ldots,\widetilde L_j\}, l'\in \cL_j}$.
\end{definition}

For ease of notation, to avoid double indexing, we re-index the functions $\FFF_{l,l'}$ to obtain the bijections $\tFFF_k:[0,1]\to \bbJ_k$, $k=1,\ldots,K$, with $K=\widetilde L_j \cdot \#(\cL_j)$. Note that
\begin{equation}\label{*psiljvjlv}
\text{the intervals $\bbJ_k$ are of the form $\Psi_l(\bbJ'_{j,l'})$ with suitable values of $l$ and $l'$.}
\end{equation}
The contraction ratio of $\tFFF_k$ is denoted by $r_k$.

\begin{lemma}\label{lemma-dimF}
Let $d_5$ be the unique solution to
\begin{equation}\label{*IFSsd}
\sum_{l=1}^{\widetilde L_j}\sum_{l'\in \cL_j} r_{l,l'}^{d_5} =  \sum_{k=1}^K r_k^{d_5} =  1.
\end{equation}
Then $\dimh(\bbF)=d_5$, and for every $\alpha\in(\alpha_{\min},\alpha_0+\eps_0)$ one has $\dimh(\bbF)>d_2> \tau^*_\mula(\alpha+\ga-1)$.
\end{lemma}

\begin{proof}
It is clear that the IFS $(\tFFF_k)_{k=1,\ldots,K}$ satisfies the Open Set Condition, so, for example, by Theorem $9.3$ in \cite{Fal97}, it follows from \eqref{*IFSsd} that $\dimh(\bbF)=d_5$. In particular, recalling \eqref{*epso}, and \eqref{*rllr} $\dimh(\bbF)=d_5>d_2$. By construction one has $d_2>\tau^*_\mula(\alpha+\ga-1)$ for $\alpha\in(\alpha_{\min},\alpha_0+\eps_0$).
\end{proof}

\subsection{Regularity of $\mula$ and $\fla$ at points in $\bbF$}\label{subsectionregF}

\noindent

Let us introduce the symbolic dynamical system associated with the IFS $(\tFFF)_{k=1,\ldots,K}$. Take the symbolic space $\ds\tOOO=\{\omega=(\omega_1\omega_2\cdots) : \omega_i\in\{ 1,\ldots,K \} \}$, with the one-sided shift $\sigma$. The distance of the points in $\widetilde\Omega$ is given by $d_{\widetilde \Omega}(\omega,\omega')=\sum_{i\geq1}|\omega_i-\omega'_i| 2^{-i}$. One can define the projection $\pi$ of $\widetilde\Omega$ onto $\bbF$ by
\begin{equation*}
\pi(\omega)= \bigcap_{p=1}^{\infty} (\tFFF_{\omega_1}\circ\cdots\circ\tFFF_{\omega_p})([0,1]).
\end{equation*}
Since the intervals $\bbJ_k$ are non-overlapping there is a countable set $\cE'\subset \bbF$ such that if $x\in\bbF\setminus\cE'$, then there is a unique $\widetilde\omega(x)=(\widetilde\omega_p(x))_{p\geq 1} \in \widetilde\Omega$ such that
\begin{equation*}
\{x\}= \bigcap_{p=1}^{\infty} (\tFFF_{\widetilde\omega_1(x)}\circ\cdots\circ\tFFF_{\widetilde\omega_p(x)})([0,1]).
\end{equation*}
One easily sees that $\cE'\subset\cE.$ Since we work with non-atomic measures, the countable exceptional set will be of measure zero and has no influence on the dimensions considered hereafter. On $\bbF\setminus\cE'$ one can define the ``inverse dynamical system" to the IFS. Indeed, for $x\in\bbF\setminus\cE'$ there is a unique interval $\bbJ_k$ containing $x$ and we can define $S(x)=\tFFF_k^{-1}(x)$. For this $x$ we have $\widetilde\omega(S(x))=\sigma\widetilde\omega(x)$ and in general $S^p(x)\in\bbJ_{\widetilde\omega_{p+1}(x)}$ for $p\geq0$. Given a finite word $\overline{\omega}_1\cdots\overline{\omega}_j$, we introduce the cylinder sets $[\overline{\omega}_1\cdots\overline{\omega}_j]=\{\widetilde\omega\in\widetilde\Omega: \widetilde\omega_i=\overline{\omega}_i$ for $i=1,\ldots,j\}$. Clearly $\pi([k])=\bbJ_k$.

\medskip

First we prove that $\bbF$ contains essentially regular points, i.e. points whose digit frequencies in \eqref{*F1*a} all converge.

\begin{lemma}\label{lem-as-F}
$\cH^{d_5}$-almost every point of $\bbF$ is regular.
\end{lemma}

\begin{proof}
The idea here is to introduce a Borel probability measure $\mu_\bbF$ living naturally on $\bbF$. We use the dynamical system introduced at the beginning of this subsection \ref{subsectionregF}. For each integer $k$, let $\bbn_k$ be the integer such that $\oI_{\bbn_k}(x)=\bbJ_k$ for $x\in \bbJ_k\setminus\cE'$. Observe that for $i=0,\ldots,3$ the quantity
\begin{equation}\label{*bbbik}
\text{$\beta_{i,k}:=\beta_i(x,\bbn_k)$ does not depend on which $x\in\bbJ_k\setminus\cE'$ is chosen.}
\end{equation}

Let us define $\xi_1(x)=x$ and, for $p\geq 2$,
\begin{equation*}
\xi_p(x)=
\tFFF_{\widetilde\omega_{p-1}(x)}^{-1}\circ \cdots \circ \tFFF_{\widetilde\omega_1(x)}^{-1}(x)=\pi(\sigma^{p-1}(\pi^{-1}(x)))=S^{p-1}(x).
\end{equation*}
By definition, for every $p$, $\xi_p(x)\in \bbJ_{\widetilde\omega_p(x)}$. In addition, for every $x\in\bbJ_k$, set
\begin{equation}\label{*defnk}
\tI_1(x)=\bbJ_k \text{ and } \bbn(x)=\bbn_k,
\end{equation}
that is $\bbJ_k=\tI_1(x)=\oI_{\bbn(x)}(x)= \oI_{\bbn(\xi_1(x))}(\xi_1(x))$. Next, set
\begin{equation}\label{*bbnpdef}
\tI_p(x)= \oI_{\bbn_p(x)}(x)  \ \text{ where } \ \bbn_p(x)=\bbn(\xi_1(x))+\cdots+\bbn(\xi_p(x)).
\end{equation}
Since our IFS satisfies the Open Set Condition, it holds $0<\cH^{d_5}(\bbF)<+\infty$. Then, for every Borel set $A\subset\R$, define
\begin{equation*}
\mu_\bbF(A)=\displaystyle \frac{\cH^{d_5}(A\cap\bbF)}{\cH^{d_5}(\bbF)}.
\end{equation*}
In particular, $\mu_\bbF(\bbJ_k)=\mu_\bbF(\pi[k])=r_k^{d_5}$.

Since the IFS has no overlaps we can define the measure $\nu$ on $\widetilde\Omega$ by $\nu(A)=\mu_\bbF(\pi(A))$ for Borel sets $A\subset\widetilde\Omega$. By construction, it satisfies $\nu([\overline{\omega}_1\cdots\overline{\omega}_j])=r_{\overline{\omega}_1}^{d_5}\cdots r_{\overline{\omega}_j}^{d_5}$. Hence, it is a Bernoulli measure on $\widetilde\Omega$, and so a fortiori an ergodic measure on $(\widetilde\Omega,\sigma)$, likewise $\mu_\bbF$ is ergodic for the dynamical system $(\bbF\setminus\cE',S)$. Applying the ergodic theorem, we obtain that for $\mu_\bbF$-almost every $x\in\bbF$ and every $l=1,\ldots,K$, the limit
\begin{equation*}
\widetilde\beta_l :=\lim_{p\to\infty}\frac1{p}\#\big\{j\leq p: \widetilde\omega_j(x)=l\big\}
\end{equation*}
exists and $\sum_{l=1}^K\widetilde\beta_l =1$.

Recalling \eqref{*defnk}, the limit $\kappa=\lim_{p\to\infty}\frac{\bbn_p(x)}{p}=
\sum_{l=1}^K\widetilde\beta_l\bbn_l$ also exists. This implies that, for $i=0,\ldots,3$ and for $\mu_\bbF$-almost every $x$,
\begin{align}\label{*bbbib}
\frac{\beta_i(x,\bbn_p(x))}{\bbn_p(x)}
= \frac1{\bbn_p(x)}\sum_{s=1}^p\bbb_i(\xi_s(x),\bbn(\xi_s(x)))
& = \frac{p}{\bbn_p(x)}\frac1{p}\sum_{s=1}^p\bbb_i(\xi_s(x),\bbn(\xi_s(x))) \nonumber \\
& \underset{p\to\infty}{\longrightarrow}\frac1{\kappa}\sum_{l=1}^K\widetilde\beta_l\beta_{i,l}=:\beta_i,
\end{align}
where $\beta_{i,l}$ was given at \eqref{*bbbik}.

By \eqref{*defnk} and \eqref{*bbnpdef}, the difference $\bbn_{p+1}(x)-\bbn_p(x)$ can take only finitely many different values and hence it is bounded while $\bbn_p(x)\to\infty$.
This implies that in \eqref{*bbbib} instead of taking limit along the subsequence $\bbn_p(x)$ we can take a similar limit along the natural numbers.

Therefore, we obtain that $\mu_\bbF$-almost every point of $\bbF$ is regular.
\end{proof}

In order to compute the pointwise H\"older exponents of $\fla$ at the typical points in $\bbF$, we first prove that these points $x$ belong to $\cI$.

From the definition of $\mu_\bbF$ it is clear that it is equivalent that a property holds for $\mu_\bbF$-almost every $x\in\bbF$, or this property holds for $\cH^{d_5}$-almost every point of $\bbF$.

\begin{lemma}\label{lem-infty}
For $\mu_\bbF$-almost every $x\in \bbF$,
\begin{equation}\label{*eqmozgatott}
\lim_{n\to\infty}\frac{\log\tm_n(x)}{|\log\ell_n(x)|} =\lim_{p\to\infty}\frac{\log\tm_{\bbn_p(x)}(x)}{|\log \ell_{\bbn_p(x)}(x)|}\geq (1-3\eps_2)(1-\alpha_0-\eps_1)>0.
\end{equation}
In particular, $\lim_{n\to\infty} |m_n(x)|=\lim_{n\to\infty} \tm_n(x)=+\infty$.
\end{lemma}

\begin{figure}[ht!]
\centering
\includegraphics[width=1\textwidth]{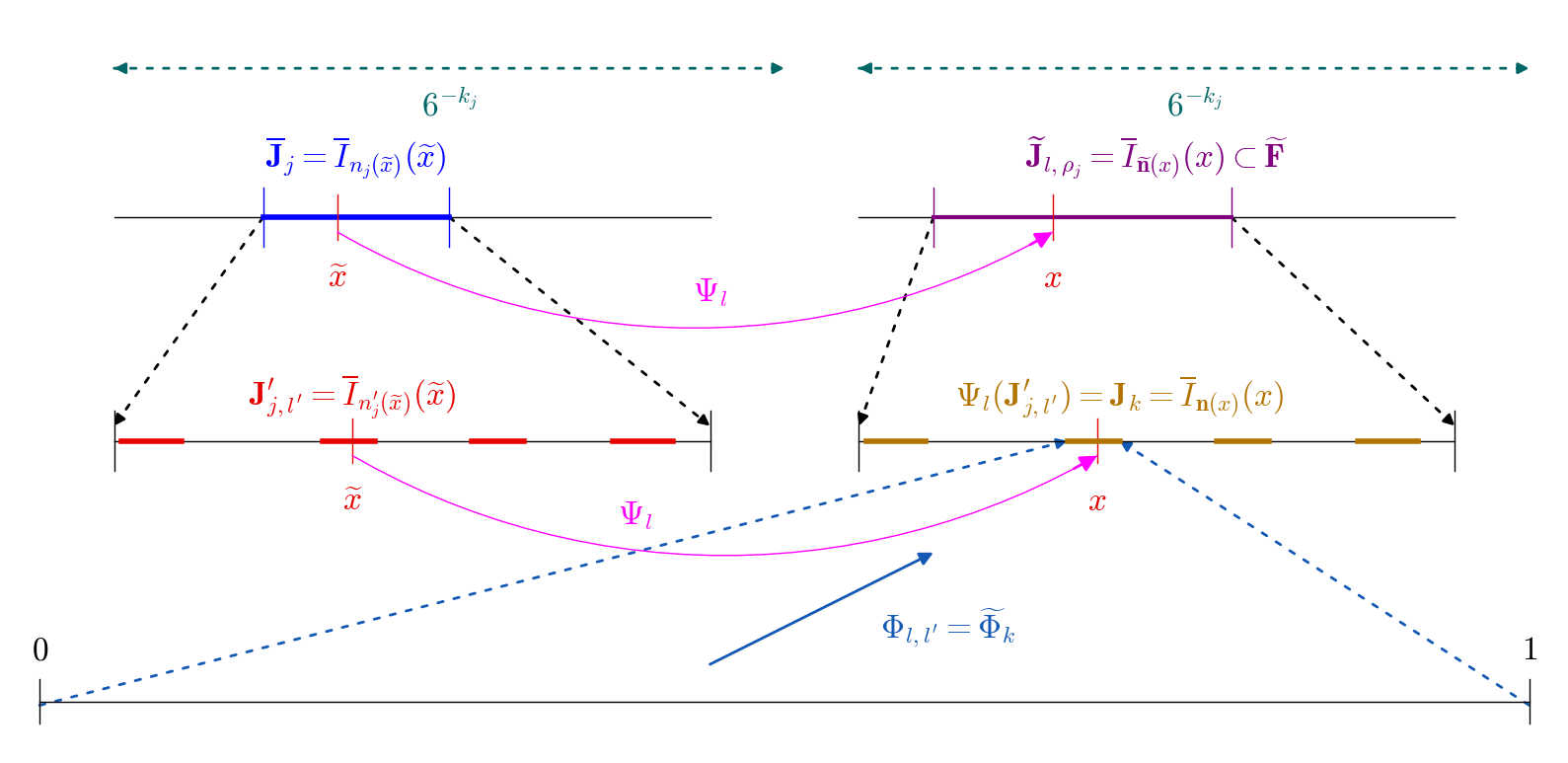}
\captionsetup{width=1\textwidth}
\caption{Representation of the construction of the components of $\widetilde\bbF$ by copying by $\Psi_l$ the construction at point $\tx:=\PPP_l^{-1}(x)\in\obbJ_j=\oI_{n_{j}(\tx)}(\tx)$ to the one at point $x\in\widetilde\bbJ_{l,\rho_j}$.}
\label{fig:fig79}
\end{figure}

\begin{proof}
Observe that \eqref{*eqmozgatott} is stronger than the statement $\lim_{n\to \infty} \tm_n(x)=+\infty
$, and the last statement is immediate since   $|m_n(x)|\geq \tm_n(x)$.

Recall that $\bbJ_{\widetilde\omega_s(x)}=\oI_{\bbn(\xi_s(x))}(\xi_s(x))$, \eqref{*bbnAdef} and the definition \eqref{def-FB} of $\widetilde\bbG$. When $x\in\widetilde\bbG$ one has $x\in\widetilde\bbJ_{l,\rho_j}$ with a suitable $l$ and $\oI_{\widetilde\bbn(x)}(x)=\widetilde\bbJ_{l,\rho_j}$. Set $\tx=\PPP_l^{-1}(x)$, see Figure \ref{fig:fig79}.

The monotonicity of the affine maps $\Psi_l$ defined in \eqref{*defpsil} is compatible with our construction, that is we copied the construction at $\tx\in\obbJ_j=\oI_{n_{j}(\tx)}(\tx)$ by $\Psi_l$ to the one at $x\in\widetilde\bbJ_{l,\rho_j}$, recall also \eqref{*defJjv}. Observe that this implies that the digits $u_{\widetilde\bbn(x)+k}(x)= u_{n_j(\widetilde x)+k}(\widetilde x)$ for $k=0,\ldots, n'_j(\widetilde x)- n_j(\widetilde x)-1$. By \eqref{iterfuncfg} and \eqref{itermnideal},
\begin{equation*}
\frac{\tm_{\widetilde\bbn(x)+k+1}(x)}{\tm_{\widetilde\bbn(x)+k}(x)} = \frac{\tm_{n_j(\widetilde x)+k+1}(\widetilde x)}{\tm_{n_j(\widetilde x)+k}(\widetilde x)}
\end{equation*} 
for $k=0,\ldots,n'_j(\widetilde x)- n_j(\widetilde x)-1$ since they are determined by $u_{\widetilde\bbn(x)+k}(x)= u_{n_j(\widetilde x)+k}(\widetilde x)$, hence $\tm_n(x)$ changes for $n$ between $\widetilde\bbn(x)$ and $\bbn(x)$ the same way as
$\tm_n(\tx)$ changes when $n$ varies between $n_j(\tx)$ and $n_j'(\tx)$, see \eqref{*mnx38aa} and \eqref{*31A*a} applied to $\tx$. By \eqref{*defnk}
\begin{equation}\label{*tmbx}
\tm_{\bbn(x)}(x)=\tm_{\widetilde\bbn(x)}(x)\frac{\tm_{\bbn(x)}(x)}{\tm_{\widetilde\bbn(x)}(x)}
=\tm_{\widetilde\bbn(x)}(x)\frac{\tm_{\widetilde\bbn(x)+n'_j(\tx)-n_j(\tx)}(x)}{\tm_{\widetilde\bbn(x)}(x)} =\tm_{\widetilde\bbn(x)}(x)\frac{\tm_{n'_j(\tx)}(\tx)}{\tm_{n_j(\tx)}(\tx)}.
\end{equation}

In Definition \ref{*defwtB} we have $k_j> N_\cL$, from \eqref{*nlest}  we obtain that $\widetilde\bbn(x)\geq k_j> N_\cL$.

By \eqref{*mnx35a}, $\tm_{\widetilde\bbn(x)}(x)>100$. Thus
\begin{equation}\label{*mtnt}
\log\tm_{\bbn(x)}(x)>\log100 + \log\tm_{n'_j(\tx)}(\tx) - \log \tm_{n_j(\tx)}(\tx).
\end{equation}
Also, by \eqref{*mnx38a} and \eqref{*mnx39a}, one knows that
\begin{equation}\label{*mnx38aaa}
\log |m_{n'_j(\tx)}(\tx)|-\log|m_{n_j(\tx)}(\tx)|\geq
(1-\eps_2)(1-\alpha_0-\eps_1)|\log\ell_{n'_j(\tx)}(\tx)|.
\end{equation}
This provides us with information about the growth
of $\log|m_{n'_j(\tx)}(\tx)|-\log |m_{n_j(\tx)}(\tx)|$, but additional asymptotic information is needed about $\log \tm_{n'_j(\tx)}(\tx)-\log\tm_{n_j(\tx)}(\tx)$. To obtain this, recall \eqref{itermnbis} and \eqref{iterfuncfg}. With a suitable constant $C>0$, for all $|t|>1$ and for any $i\in\{0,\ldots,3\}$,
\begin{equation}\label{*Cgx}
0\leq \log g_i(t) \leq \frac{C}{t^2}.
\end{equation}
Hence, with another suitable constant $C'>0$, using \eqref{*31A*a} and recalling the definition of the digits $u_k(\tx)$ of $\tx$ in \eqref{def-xn},
\begin{align}\label{*mntdiff}
\big|\log |m_{n'_j(\tx)}(\tx)| & -\log |m_{n_j(\tx)}(\tx)|
- \big(\log \tm_{n'_j(\tx)}(\tx)-\log \tm_{n_j(\tx)}(\tx)\big)\big| \nonumber \\
& = \sum_{k=n_j(\tx)}^{n'_j(\tx)-1}\log g_{u_k(\tx)}(|m_k(\tx)|)\leq C'(n'_j(\tx)-n_j(\tx))3^{-2\eps_2(1-\alpha_0-\eps_1)n'_j(\tx)} \\
& \leq C' n'_j(\tx) 3^{-2\eps_2(1-\alpha_0-\eps_1)n'_j(\tx)}\leq 1 \leq
\eps_2(1-\alpha_0-\eps_1)|\log \ell_{n'_j(\tx)}(\tx)|\nonumber
\end{align}
if $n_j(\tx)$ (and then $n_j'(\tx)$ since $n_j'(\tx)\geq n_j(\tx)$) is sufficiently large.

Using \eqref{*tmbx}, after omission of the positive number $\log 100$, by \eqref{*mtnt}, and \eqref{*mntdiff} we obtain
\begin{equation}\label{*esttmn}
\log\tm_{\bbn(x)}(x) > (1-2\eps_2)(1-\alpha_0-\eps_1)|\log \ell_{n'_j(\tx)}(\tx)|.
\end{equation}
Observe that by \eqref{*slopepsil} and \eqref{*psiljvjlv}, for every $x\in\widetilde\bbG$,
\begin{equation*}\label{*Inxx*}
\frac16 < \frac{|\bbJ'_{j,l'}|}{|\bbJ_k|} =\frac{\ell_{n'_j(\tx)}(\tx)}{|\oI_{\bbn(x)}(x)|}<6 \ \text{ and } \ \ell_{n'_j(\tx)}(\tx)\leq 3^{-n'_j(\tx)}.
\end{equation*}
This implies that
\begin{equation}\label{*Inx*}
|\log |\oI_{\bbn(x)}(x)|| \leq |\log \ell_{n'_j(\tx)}(\tx)| + \log6 \leq (1-\eps_2)^{-1} |\log \ell_{n'_j(\tx)}(\tx)|,
\end{equation}
if $\log 6 \leq \eps_2 n'_j(\tx) \log3 \leq \eps_2 |\log \ell_{n'_j(\tx)}(\tx)|$, which can be achieved by choosing a sufficiently large $n'_j(\tx)$. It can be done since $n'_j(\tx)\to\infty$ and we can assume that the fixed $j$ we work with is sufficiently large. We also remark that this choice can be made uniform for all $\tx$s since from its definition $n_{j+1}(\tx)>n'_j(\tx)>n_j(\tx)$, therefore $n_{j+1}(\tx)>n_j(\tx)+2$.

Hence for $x\in\widetilde\bbG$ we have \eqref{*esttmn} and \eqref{*Inx*}. Therefore
\begin{equation}\label{*tmi}
\frac{\log\tm_{\bbn(x)}(x)}{|\log |\oI_{\bbn(x)}(x)||}\geq
(1-\eps_2)(1-2\eps_2)(1-\alpha_0-\eps_1)> (1-3\eps_2)(1-\alpha_0-\eps_1).
\end{equation}
In addition, recalling that $\xi_s(x)\in \widetilde\bbG$ for all $s$, and using notation from the proof of Lemma \ref{lem-as-F}, we have
\begin{equation}\label{*32F*a}
\ell_{\bbn_p(x)}(x)=\prod_{s=1}^p |\bbJ_{\widetilde\omega_s(x)}| \ \text{ and } \ \tm_{\bbn_p(x)}(x)=\prod_{s=1}^p \tm_{\bbn(\xi_s(x))}(\xi_s(x)).
\end{equation}
Therefore, by using \eqref{*tmi} we obtain
\begin{align*}
\log \tm_{\bbn_p(x)}(x)  = \sum_{s=1}^p\log \tm_{\bbn(\xi_s(x))}(\xi_s(x))
& \geq  (1-3\eps_2)(1-\alpha_0-\eps_1)\sum_{s=1}^p|\log |I_{\bbn(\xi_s(x))}(\xi_s(x))|| \\
& = (1-3\eps_2)(1-\alpha_0-\eps_1)\sum_{s=1}^p|\log |\bbJ_{\widetilde\omega_s(x)}||.
\end{align*}
Hence, combining \eqref{*32F*a} with the fact that $\bbn_{p+1}(x)-\bbn_p(x)$ is uniformly bounded in $p$ and $x$, one deduces that, for all the regular points of $\bbF$ and hence for $\mu_\bbF$-almost every $x\in\bbF$, we have \eqref{*eqmozgatott}. The limits exist in \eqref{*eqmozgatott} by the regularity of $x$, see also \eqref{seqideal}.
\end{proof}

We now show that most points of $\bbF$ share the same pointwise H\"older exponent for $\fla$.

\begin{lemma}\label{lem-estimF}
For $\cH^{d_5}$-almost all points $x\in\bbF$, $\ldimloc(\mula,x)$
takes the same value and it satisfies $\ldimloc(\mula,x)+1-\ga \leq \alpha_0+\eps_0$. The common value of $\ldimloc(\mula,x)+1-\ga$ is denoted by $\alpha_\bbF$. It also holds that $\alpha_\bbF=h_\fla(x)$ for these $\cH^{d_5}$-almost all points $x\in\bbF$.
\end{lemma}

\begin{proof}
Using Lemma \ref{lem-as-F}, for $\cH^{d_5}$-almost every $x\in\bbF$, we derive from \eqref{value-dimloc} that
\begin{align}\label{*mnx314aa}
\ldimloc(\mula,x) = \ga-\frac{\beta_1\log(6\la+1)+\beta_2\log(6\la-1)}{(\beta_0+\beta_3)\log 3+(\beta_1+\beta_2)\log 6}
\end{align}
where $\beta_i$ is given by \eqref{*bbbib}. By Lemma \ref{lem-infty}, $\mu_\bbF$-almost every $x$ belongs to $\cI$, so by \eqref{hexpb}
\begin{equation}\label{*mnx314ab}
\lim_{n\to\infty}\frac{\log|m_n(x)|}{|\log\ell_n(x)|}=
\frac{\beta_1\log(6\la+1)+\beta_2\log(6\la-1)}{(\beta_0+\beta_3)\log 3+(\beta_1+\beta_2)\log 6}=\lim_{n\to\infty}\frac{\log\tm_n(x)}{|\log\ell_n(x)|}.
\end{equation}
Therefore, using Corollary \ref{coro-compare}, \eqref{*mnx314aa} and \eqref{*mnx314ab}, we obtain that for $\mu_\bbF$-almost every $x$,
\begin{equation*}
h_\fla(x)=\alpha_\bbF=\ldimloc(\mula,x)+1-\ga = 1- \lim_{n\to\infty}\frac{\log\tm_n(x)}{|\log\ell_n(x)|}.
\end{equation*}
Finally, it follows from \eqref{*eqmozgatott} that
\begin{equation}\label{*epschoose}
h_\fla(x)=\alpha_\bbF\leq (\ga-(1-3\eps_2)(1-\alpha_0-\eps_1))+1-\ga\leq \alpha_0+\eps_0
\end{equation}
if $\eps_1>0$ and $\eps_2>0$ are chosen sufficiently small from the beginning. Hence, the result.
\end{proof}

We are now able to conclude on Proposition \ref{*sppropb}.

\begin{proof}[Proof of Proposition \ref{*sppropb}.]
By Lemma \ref{lem-infty} $\mu_\bbF$-almost every point of $\bbF$ belongs to $\cI$. Hence by \eqref{coro-compare}, $\alpha_\bbF=h_\fla(x)=1-\ga+\ldimloc(\mula,x)$ for $\mu_\bbF$-almost every $x$ in $\bbF$.

Since  by \eqref{*epschoose} $\alpha_\bbF<\alpha_0+\eps_0$, using \eqref{*epso} and \eqref{form-mula} we obtain that
\begin{equation}\label{*contr}
d_5=\dimh(\bbF)\leq \tau^*_\mula(\alpha_\bbF+\ga-1)<d_2.
\end{equation}
This contradicts the result in Lemma \ref{lemma-dimF} and concludes the proof of Proposition \ref{*sppropb}.
\end{proof}

\bigskip

\noindent{\sc Acknowledgement.} The authors would like to thank the anonymous reviewer for his/her careful reading of the manuscript. We appreciate all of his/her valuable comments and suggestions, which have greatly improved the quality of the manuscript.

\addtocontents{toc}{\vspace{0.2cm}}%
\addcontentsline{toc}{section}{\protect\numberline{}\hspace{-18pt}\bf References}

\let\oldaddcontentsline\addcontentsline
\renewcommand{\addcontentsline}[3]{}

\nocite{*}
\bibliographystyle{alpha}
\bibliography{HVKFunctionPart1Biblio}

\let\addcontentsline\oldaddcontentsline
\end{document}